\documentclass[moor,accepted]{informs3} 

\usepackage{natbib}
 \NatBibNumeric
 \def\bibfont{\small}%
 \bibpunct[, ]{[}{]}{,}{n}{}{,}%

\usepackage[colorlinks=true,breaklinks=true,bookmarks=true,urlcolor=blue,
     citecolor=blue,linkcolor=blue,bookmarksopen=false,draft=false]{hyperref}
\usepackage{multirow}
\def\EMAIL#1{\href{mailto:#1}{#1}}

\usepackage[mathscr]{eucal}
\usepackage{epsfig,epsf,psfrag}
\usepackage{amssymb,amsfonts,latexsym}
\usepackage{amsmath}
\usepackage{graphicx}
\usepackage{epstopdf}
\usepackage{bm,xcolor,url}
\usepackage{fixltx2e}
\usepackage{array}
\usepackage{verbatim}
\usepackage{algpseudocode, cite}
\usepackage{algorithm}
\usepackage{verbatim}
\usepackage{textcomp}
\usepackage{mathrsfs}
\usepackage[referable]{threeparttablex}
\usepackage{subfig}
\usepackage{enumerate}
\usepackage{footnote}
\usepackage{enumitem}
\usepackage{xr}
\usepackage{mathtools}
\catcode`~=11 \def\UrlSpecials{\do\~{\kern -.15em\lower .7ex\hbox{~}\kern .04em}} \catcode`~=13 

\allowdisplaybreaks[3]

\newcommand{\tnorm}[1]{{\left\vert\kern-0.25ex\left\vert\kern-0.25ex\left\vert #1 
    \right\vert\kern-0.25ex\right\vert\kern-0.25ex\right\vert}}
\newcommand{\tnormt}[1]{{\vert\kern-0.25ex\vert\kern-0.25ex\vert #1 
    \vert\kern-0.25ex\vert\kern-0.25ex\vert}}

\newcommand{\where}{\mbox{where}}
\newcommand{\ipt}{\lranglet}

\newcommand{\clconv}{\mathsf{clconv}\,}

\newcommand{\norm}[1]{\left\Vert#1\right\Vert}
\newcommand{\normt}[1]{\Vert#1\Vert}
\newcommand{\abst}[1]{\vert#1\vert}
\newcommand{\abs}[1]{\left\lvert#1\right\rvert}

\newcommand{\nn}{\nonumber}

\newcommand{\defeq}{:=}
\newcommand{\eqcst}{\stackrel{{\rm c}}{=}}
 
\newcommand{\nt}{\addtocounter{equation}{1}\tag{\theequation}} 
\newcommand{\prox}{\mathsf{prox}}
\newcommand{\dom}{\mathsf{dom}\,}

\newcommand{\inter}{\mathsf{int}\,}
\newcommand{\bdry}{\mathsf{bd}\,}

\newcommand{\dist}{\mathsf{dist}\,}

\newcommand{\barDelta}{\bar{\Delta}}
\newcommand{\tilnabla}{\widetilde{\nabla}}

\newcommand{\bareps}{\bar{\epsilon}}

\newcommand{\calA}{\mathcal{A}}
\newcommand{\calB}{\mathcal{B}}

\newcommand{\calP}{\mathcal{P}}

\newcommand{\calU}{\mathcal{U}}
\newcommand{\calV}{\mathcal{V}}
\newcommand{\calW}{\mathcal{W}}
\newcommand{\calX}{\mathcal{X}}
\newcommand{\calY}{\mathcal{Y}}

\newcommand{\barcalX}{\bar{\calX}}

\newcommand{\rmd}{\mathrm{d}}
\newcommand{\rmD}{\mathrm{D}}
\newcommand{\rme}{\mathrm{e}}

\newcommand{\rmp}{\mathrm{p}}
\newcommand{\rmP}{\mathrm{P}}

\newcommand{\bbE}{\mathbb{E}}

\newcommand{\bbN}{\mathbb{N}}

\newcommand{\bbR}{\mathbb{R}}

\newcommand{\bbU}{\mathbb{U}}

\newcommand{\bbX}{\mathbb{X}}
\newcommand{\bbY}{\mathbb{Y}}

\newcommand{\barbbR}{\overline{\bbR}}

\DeclareMathAlphabet{\mathbsf}{OT1}{cmss}{bx}{n}

\newcommand{\rvA}{\mathsf{A}}

\newcommand{\hath}{\hat{h}}

\newcommand{\hatL}{\widehat{L}}

\newcommand{\hatv}{\widehat{v}}

\newcommand{\hatx}{\hat{x}}

\newcommand{\tily}{\widetilde{y}}

\newcommand{\barf}{\bar{f}}

\newcommand{\barp}{\overline{p}}

\newcommand{\barr}{\overline{r}}

\newcommand{\bart}{\overline{t}}
\newcommand{\baru}{\overline{u}}
\newcommand{\barv}{\overline{v}}

\newcommand{\barx}{\overline{x}}
\newcommand{\bary}{\overline{y}}

\newcommand{\barC}{\overline{C}}

\newcommand{\barG}{\overline{G}}

\newcommand{\barK}{\bar{K}}
\newcommand{\barL}{\bar{L}}

\newcommand{\barxi}{\overline{\xi}}
\newcommand{\barPsi}{\bar{\Psi}}

\newcommand{\lrangle}[2]{\left\langle{#1},{#2}\right\rangle}
\newcommand{\lranglet}[2]{\langle{#1},{#2}\rangle}

\newcommand{\eqa}{\stackrel{\rm(a)}{=}}

\newcommand{\eqc}{\stackrel{\rm(c)}{=}}

\newcommand{\lea}{\stackrel{\rm(a)}{\le}}
\newcommand{\leb}{\stackrel{\rm(b)}{\le}}
\newcommand{\lec}{\stackrel{\rm(c)}{\le}}

\newcommand{\gea}{\stackrel{\rm(a)}{\ge}}
\newcommand{\geb}{\stackrel{\rm(b)}{\ge}}
\newcommand{\gec}{\stackrel{\rm(c)}{\ge}}
\newcommand{\ged}{\stackrel{\rm(d)}{\ge}}

\newcommand{\qednew}{\nobreak \ifvmode \relax \else
      \ifdim\lastskip<1.5em \hskip-\lastskip
      \hskip1.5em plus0em minus0.5em \fi \nobreak
      \vrule height0.75em width0.5em depth0.25em\fi}

\TheoremsNumberedBySection  

\EquationsNumberedBySection 

\begin{document}


\RUNAUTHOR{Zhao}

\RUNTITLE{A Primal-Dual Smoothing Framework for Max-Structured Non-Convex Optimization}

\TITLE{A Primal-Dual Smoothing Framework for Max-Structured Non-Convex Optimization}

\ARTICLEAUTHORS{
\AUTHOR{Renbo Zhao}
\AFF{Department of Business Analytics, Tippie College of Business, 
University of Iowa \EMAIL{renbo-zhao@uiowa.edu}}
} 

\ABSTRACT{
We propose a primal-dual smoothing framework for finding a near-stationary point of a class of non-smooth non-convex optimization problems  with max-structure. We analyze the primal and dual gradient complexities of the framework via two approaches, i.e., the dual-then-primal and primal-the-dual smoothing approaches. Our framework improves the best-known oracle complexities of the existing method, even in the restricted problem setting. As an important part of our framework, we propose a first-order method for solving a class of (strongly) convex-concave saddle-point problems,  which is based on a newly developed non-Hilbertian inexact accelerated proximal gradient algorithm for strongly convex composite minimization that enjoys duality-gap convergence guarantees. Some variants and extensions of our framework are also discussed.
}

\MSCCLASS{Primary: 90C47; secondary: 90C26}
\ORMSCLASS{Programming: Nonlinear: Algorithms}

\maketitle

 
\vspace{-1cm}

\section{Introduction.}
We consider a class of non-convex non-smooth optimization problems, where the non-convex function has a max-structure. Let us first formally state the problem. 

\subsection{Problem statement.}\label{sec:Problem}
Let $(\bbX,\normt{\cdot}_\bbX)$ and $(\bbY,\normt{\cdot}_\bbY)$ be finite-dimensional real normed spaces, with dual spaces denoted by $(\bbX^*,\normt{\cdot}_{\bbX^*})$ and $(\bbY^*,\normt{\cdot}_{\bbY^*})$, respectively. 
Let us consider the following optimization problem: 
\begin{equation}
q^*\defeq {\min}_{x\in\bbX}\; \big\{q(x)\defeq f(x) + r(x)\big\}, \quad\mbox{where}\quad f(x)\defeq\; {\max}_{y\in\calY} \;\Phi(x,y) - g(y).\label{eq:main}
\end{equation}
In~\eqref{eq:main}, the function $r:\bbX\to\bbR\cup\{+\infty\}$ is a closed and convex function 
with domain $\dom r:= \{x\in\bbX:r(x)<+\infty\}$. 
For convenience, define $\calX:= \dom r$, and we assume that both $\calX$ and $\calY\subseteq\bbY$ are nonempty, closed and convex sets, and additionally, $\calY$ is bounded. 
The function $g:\calY\to\bbR$ is convex and continuous on $\calY$. We do not require $r$ or $g$ to be 
differentiable, but instead assume that both $r$ and $g$ are ``simple'' in the sense that certain associated Bregman proximal projection (BPP) problems (formally introduced in Section~\ref{sec:BPP}) are easily solvable. 
Let $\calX'\subseteq\bbX$ and $\calY'\subseteq\bbY $ be some open sets that contain $\calX$ and $\calY$, respectively. We let the function $\Phi:\calX'\times\calY'\to\bbR$ be jointly continuous on $\calX'\times\calY'$ and $\Phi(x,\cdot)$ be concave on $\calY$, for any $x\in\calX$. In addition, we let $\Phi$ satisfy the following assumptions. 
\begin{assumption}[\rm Smoothness of $\Phi(\cdot,y)$]\label{assump:smooth_Phi}
 For any $y\in\calY$, $\Phi(\cdot,y)$ is (Fr\'echet) differentiable on $\calX'$, with the 
 gradient at $x\in\calX'$ denoted by $\nabla_x \Phi(x,y)$. 
 Furthermore, there exist Lipschitz parameters $L_{xx},\hatL_{xy}<+\infty$ such that for any $x,x'\in\calX$ and any $y,y'\in\calY$, 
\begin{align}
&\normt{\nabla_x \Phi(x,y) - \nabla_x \Phi(x',y)}_{\bbX^*}\le L_{xx}\normt{x-x'}_\bbX,\label{eq:L_xx}\\
&\normt{\nabla_x \Phi(x,y) - \nabla_x \Phi(x,y')}_{\bbX^*}\le \hatL_{xy}\normt{y-y'}_\bbY.\label{eq:L_xy}
\end{align} 
\end{assumption}

\begin{assumption}[\rm Weak Convexity of $\Phi(\cdot,y)$]\label{assump:weak_convex_x}
For any $y\in\calY$, $\Phi(\cdot,y)$ is $\gamma$-weakly convex on $\calX$ for some $\gamma\in(0,L_{xx}]$, 
i.e., for any $x,x'\in\calX$,
\begin{equation}
\Phi(x',y) - \Phi(x,y) - \lranglet{\nabla_x \Phi(x,y)}{x'-x}\ge -(\gamma/2)\norm{x'-x}_\bbX^2. \label{eq:weak_cvx}
\end{equation} 
\end{assumption}

\begin{assumption}[\rm Smoothness of $\Phi(x,\cdot)$]\label{assump:smooth_Phi2}
 For any $x\in\calX$, $\Phi(x,\cdot)$ is (Fr\'echet) differentiable on some open set $\calY'\supseteq\calY$, with the gradient at $y\in\calY$ denoted by $\nabla_y \Phi(x,y)$. 
 Furthermore, there exist Lipschitz parameters  $\hatL_{yx},L_{yy}<+\infty$ such that for any $x,x'\in\calX$ and $y,y'\in\calY$, 
\begin{align}
&\normt{\nabla_y \Phi(x,y) - \nabla_y \Phi(x',y)}_{\bbY^*}\le \hatL_{yx}\normt{x-x'}_\bbX,\label{eq:L_yx}\\
&\normt{\nabla_y \Phi(x,y) - \nabla_y \Phi(x,y')}_{\bbY^*}\le L_{yy}\normt{y-y'}_\bbY.\label{eq:L_yy}
\end{align} 
\end{assumption}


Before introducing some applications of the problem in~\eqref{eq:main}, we make some remarks. First, note that in Assumption~\ref{assump:weak_convex_x}, the reason we have $\gamma\le L_{xx}$ is due to~\eqref{eq:L_xx} in Assumption~\ref{assump:smooth_Phi}. Specifically, from~\eqref{eq:L_xx}, by the descent lemma (see e.g.,~\citet[Lemma~1.30]{Peyp_15}), we see that
\begin{equation}
\abst{\Phi(x',y) - \Phi(x,y) - \lranglet{\nabla_x \Phi(x,y)}{x'-x}} \le (L_{xx}/2)\normt{x'-x}_\bbX^2,\quad \forall\,x,x'\in\calX, \;\;\forall\, y \in\calY,  
\end{equation}
which implies~\eqref{eq:weak_cvx} with $\gamma=L_{xx}$. 
However, as we will see later, our proposed smoothing framework can take advantage of the situation where $\gamma\ll L_{xx}$. 
Second, we unify the ``cross'' Lipschitz parameters $\hatL_{xy}$ and $\hatL_{yx}$ into a single one, by defining a new parameter $$L_{xy}:= \max\{\hatL_{xy},\hatL_{yx}\},$$ so that both~\eqref{eq:L_xy} and~\eqref{eq:L_yx} hold with $L_{xy}$. (Indeed, under certain regularity conditions of $\Phi$, the tightest choices of $\hatL_{xy}$ and $\hatL_{yx}$ coincide, and we can set the value to be $L_{xy}$.)
Finally, for well-posedness, we will always assume that the optimal value of~\eqref{eq:main} is finite, namely $q^*>-\infty$. 
 
%
%
%

\subsection{Applications.}\label{sec:applications} 

The problem in~\eqref{eq:main} has many applications, from which we detail three.  

\begin{example}[Distributionally robust learning]\label{eg:DRO}
In learning theory, an important problem is population risk minimization (PRM), which reads
\begin{equation}
{\min}_{x\in\calX}\; \bbE_{\xi\sim p}[\ell(x,\xi)].
\label{eq:PRM}
\end{equation} 
In~\eqref{eq:PRM}, the optimization variable $x$ represents (the coefficients of) the model that one intends to learn, and typically constrained in some closed and convex set $\calX\ne \emptyset$. Since we have some uncertainty about the problem data, we model it as a random variable $\xi$ with distribution $p$ and support $\Xi$. Given a model $x$ under a realization of the problem data $\xi$, the loss function $\ell:\bbX\times\Xi\to\bbR$ returns loss $\ell(x,\xi)$, which then gives the {\em population risk} $\bbE_{\xi\sim p}[\ell(x,\xi)]$ after taking expectation. The optimal solution of~\eqref{eq:PRM} then represents the learned model. To improve the statistical properties of the learned model (e.g., unbiasedness or sparsity), one typically adds a regularizer $\barr:\calX\to\bbR$ to the objective function in~\eqref{eq:PRM}, and solves the regularized PRM problem instead: 
\begin{equation}
{\min}_{x\in\calX}\; \bbE_{\xi\sim p}[\ell(x,\xi)] + \barr(x).\label{eq:PRM_reg}
\end{equation}

For simplicity, we will consider the case where $\xi$ can only take finitely many values, 
and denote the support $\Xi\defeq \{\xi_1,\ldots,\xi_n\}$. The formulation in~\eqref{eq:PRM} or~\eqref{eq:PRM_reg} implicitly assumes that the distribution $p$ is known exactly, which is not the case 
 in many circumstances.  
 However, for most of these situations, 
we know an uncertainty set $\calP$ that $p$ belongs to (by using either prior knowledge or certain estimation procedures). Some typical examples of $\calP$ include 
\begin{align*}
\calB_{\rm TV}(\barp,\alpha)\defeq \{p\in\Delta_n:d_{\rm TV}(p,\barp)\le \alpha\} \quad \mbox{or}\quad \calB_{\rm W_2}(\barp,\alpha)\defeq \{p\in\Delta_n:d_{\rm W_2}(p,\barp)\le \alpha\},
\end{align*}
where $\Delta_n\defeq \{p\in\bbR^n:p\ge 0,\sum_{i=1}^np_i = 1\}$ denotes the (standard) unit simplex in $\bbR^n$, $\barp\in\Delta_n$ denotes the nominal distribution, $\alpha>0$ denotes the ``radius'' of $\calP$,    and $d_{\rm TV}$ and $d_{\rm W_2}$ denote the total variation and 2-Wasserstein distances, respectively. Based on $\calP$, we solve the following distributionally robust regularized PRM instead:
\begin{equation}
{\min}_{x\in\calX}\;{\max}_{p\in\calP} \; \textstyle\sum_{i=1}^n p_i\ell(x,\xi_i)+\barr(x). \label{eq:DRO}
\end{equation}

We note that~\eqref{eq:DRO} fits into the template in~\eqref{eq:main} if we define $\Phi:(x,p)\mapsto\textstyle\sum_{i=1}^n p_i\ell(x,\xi_i)$, for $x\in\calX$ and $p\in\calP$, and $r:= \barr + \iota_\calX$, where $\iota_\calX$ denotes the indicator function of $\calX$, namely
\begin{equation}
\iota_\calX(x)=\begin{cases}
0,&\quad x\in\calX\\[1ex]
+\infty, & \quad x\not\in\calX
\end{cases}. 
\end{equation}
 In the learning problems, we usually assume that $\ell(\cdot,\xi)$ is  $L_\xi$-smooth on $\calX$, namely $\ell(\cdot,\xi)$ is differentiable on some open set $\calX'\supseteq \calX$ with $L_\xi$-Lipschitz gradient on $\calX$. As such, $\Phi(\cdot,p)$ satisfies~\eqref{eq:L_xx} with $L_{xx}\defeq\sup_{p\in\calP}\sum_{i=1}^n p_i L_{\xi_i}<+\infty$. In addition, we see that 
\begin{equation}
\normt{\nabla_x \Phi(x,p) - \nabla_x \Phi(x,p')}_*\le L_{xp}\normt{p-p'}_1,\quad\mbox{where}\quad
L_{xp}\defeq {\sup}_{x\in\calX}\; {\max}_{i\in[n]}\; \normt{\nabla_x \ell({x},\xi_i)}_*.\nn
\end{equation}
Clearly, $L_{xp}<+\infty$ if $\calX$ is bounded. If $\calX$ is unbounded, 
we may still 
have $L_{xp}<+\infty$ --- for example, this happens if $\ell(x,\xi)=l\left(\lrangle{\xi}{x}\right)$, where $l:\bbR\to\bbR$ has bounded derivative on $\bbR$. Examples of such $l$ include the Huber loss function or the quadratically smoothed hinge loss function~\citep{Zhang_04}. The above shows that~\eqref{eq:DRO} satisfies Assumption~\ref{assump:smooth_Phi}. Using similar reasoning, we see that Assumptions~\ref{assump:weak_convex_x} and~\ref{assump:smooth_Phi2} are satisfied as well. 
\end{example}

\begin{example}[Minimizing maximum of smooth functions] \label{eg:max_smooth}
Given 
$n$ functions $\{f_i:\calX\to\bbR\}_{i=1}^n$ such that each $f_i$ is $L_i$-smooth on $\calX$, we aim to minimize their point-wise maximum, namely
\begin{equation}
{\min}_{x\in\calX}\,\big[f(x):= {\max}_{i\in[n]} \, f_i(x) = {\max}_{p\in\Delta_n} \textstyle\sum_{i=1}^n p_i\ell_i(x)\big], \quad \forall\,x\in\calX.  \label{eq:min_smooth_func}
\end{equation}
If we let $\Phi(x,p)\defeq \sum_{i=1}^n p_i\ell_i(x)$ for $x\in\calX$ and $p\in\Delta_n$ and $r:= \iota_\calX$, then~\eqref{eq:min_smooth_func} fits into the template in~\eqref{eq:main}. 
Therefore, following the discussions  in Example~\ref{eg:DRO}, we see that~\eqref{eq:min_smooth_func} also satisfies Assumptions~\ref{assump:smooth_Phi} to~\ref{assump:smooth_Phi2}.
\end{example}

\begin{example}[Dual problem of composite optimization]\label{eg:Dual_CO}
Consider the  following composite optimization problem, where one wishes to solve 
\begin{align}
{\min}_{x\in\calX}\; h(c(x)) + \barr(x),\label{eq:CO}
\end{align}
where $c:\bbX\to\bbR^n$ is (Fr\'echet) differentiable on $\calX$ with Jocobian at $x\in\calX$ denoted by $J_c(x)$,  and $h:\calX\to\bbR$ is a closed and convex function and is Lipschitz on $\calX$. We assume that $J_c$ is Lipschitz on $\calX$, namely there exists $b<+\infty$ such that $\normt{J_c(x) - J_c(x')}\le b\normt{x-x'}$, for any $x,x'\in\calX$. In addition, $\barr:\calX\to\bbR$ is some regularizer.  
As detailed in~\citet{Davis_19}, the problem in~\eqref{eq:CO} has many applications, e.g., robust phase retrieval~\citep{Candes_13} and covariance matrix estimation~\citep{Chen_15}. Using Fenchel duality, we can rewrite~\eqref{eq:CO} as 
\begin{equation}
{\min}_{x\in\calX}\;{\max}_{y\in\dom h^*}\;\textstyle\sum_{i=1}^n y_i{c_i(x)} - h^*(y) + r(x), \label{eq:dual_CO}
\end{equation}
where for each $i\in[n]$, $c_i:\bbX\to\bbR$ denotes the $i$-th component of the (vector-valued) function $c$. 
Additionally, by the Lipschitz continuity of $h$, we see that $\dom h^*\subseteq\bbR^n$ is nonempty, convex and bounded, and is also closed in many cases of interest, e.g., $h$ is any norm function. 
Note that by defining $\Phi(x,y) := \sum_{i=1}^n y_i{c_i(x)}$, $\barr:= r+\iota_\calX$ and $g:=h^*$,~\eqref{eq:dual_CO} fits into the template in~\eqref{eq:main}, and since $\Phi(\cdot,\cdot)$ takes a similar form to that in Example~\ref{eg:DRO}, \eqref{eq:dual_CO} also satisfies Assumptions~\ref{assump:smooth_Phi} to~\ref{assump:smooth_Phi2}. 
\end{example}

\subsection{Convergence criterion: near-stationary point.}\label{sec:NS-point}
Since the objective function in~\eqref{eq:main} is non-convex, without any additional assumptions on the problem structure, it is generally NP-hard to obtain an approximate optimal solution of~\eqref{eq:main} for any desired accuracy.  Therefore, recent research (see e.g.,~\citet{Davis_18,Davis_19,Davis_19b}) has been focusing on finding an $\varepsilon$-near-stationary ($\epsilon$-NS) point of~\eqref{eq:main}, which we introduce informally below under the simpler setting that 
$\bbX$ is Hilbertian, namely the norm $\normt{\cdot}=\normt{\cdot}_\bbX$ can be induced by some inner product. The formal definition of $\varepsilon$-NS point for the general normed space involves the notion of Bregman divergence and is deferred to Section~\ref{sec:eps_near_stat}. 
 Let us first define the proximal point at $x$ with function $q$ and step-size $0<\lambda<\gamma^{-1}$ as
\begin{equation}
\prox(q,x,\lambda)\defeq {\argmin}_{x'\in\bbX}\,\; q(x') + (2\lambda)^{-1}\normt{x'-x}^2. \label{eq:prox_Euc}
\end{equation}
(Note that as we will see in Section~\ref{sec:pd_smoothing_framework}, the condition $0<\lambda<\gamma^{-1}$ ensures that the minimization problem in~\eqref{eq:prox_Euc} is strongly convex and hence admits a unique optimal solution.) 
We call $x\in\calX$ an $\varepsilon$-NS point of~\eqref{eq:main}, if there exists some $0<\lambda<\gamma^{-1}$ such that 
\begin{align}
\normt{\lambda^{-1}(x-\prox(q,x,\lambda))}\le \varepsilon. \label{eq:eps_near_stat_pt}
\end{align}
As we will see later (cf.\ Lemma~\ref{lem:eps_subdiff}), if~\eqref{eq:eps_near_stat_pt} holds, then $\prox(q,x,\lambda)$ is an $\varepsilon$-approximate-stationary ($\varepsilon$-AS) point of~\eqref{eq:main}, meaning that there exists a Fr\'echet subgradient (defined in Section~\ref{sec:Frechet_subdiff}) of $q$ at $\prox(q,x,\lambda)$ with norm no larger than $\varepsilon$. 
Since the normalized distance from $x$ to $\prox(q,x,\lambda)$, namely $\normt{\lambda^{-1}(x-\prox(q,x,\lambda))}$, is no larger than $\varepsilon$, 
we call $x$  an $\varepsilon$-NS point of~\eqref{eq:main}.

\subsection{Measure of computational cost.}\label{sec:complexity_measure}
In this work, we will develop a first-order method to find an $\varepsilon$-NS point of~\eqref{eq:main}. 
The main computational cost of our method occurs in two aspects: 
\begin{enumerate}[label=\roman*)]
\item computing the {\em primal gradient} $\nabla_x \Phi(x,y)$ and the {\em dual gradient}  $\nabla_y \Phi(x,y)$ at $(x,y)\in\calX\times\calY$,
\item solving certain BPP problems involving the ``simple'' non-smooth functions $r$ and $g$. 
\end{enumerate}
Indeed, due to the proximal-gradient nature of our method, the numbers of solved BPP problems involving $r$ and $g$ are  constant multiples (in fact, at most two) of the numbers of computed primal and dual gradients in our method, respectively. Due to this reason, we measure the computational cost of our method by the complexities of the computed primal and dual gradients, which we call {\em primal gradient complexity} and {\em dual gradient complexity}, respectively. Note that we distinguish between the primal and dual gradient complexities, instead of 
combining them together, mainly because in certain scenarios,  the cost of computing the primal and dual gradients can be different, and/or the cost of solving the BPP problems involving $r$ and involving $g$ can also be different. In these situations, distinguish between the primal and dual gradient complexities allows a more accurate characterization of the computational cost of certain first-order method designed to find an $\varepsilon$-NS point of~\eqref{eq:main}. 



\subsection{Related work.}

Let us review the representative works in the literature, all of which focus on the Hilbertian setting, namely both $\bbX$ and $\bbY$ are finite-dimensional real Hilbert spaces. 

\underline{\em Weakly convex optimization (WCO).} Note that from Assumption~\ref{assump:weak_convex_x}, we can easily show that $f$ is $\gamma$-weakly convex  on $\calX$ (cf.\ Lemma~\ref{lem:weakly_cvx}). Therefore, the problem in~\eqref{eq:main} belongs to the class of WCO problems, which has been studied in several works recently.~\citet{Davis_19b} propose a proximal point method (PPM) for finding an $\varepsilon$-NS point of~\eqref{eq:main} (with $r=\iota_\calX$), where each proximal sub-problem is solved inexactly by the  subgradient method (where the subgradient is in the sense of Fr\'echet; see Section~\ref{sec:Frechet_subdiff}). As another approach,~\citet{Davis_19} propose to find an $\varepsilon$-NS point of~\eqref{eq:main} using the proximal subgradient method directly, without leveraging the PPM framework. Additionally, both works consider the stochastic setting, where the subgradient of $f$ can only be accessed through its unbiased stochastic estimator with finite second moment. Despite the ingenuity and success of these methods, a standing assumption is that at any $x\in\calX$, a subgradient of $f$ can be easily obtained. However, in the case of~\eqref{eq:main}, as we shall see in Lemma~\ref{lem:Lipschitz}, computing a subgradient of $f$ generally requires solving the dual maximization problem in the definition of $f$  {\em exactly}, which may not be possible or accomplished easily at least. 
As such, these methods may not be readily applicable to our problem in~\eqref{eq:main}.


 \underline{\em WCO with max-structure.} In the case where $f$ has the max-structure as in~\eqref{eq:main} with $g\equiv 0$,~\citet{Kong_19} propose an accelerated inexact PPM to find an $\varepsilon$-NS point of the smoothed version of~\eqref{eq:main}, i.e., $\min_{x\in\calX} \,f_\rho(x) + r(x)$, where $f_\rho$  is a smooth approximation of $f$ and is given by $f_\rho:x\mapsto {\max}_{y\in\calY}\;\; \Phi(x,y) - (2\rho)^{-1} \normt{y-\bary}^2$ for all $x\in\calX'$ and some $\bary\in\calY.$
 However, similar to the works~\citep{Davis_19,Davis_19b} reviewed above, the authors assume that the gradient of $f_\rho$ can be obtained easily. As we will see in Lemma~\ref{lem:smooth_f_rho}, this amounts to assume that the maximization problem in the definition of $f_\rho$ can be solved exactly and easily, which may not be the case in general, especially when the structure of either $\Phi(x,\cdot)$ or $\calY$ (or both) is complicated. In another work,~\citet{Theku_19} propose a PPM-based approach to find an $\varepsilon$-NS point of~\eqref{eq:main} by further assuming that $r\equiv 0$.  
 Unlike the aforementioned works, they do not assume certain dual maximization problem can be solved exactly. Instead, they assume the smoothness properties of $\Phi(\cdot,\cdot)$ as in Assumptions~\ref{assump:smooth_Phi} and~\ref{assump:smooth_Phi2}, and solve the proximal sub-problem by combining the Mirror-Prox method~\citep{Nemi_05} and the accelerated gradient method~\citep{Nest_83}. 
However, the analysis in this work critically leverage the inner-product-inducibility of the norm $\normt{\cdot}$ (namely $\normt{x}^2 = \ipt{x}{x}$), and it is not clear (at least to us) how to generalize this approach to the non-Hilbertian setting. 

\underline{\em $\delta$-saddle-stationary ($\delta$-SS) point.} As a final note, since~\eqref{eq:main} can also be viewed as a non-convex-concave minimax optimization problem, i.e., $\min_{x\in\calX}\max_{y\in\calY} \,f(x) + \Phi(x,y) -g(y)$, there exist several works (e.g.,~\citet{Nouie_19,Lu_19,Ostrov_20}) that aim to find a $\delta$-SS point of this minimax problem for any given $\delta>0$ (see~\citet[Section III-A]{Lu_19} for the definition of $\delta$-SS point). In the Hilbertian setting, it can be shown  that one can obtain this point from an $\varepsilon$-NS point of~\eqref{eq:main}, and {vice versa}, by properly choosing $\varepsilon$ and $\delta$ (see e.g.,~\citet[Proposition~4.12]{Lin_19}). However, since we are mainly interested in the minimization problem in~\eqref{eq:main}, rather than the above-mentioned minimax problem, we only focus on the $\varepsilon$-NS point in this work. 

\subsection{Main contributions.}
Our main contributions are threefold. 

First, we propose a primal-dual smoothing framework (namely Algorithm~\ref{algo:PD_smoothing}) for finding an $\varepsilon$-NS point of~\eqref{eq:main}. We analyze the primal and dual gradient complexities of our framework using two approaches: 
(i) the dual-then-primal smoothing approach and (ii) the primal-then-dual smoothing approach. To the best of our knowledge, our framework is the first one that finds an $\varepsilon$-NS point of~\eqref{eq:main} under the non-Hilbertian  setting. Even under the Hilbertian setting and the restrictive case where both $r\equiv 0$ and $g\equiv 0$, 
 the primal and dual gradient complexities of our framework are better than the those of~\citet{Theku_19}, and the improvement is especially significant in the regime where $\gamma\ll L_{xx}$ (recall that $0<\gamma\le L_{xx}$) --- see Table~\ref{table:eps_near_stat} for details. 
 

\renewcommand{\arraystretch}{1.5}
\setlength{\tabcolsep}{8pt}
\begin{table}[t]\centering
\caption{Comparison of 
primal and dual gradient complexities with~\citet{Theku_19} to find an $\varepsilon$-NS point of~\eqref{eq:main}, in the restricted case where $r\equiv 0$, $g\equiv 0$ and both $\bbX$ and $\bbY$ are Hilbert spaces. Note that $0<\gamma\le L_{xx}$.}\vspace{-0pt}\label{table:eps_near_stat}
\footnotesize{
\begin{tabular}{|c|c|c|}\hline
 Algorithms & Primal Oracle Comp. & Dual Oracle Comp.\\\cline{1-3}
\begin{tabular}{@{}c@{}}
Theku.~et~al.~\citep{Theku_19}\\
\end{tabular} & $O\big((L_{xx}+L_{xy}+L_{yy})^2\varepsilon^{-3}\log^2(\varepsilon^{-1})\big)$ & $O\big((L_{xx}+L_{xy}+L_{yy})^2\varepsilon^{-3}\log^2(\varepsilon^{-1})\big)$ \\\hline
Our framework (Algo.~\ref{algo:PD_smoothing}) & $O\big(\textstyle{\sqrt{\gamma L_{xx}}}\big(\textstyle{\sqrt{L_{yy}\gamma}}+L_{xy}\big)\varepsilon^{-3}\log^2(\varepsilon^{-1})\big)$ & $O\big(\gamma\big(\textstyle{\sqrt{L_{yy}\gamma}}+L_{xy}\big)\varepsilon^{-3}\log(\varepsilon^{-1})\big)$\\\hline
\end{tabular}}
\end{table} 

Second, as an important part of our framework, we propose an efficient method for solving a class of (strongly) convex-concave saddle-point problems (SPPs) with primal strong convexity (cf.\ Section~\ref{sec:SPP}). As the workhorse of this method,  we develop a {non-Hilbertian inexact} accelerated proximal gradient (APG) method for strongly convex composite optimization (cf.\ Section~\ref{sec:inexact_APG}) that enjoys certain duality-gap convergence guarantees, and appears to be the first of its kind in the literature. We believe that this inexact APG method may be of independent interest.  

Third, we provide a  variant and an extension of our framework (cf.\ Section~\ref{sec:extension}). We first consider the case where $f$ has a ``simple'' dual structure, and corresponds to the same assumption in~\citet{Kong_19}.  In this case, we show that the primal gradient complexity for finding an $\varepsilon$-NS point of~\eqref{eq:main} has order $O(\varepsilon^{-3}\ln(\varepsilon^{-1}))$, which recovers the result in~\citet{Kong_19} up to a logarithmic factor.  Secondly, we extend our framework to the stochastic case, where the gradients $\nabla_x\Phi(x,y)$ and $\nabla_y\Phi(x,y)$ are only accessible through their stochastic unbiased estimators. Although such an extension is rather straightforward, we show that primal and dual gradient complexities of this extension indeed match the best-known in the literature (see e.g.,~\citet{Raf_18}).

\subsection{Organization.}\label{sec:org}

Our work mainly consists of three parts. 

The first part includes Sections~\ref{sec:prelim} and~\ref{sec:key_lemmas}, and lays the foundation of the whole work. Specifically, in Section~\ref{sec:prelim}, we introduce  several important notions in non-convex analysis and non-Hilbertian optimization (such as Bregman divergence and Bregman proximal projection). In Section~\ref{sec:key_lemmas}, we develop several important lemmas that characterize the (sub-)differential and convexity properties of the function $f$ and its smooth approximation.

\begin{figure}[t]\centering
  \includegraphics[width=.8\textwidth]{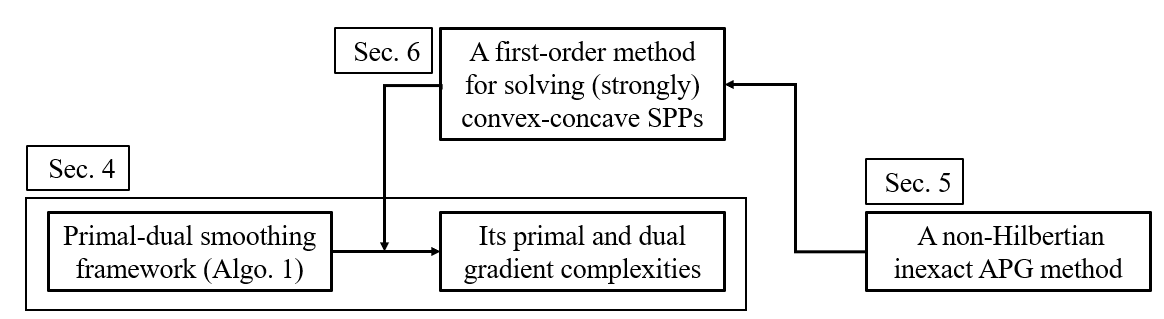}
  \caption{Organization of the main body of this work.}
  \label{fig:org}
\end{figure}

The second part consists of Sections~\ref{sec:pd_smoothing_framework} to~\ref{sec:SPP}, and forms the main body of this work. Its organization is illustrated in Figure~\ref{fig:org}.  Specifically, in Section~\ref{sec:pd_smoothing_framework}, we propose our primal-dual smoothing framework  and analyze its number of iterations. We also analyze the primal and dual gradient complexities of this framework based on the the complexity results of a sub-problem solver that we will develop in Section~\ref{sec:SPP}. This sub-problem solver is essentially an efficient first-order method for solving a class of (strongly) convex-concave SPPs. This method critically leverage a non-Hilbertian inexact APG method that is developed in Section~\ref{sec:inexact_APG}. 

The last part consists of Sections~\ref{sec:extension} and~\ref{sec:conclusion}, wherein we discuss some variants and extensions of our framework, and conclude by pointing out some open problems and future research directions.

\section{Preliminaries.}\label{sec:prelim}
We introduce several important notions that will be used in our analysis. Throughout this work, for any nonempty set $\calX\subseteq\bbX$, we denote its interior by $\inter \calX$, its boundary by $\bdry \calX$ and its closed convex hull by $\clconv\calX$. Also, define the extended real line $\barbbR:= \bbR\cup\{+\infty\}$.

\subsection{Directional derivative, Fr\'echet subdifferential and gradient.} \label{sec:Frechet_subdiff}
Following~\citet{Kruger_03}, given any function $f:\calW\to\bbR$, where $\calW\subseteq\bbX$ is a nonempty open set, 
define its 
(Hadamard) directional derivative  at $x\in\calW$ in the direction of $d\in\bbX$ as 
\begin{equation}
f'(x;d) = \lim_{t\downarrow 0,\;d'\to d}\frac{f(x+td')-f(x)}{t},\label{eq:Hadamard_dir_deriv}
\end{equation}
whenever the limit exists. If the limit in~\eqref{eq:Hadamard_dir_deriv} exists for every $d\in\bbX$, then $f$ is  directionally differentiable at $x$. 
If $f$ is  directionally differentiable at each $x\in\calW$, then we say that $f$ is  directionally differentiable on $\calW$. 

We define the Fr\'echet subdifferential of $f$ at $x\in\calW$ as 
\begin{equation}
\partial f(x) \defeq \left\{\xi\in\bbX^*: \liminf_{h\in\bbX,\;h\to 0} \frac{f(x+h)-f(x)-\lranglet{\xi}{h}}{\normt{h}}\ge 0\right\}.
\end{equation}
In other words, $\xi\in\partial f(x)$ if and only if $f(x+h)\ge f(x) + \lranglet{\xi}{h} + o(\normt{h})$. 
We call the elements in $\partial f(x)$ the Fr\'echet subgradients of $f$ at $x$. 
Note that  $\partial f(x)$ is closed and convex, and  if $f$ is convex, then $\partial f(x)$ is the (convex) subdifferential of $f$ at $x$. We say $f$ is Fr\'echet subdifferentiable at $x$ if $\partial f(x)\ne \emptyset$, and Fr\'echet subdifferentiable on $\calW$ if $\partial f(x)\ne \emptyset$ for all $x\in\calW$.

 We define the  
gradient (or Fr\'echet derivative) of $f$ at $x\in\calW$, denoted by $\nabla f(x)$, as the unique element in $\bbX^*$ that satisfies
\begin{equation}
\lim_{h\in\bbX,\;h\to 0}\frac{f(x+h) - f(x) - \lranglet{\nabla f(x)}{h}}{\normt{h}} = 0. \label{eq:Frechet_deriv}
\end{equation}
In other words, $f(x+h) = f(x) + \lranglet{\nabla f(x)}{h} + o(\normt{h})$.
From the definitions in~\eqref{eq:Hadamard_dir_deriv} and~\eqref{eq:Frechet_deriv}, we see that $f$ is differentiable at $x$ if and only if $d\mapsto f'(x;d)$ is a linear function on $\bbX$, and in this case, we have 
$f'(x;d)= \ipt{\nabla f(x)}{d}$ for all $d\in\bbX$. 
 We say that $f$ is (Fr\'echet) differentiable at $x\in\bbX$ if $\nabla f(x)$ exists, and (Fr\'echet) differentiable on $\calW$ if $\nabla f(x)$ exists for all $x\in\calW$.


\subsection{Distance generating function and Bregman divergence.}\label{sec:DGF_BD}
Let  $\calU$ be a nonempty, convex and closed set in a finite-dimensional real normed space $\bbU$. 
We call $\omega_\calU:\bbU\to\barbbR$ a distance generating function (DGF) on $\calU$ if it is continuous   and 1-strongly-convex on $\calU$ and essentially smooth, i.e., it is continuously differentiable on the interior of its domain (denoted by $\inter \dom \omega_\calU$) and for any sequence $\{u_k\}_{k\ge 0}\subseteq \inter \dom \omega_\calU$ that converges to a boundary point $u\in\bdry \dom \omega_\calU$, we have 
$\normt{\nabla \omega_\calU(u_k)}_*\xrightarrow{k\to +\infty}+\infty$. 
 Based on $\omega_\calU$, let us define its induced Bregman divergence as 
 \begin{equation}
D_{\omega_\calU}(u',u)\defeq \omega_\calU(u') - \omega_\calU(u) - \lranglet{\nabla \omega_\calU(u)}{u'-u}, \quad \forall\,u'\in\dom \omega_\calU, \;\; \forall\,  u\in\inter\dom \omega_\calU. 
\label{eq:Bregman_div}
\end{equation}
Define $\calU^o\defeq \calU\cap \inter\dom \omega_\calU$. Since $\omega_\calU$ is 1-strongly-convex on $\calU$, we have 
 \begin{equation}
 D_{\omega_\calU}(u',u)\ge (1/2)\norm{u'-u}^2, \quad \forall\, u'\in\calU, \;\;\forall\,u\in\calU^o. \label{eq:lb_BD}
 \end{equation}

\subsection{Bregman proximal projection (BPP).} \label{sec:BPP}
For any $u\in\calU^o$ and any convex and closed function $\varphi:\bbU\to\barbbR$, define the BPP of $u$ onto $\calU$ under $\varphi$ and the DGF $\omega_\calU$ (associated with dual vector $\xi\in\bbU^*$ and step-size $\lambda>0$) as the following mapping:
\begin{align}
u\mapsto u^+&\defeq {\argmin}_{u'\in\calU}\;\varphi(u') + \lrangle{\xi}{u'} + \lambda^{-1}D_{\omega_\calU}(u',u)\\
&= {\argmin}_{u'\in\calU}\;\varphi(u') + \lranglet{\barxi}{u'} + \lambda^{-1}\omega_\calU(u') \quad\mbox{where}\quad  \barxi:=\xi-\lambda^{-1}\nabla \omega_\calU(u).
\label{eq:Bregman_projection}
\end{align}
Note that if $\inf_{u\in\calU}\varphi(u)>-\infty$, 
then the minimization problem in~\eqref{eq:Bregman_projection} always has a unique solution in $\calU^o$  (cf.\ \citet[Lemma~A.1]{Zhao_19}). 
We say that the function $\varphi$ has an {\em easily computable} BPP on $\calU$ if there exists a DGF $\omega_\calU$ on $\calU$ such that the minimization problem in~\eqref{eq:Bregman_projection} has a unique {and easily computable} solution in $\calU^o\cap\dom \varphi$, for any $\xi\in\bbU^*$ and $\lambda>0$. For further discussions on BPP, we refer readers to~\citet{Nest_05} and~\citet{Juditsky_12a}.

\subsection{$\varepsilon$-NS point.} \label{sec:eps_near_stat}
Let us formally define the $\varepsilon$-NS point when $\bbX$ is a normed space, using the notions of DGF and BPP in Sections~\ref{sec:DGF_BD} and~\ref{sec:BPP}, respectively.
  Let $\omega_\calX:\bbX\to\barbbR$ be a DGF on $\calX$ (cf.\ Section~\ref{sec:DGF_BD}). Throughout this work, we assume that $\omega_\calX$ satisfy the following additional properties: 
\begin{enumerate}[label=(\roman*)]
\item \label{item:twice_cont} it is twice continuously differentiable  on the interior of its domain (i.e., $\inter\dom \omega_\calX$),
\item \label{item:X_domain} $\calX\subseteq \inter\dom \omega_\calX$ (so that $\calX^o:= \calX\cap \inter\dom \omega_\calX=\calX$),
\item \label{item:grad_lips} its gradient $\nabla \omega_\calX$ is $\beta_\calX$-Lipschitz on $\calX$, where $ \beta_\calX\in [1,+\infty)$.
\end{enumerate}  
(Note that if $\calX$ is bounded, then property~\ref{item:grad_lips} is implied by properties~\ref{item:twice_cont} and~\ref{item:X_domain}, together with the 1-strong-convexity of $\omega_\calX$ on $\calX$.)
 Using a simpler form of BPP in~\eqref{eq:Bregman_projection}  with 
 $\xi=0$, we can define the proximal point at $x\in\inter\dom \omega_\calX$ with function $q$ and step-size $0<\lambda<\gamma^{-1}$ as
 \begin{align}
 \prox(q,x,\lambda)&\defeq {\argmin}_{x'\in\calX}\;\; q(x') + \lambda^{-1}D_{\omega_\calX}(x',x)\label{eq:def_prox_x_BPP}\\
 &= {\argmin}_{x'\in\bbX}\;\; q(x') + \lambda^{-1}D_{\omega_\calX}(x',x), \label{eq:def_prox_x}
 \end{align}
 where~\eqref{eq:def_prox_x} follows from that $\dom q=\dom r= \calX$. (Note that as we will see in Section~\ref{sec:pd_smoothing_framework}, the condition $0<\lambda<\gamma^{-1}$ ensures that the minimization problem in~\eqref{eq:def_prox_x} is strongly convex and hence admits a unique optimal solution.)  
 We call $x\in\calX$ an $\varepsilon$-NS point of~\eqref{eq:main}, if there exists some $0<\lambda<\gamma^{-1}$ such that  
\begin{align}
\normt{\lambda^{-1}(x-\prox(q,x,\lambda))}\le \varepsilon/\beta_\calX. 
\label{eq:eps_near_stat_pt_normed}
\end{align}
Note that the above definition recovers the one when $\bbX$ is Hilbertian, as introduced in Section~\ref{sec:complexity_measure}. Specifically, if we let $\calX = \bbX$ and  $\omega_\calX = (1/2)\normt{\cdot}^2$, 
then  $D_{\omega_\calX}(x',x) = (1/2)\normt{x'-x}^2$ and $\beta_\calX = 1$,  
and hence~\eqref{eq:eps_near_stat_pt_normed} reduces to~\eqref{eq:eps_near_stat_pt}. Similar to the discussions in Section~\ref{sec:complexity_measure}, if~\eqref{eq:eps_near_stat_pt_normed} holds, then $\prox(q,x,\lambda)$ is an $\varepsilon$-AS point of~\eqref{eq:main}. Since $\beta_\calX\ge 1$, the normalized distance  $\normt{\lambda^{-1}(x-\prox(q,x,\lambda))}$ is no larger than $\varepsilon$. As a result, $x$  is called an $\varepsilon$-NS point of~\eqref{eq:main}. 

\section{Important Lemmas.}\label{sec:key_lemmas}
Our algorithmic framework, both in terms of its development and analysis, critically leverage the following lemmas that characterize certain (sub-)differential and convexity properties of the function $f$ (defined in~\eqref{eq:main}) and its smooth approximation. Recall from Assumption~\ref{assump:smooth_Phi} that $\Phi(\cdot,y)$ is differentiable on the open set $\calX'\supseteq \calX$ for all $y\in\calY$.

\subsection{Lemmas on $f$.}
Let us begin by characterizing the directional derivative of $f$ (recall that $f$ is a non-smooth and non-convex function). Our characterization can be regarded as a particular version of Danskin's Theorem~\citep{Daskin_67}. 

\begin{lemma}
\label{lem:Danskin}
The function $f$ in~\eqref{eq:main} is  directionally differentiable on the open set $\calX'$. For any $x\in\calX'$ and $d\in\bbX$, its directional derivative $f'(x;d)$ can be characterized as 
\begin{equation}
f'(x;d)= {\sup}_{y\in\calY^*(x)} \lranglet{\nabla_x \Phi(x,y)}{d}, \quad\mbox{where}\quad \calY^*(x)\defeq {\argmax} _{y\in\calY}\; \Phi(x,y)-g(y)\ne \emptyset. 
\label{eq:dir_deriv_f}
\end{equation}
 In particular, if $\calY^*(x)$ is a singleton, i.e., $\calY^*(x)=\{y^*(x)\}$, 
 then $f$ is  
 differentiable at $x$ and 
\begin{equation}
\nabla f(x)=\nabla_x \Phi(x,y^*(x)).
\end{equation}
\end{lemma}

\proof{Proof.}
See Appendix~\ref{app:proof_Danskin}.
\endproof

Next, we show the local Lipschitz continuity of $f$, and 
characterize its  Fr\'echet subdifferential. 

\begin{lemma}\label{lem:Lipschitz}
The function $f$ is locally Lipschitz on $\calX$. 
In addition, for any $x\in\calX'$, the  Fr\'echet subdifferential $ \partial f(x)=\clconv\{\nabla_x \Phi(x,y):y\in\calY^*(x)\}$. 
\end{lemma}

\proof{Proof.}
See Appendix~\ref{app:proof_local_Lips}.
\endproof

Finally, let us make a simple observation about the weak convexity of $f$, based on Lemma~\ref{lem:Lipschitz}.

\begin{lemma}\label{lem:weakly_cvx}
The function $f$ is $\gamma$-weakly convex on $\calX$, namely, for any $x\in\calX$, we have 
\begin{align}
f(x+d) - f(x) - \lranglet{\xi}{d}\ge -({\gamma}/{2})\norm{d}^2, \quad \forall\, \xi\in\partial f(x).   
\end{align}
\end{lemma}
\proof{Proof.}
Fix any $x\in\calX$ and any $d\in\bbX$. By Assumption~\ref{assump:weak_convex_x}, 
we have that for any $y\in\calY^*(x)$, 
\begin{align*}
f(x+d) - f(x) - \lranglet{\nabla_x \Phi(x,y)}{d} \ge \Phi(x+d,y)  - \Phi(x,y)  - \lranglet{\nabla_x \Phi(x,y)}{d}\ge -({\gamma}/{2})\normt{d}^2.\nt\label{eq:wcvx_grad_Phi}
\end{align*}
By taking convex combination and limit on~\eqref{eq:wcvx_grad_Phi} if necessary, we see that 
\begin{equation}
f(x+d) - f(x) - \lranglet{\xi}{d}\ge -({\gamma}/{2})\normt{d}^2,  \quad \forall\, \xi\in \clconv\{\nabla_x \Phi(x,y):y\in\calY^*(x)\}\eqa \partial f(x), \nn
\end{equation}
where (a) follows from Lemma~\ref{lem:Lipschitz}. This completes the proof. \Halmos
\endproof

\subsection{Lemmas on the dually smoothed $f$.}\label{sec:lemma_f_rho}
Let $\omega_\calY:\bbY\to\barbbR$ be a DGF on $\calY$ (cf.\ Section~\ref{sec:DGF_BD}), and define the {\em $\rho$-dually-smoothed} $f$ as 
\begin{equation}
f_{\rho}(x)= {\max}_{y\in\calY} \;\left[\phi^\rmD_\rho(x,y) \defeq  \Phi(x,y) - g(y) - \rho \omega_\calY(y)\right], \quad \forall\,x\in\calX',\label{eq:dual_reg}
\end{equation}
where $\rho>0$ is called the dual smoothing parameter.  
Define the range of $\omega_\calY$ on $\calY$ as 
\begin{equation}
R_\calY(\omega_\calY)\defeq {\sup}_{y\in\calY}\; \abst{\omega_\calY(y)},  \label{eq:Omega_calY}
\end{equation}
and we have $R_\calY(\omega_\calY)<+\infty$ since $\calY$ is compact and $\omega_\calY$ is continuous on $\calY$. 
Clearly, with this parameter, 
we can uniformly bound the point-wise difference between $f_\rho$ and $f$: for any $x\in\calX$, 
\begin{equation}
\abst{f_\rho(x) - f(x)}\le {\sup}_{y\in\calY}\; \abst{\phi^\rmD_\rho(x,y) - (\Phi(x,y) - g(y))} = {\sup}_{y\in\calY}\; \abst{\rho\omega_\calY(y)} = \rho R_{\calY}(\omega_\calY).\label{eq:bound_f_rho_f}
\end{equation}
In addition, let us define the unique solution to the maximization problem in~\eqref{eq:dual_reg} as $y^*_\rho(x)$, namely
\begin{equation}
y_\rho^*(x)\defeq{\argmax}_{y\in\calY}\; \phi^\rmD_\rho(x,y), \quad \forall\,x\in\calX'. \label{eq:solution_rho}
\end{equation}

Based on these definitions, let us first show that the mapping $y_\rho^*:\calX'\to\calY$ is Lipschitz on $\calX$, {\em even if $\calY$ is unbounded}. (To be clear, all the other results in this paper still assume the boundedness of $\calY$, unless otherwise mentioned.)

\begin{lemma}\label{lem:lips_y*_rho}
Regardless of whether $\calY$ is bounded, the mapping $y_\rho^*$ is $(L_{xy}/\rho)$-Lipschitz on $\calX$.
\end{lemma}

\proof{Proof.}
See Appendix~\ref{app:proof_lips_y*_rho}. \Halmos
\endproof

Based on Lemma~\ref{lem:lips_y*_rho}, we can show the smoothness of $f_\rho$ on $\calX$.

\begin{lemma}
\label{lem:smooth_f_rho}
The function $f_\rho$ is differentiable on $\calX'$ and $\nabla f_\rho (x) = \nabla_x \Phi(x,y^*_\rho(x))$ for all $x\in \calX'$. 
In addition, the gradient $\nabla f_\rho:\calX'\to\bbX^*$ is $L_\rho$-Lipschitz on $\calX$ with $L_\rho \defeq L_{xx}+L_{xy}^2/\rho$. 
\end{lemma}

\proof{Proof.}
See Appendix~\ref{app:proof_sm_f_rho}. 
\endproof

\begin{remark}
Two remarks are in order. First, when the function $(x,y)\mapsto\Phi(x,y)$ is bilinear, i.e., $\Phi(x,y)=\lranglet{\rvA x}{y}$ for some linear operator $\rvA:\bbX\to\bbY^*$, we have $L_{xx}=0$ and $L_{xy}=\normt{\rvA}_{\rm op}$, i.e., the operator norm of $\rvA$. Thus we exactly recover the result in~\citet[Theorem~1]{Nest_05}. Second, note that compared to similar statements about the Lipschitz continuities of $y^*_\rho$ and $\nabla f_\rho$, e.g.,~\citet[Lemma~1]{Sinha_17}, Lemma~\ref{lem:smooth_f_rho} does not require any differentiability assumptions on the function 
$\Phi(x,\cdot):\calY'\to\bbR$ for any $x\in\calX$. 
\end{remark}

Next, using the same reasoning as in Lemma~\ref{lem:weakly_cvx}, we have the following lemma. 

\begin{lemma}\label{lem:weakly_cvx_f_rho}
The function $f_\rho$ is $\gamma$-weakly convex on $\calX$. 
\end{lemma}

Finally, we prove a uniform bound on the distance between $\nabla f_\rho(x)$ and $\partial f(x)$ over $x\in\calX$, for any $\rho >0$.  Given the normed space $(\bbU,\normt{\cdot})$ as in Section~\ref{sec:DGF_BD}, for any point $u\in\bbU$ and any nonempty set $\calU\subseteq \bbU$, define the distance from $u$ to $\calU$ as 
\begin{equation}
\dist(u,\calU):= {\inf}_{u'\in\calU}\; \normt{u-u'} . 
\end{equation}
 Additionally, let us define the diameter of $\calY$ as 
\begin{equation}
D_\calY \defeq {\sup}_{y,y'\in\calY}\;\norm{y-y'}<+\infty. 
\end{equation}

\begin{lemma}\label{lem:dist_grad_subdiff}
For any $\rho>0$ and any $x\in\calX$, we have
\begin{equation}
\dist(\nabla f_\rho(x),\partial f(x))\le L_{xy}\dist(y^*_\rho(x),\calY^*(x))\le L_{xy}D_\calY. 
\label{eq:bound_dist_smoothed_grad}
\end{equation}
\end{lemma}
\proof{Proof.}
See Appendix~\ref{app:proof_dist_grad_subdiff}. \Halmos 
\endproof


\section{A primal-dual smoothing framework for finding an $\varepsilon$-NS point of~\eqref{eq:main}.} \label{sec:pd_smoothing_framework}
As the name suggests, our framework utilizes both ideas of {\em primal smoothing} and {\em dual smoothing}. The notion of dual smoothing has been introduced in Section~\ref{sec:lemma_f_rho}, and let us now introduce primal smoothing. For concreteness, we use the objective function $q$ in~\eqref{eq:main} as an example. From Lemma~\ref{lem:weakly_cvx} and the convexity of the function $r$ on $\calX$, we know that $q$ is $\gamma$-weakly convex on $\calX$. For any $x\in\calX$ and any $0<\lambda<\gamma^{-1}$, let $q^\lambda(x)$ be the optimal value of the minimization problem in~\eqref{eq:def_prox_x}, namely
\begin{align}
q^\lambda(x) &:= {\min}_{x'\in\bbX}\;\; \big[Q^\lambda(x';x):= q(x') + \lambda^{-1}D_{\omega_\calX}(x',x)\big]\label{eq:primal_sm_q} \\
& = Q^\lambda(\prox(q,x,\lambda);x),  \label{eq:q_Q_prox}
\end{align}
where~\eqref{eq:q_Q_prox} follows from the definition of $\prox(q,x,\lambda)$ in~\eqref{eq:def_prox_x}. 
(Note that since $\lambda^{-1}>\gamma$ and that $\omega_\calX$ is $1$-strong convex  on $\calX$, the function $Q^\lambda(\cdot;x)$ is ($\lambda^{-1}-\gamma$)-strongly convex on $\calX$, and hence $\prox(q,x,\lambda)$ in~\eqref{eq:def_prox_x} is indeed unique.) Using the definition of 
$D_{\omega_\calX}(x',x)$ in~\eqref{eq:Bregman_div} and invoking Lemma~\ref{lem:Danskin}, we see that the function  $q^\lambda:\inter\dom \omega_\calX\to\bbR$ is 
differentiable on $\inter\dom \omega_\calX\supseteq\calX$, and 
\begin{align}
\nabla q^\lambda(x) &= \nabla^2 \omega_\calX(x) \lambda^{-1}\big(x-\prox(q,x,\lambda)\big), \quad \forall\,x\in\inter\dom \omega_\calX.  \label{eq:nabla_q^lambda}
\end{align}
For this reason, we can regard the operation in~\eqref{eq:primal_sm_q}, which transforms $q$ to $q^\lambda$, as the primal smoothing procedure on $q$, and call the resulting function $q^\lambda$ the {\em $\lambda$-primally-smoothed} $q$.

Now, let us define $q_\rho:= f_\rho + r$ for some $\rho>0$, where $f_\rho$ is the $\rho$-dually-smoothed $f$. As such, $q_\rho$  can be regarded as the $\rho$-dually-smoothed $q$. 
Similar to the above, we can define the $\lambda$-primally-smoothed $q_\rho$ for some $0<\lambda<\gamma^{-1}$, denoted by $q_\rho^\lambda:\inter\dom \omega_\calX\to\bbR$, as 
\begin{align}
q_\rho^\lambda(x) &:= {\min}_{x'\in\bbX}\;\; \big[Q_\rho^\lambda(x';x):= q_\rho(x') + \lambda^{-1}D_{\omega_\calX}(x',x)\big]\label{eq:primal_sm_q_rho} \\
& = Q_\rho^\lambda(\prox(q_\rho,x,\lambda);x), \quad \mbox{where}\quad \prox(q_\rho,x,\lambda):= {\argmin}_{x'\in\bbX}\;\; Q_\rho^\lambda(x';x).    \label{eq:q_rho_Q_prox}
\end{align}

\begin{algorithm}[t!]
\caption{Primal dual smoothing framework 
} \label{algo:PD_smoothing}
\begin{algorithmic}
\State {\bf Input}: Accuracy parameter $\eta>0$, smoothing parameters $\lambda \in (0,\gamma^{-1})$ and $\rho = 2\eta/R_\calY(\omega_\calY)$
\State {\bf Initialize}: $k:=0$, $x_1\in\calX$
\State {\bf Repeat:} 
\begin{enumerate}
{\setlength\itemindent{25pt} \item $k:= k+1$.}
{\setlength\itemindent{25pt} \item \label{item:inexact_PPM} Find $x_{k+1}\in\calX$ such that $Q_\rho^\lambda(x_{k+1};x_k)\le q_\rho^\lambda(x_k) + \eta$.} 
\end{enumerate}
\State {\bf Until:} 
\begin{equation}
\textstyle\normt{x_{k+1}-x_k}\le \sqrt{2\eta/(\lambda^{-1}-\gamma)}. 
\label{eq:conv_crit}
\end{equation}
\end{algorithmic}
\end{algorithm}

The minimization problem in~\eqref{eq:primal_sm_q_rho} indeed suggests an iterative scheme for finding an $\varepsilon$-NS of~\eqref{eq:main}, which forms the basis of our framework shown in Algorithm~\ref{algo:PD_smoothing}. Specifically, we start with any point in $\calX$, and in each iteration, 
we {\em approximately} solve the minimization problem in \eqref{eq:primal_sm_q_rho} with some accuracy $\eta>0$, and denote this $\eta$-optimal solution 
as $x_{k+1}$ (cf.\ Step~\ref{item:inexact_PPM}). Throughout all the iterations, we fix the primal smoothing parameter $\lambda \in (0,\gamma^{-1})$ and the dual smoothing parameters $\rho = 2\eta/R_\calY(\omega_\calY)$, where  $R_\calY(\omega_\calY)$ 
is defined in~\eqref{eq:Omega_calY}. The appropriate choices of $\eta$ and $\lambda$ will become apparent after our analysis. 
We terminate the algorithm once the distance between two successive iterates falls below $\sqrt{2\eta/(\lambda^{-1}-\gamma)}$ (see the termination criterion in~\eqref{eq:conv_crit}). 

Before ending the description of our framework, let us notice that the termination criterion~\eqref{eq:conv_crit} is easily checkable --- in fact, it is solely based on  the distances between successive iterates. This is in contrast to the convergence criteria in previous works (e.g.,~\citet{Kong_19} and~\citet{Theku_19})  that involve quantities like $q(x_k)$ or $\nabla f_\rho(x_k)$, whose 
evaluation 
typically requires solving certain dual optimization problems.  Hence these convergence criteria are harder to check than ours.

In the following, we analyze Algorithm~\ref{algo:PD_smoothing} using two different approaches. 
Let $K\ge 1$ denote the iteration that Algorithm~\ref{algo:PD_smoothing} terminates, so that Algorithm~\ref{algo:PD_smoothing} outputs $x_K$. 
For either approach, we derive the choice of $\eta$ (as a function of $\lambda$, $\varepsilon$ and $\beta_\calX$) such that
\begin{equation}
\normt{x_K-\prox(q,x_K,\lambda)}\le \varepsilon\lambda/\beta_\calX, \label{eq:x_k_eps_NS}
\end{equation}
which implies that $x_K\in\calX$ is an $\varepsilon$-NS point of~\eqref{eq:main} (cf.\ Section~\ref{sec:eps_near_stat}). Before presenting our analysis, we first show that if~\eqref{eq:x_k_eps_NS} holds, then 
i) $\prox(q,x_K,\lambda)$ is an $\varepsilon$-AS point, namely, there exists a Fr\'echet subgradient of $q$ at $\prox(q,x_K,\lambda)$ with norm not exceeding $\varepsilon$ and ii) 
the gradient of $q^\lambda$ at $x_K$ has a small norm (not exceeding $\varepsilon$), where $q^\lambda$ is the $\lambda$-primally-smoothed $q$ as defined in~\eqref{eq:primal_sm_q}.  

\begin{lemma}\label{lem:eps_subdiff}
For any $\varepsilon>0$ and any $0< \lambda< \gamma^{-1}$, if $x\in\calX$ satisfies $\normt{x-\prox(q,x,\lambda)}\le \varepsilon\lambda/\beta_\calX$, then we have 
\begin{equation}
\dist\big(0,\partial q\big(\prox(q,x,\lambda)\big)\big)\le \varepsilon \qquad \mbox{and}\qquad \normt{\nabla q^\lambda(x)}_*\le \varepsilon. 
\end{equation}
\end{lemma}
\proof{Proof.}
Applying the first-order optimality condition to the definition of $\prox(q,x,\lambda)$ (cf.~\eqref{eq:def_prox_x}), we have that for any $0< \lambda< \gamma^{-1}$, 
\begin{equation}
\lambda^{-1}\big({\nabla \omega_\calX(x)-\nabla \omega_\calX\big(\prox(q,x,\lambda)\big)}\big)\in\partial q\big(\prox(q,x,\lambda)\big). 
\end{equation} 
As a result, using the $\beta_\calX$-Lipschitz continuity of $\nabla \omega_\calX$ on $\calX$, we have\begin{align*}
\dist\big(0,\partial q\big(\prox(q,x,\lambda)\big)\big) &= {\inf}_{\xi\in \partial q(\prox(q,x,\lambda))}\; \normt{\xi}_*\\
&\le \lambda^{-1}\normt{\nabla \omega_\calX(x)-\nabla \omega_\calX\big(\prox(q,x,\lambda)\big)}_*\\
&\le \lambda^{-1}\beta_\calX\normt{x-\prox(q,x,\lambda)} 
\le \varepsilon. 
\end{align*}
Using the $\beta_\calX$-Lipschitz continuity of $\nabla \omega_\calX$ on $\calX$ again, we can easily show that the operator norm of $\nabla^2 \omega(x):\bbX\to\bbX^*$ (denoted by $\normt{\nabla^2 \omega(x)}_{\rm op}$) is uniformly bounded on $\calX$ by $\beta_\calX$, namely 
\begin{equation}
\normt{\nabla^2 \omega(x)}_{\rm op}:= {\sup}_{z\in\bbX}\; \normt{\nabla^2 \omega(x)z}_* \le \beta_\calX, \quad\forall\, x\in\calX. 
\end{equation}
Therefore, by the definition of $\nabla q^\lambda(x)$ in~\eqref{eq:nabla_q^lambda}, we have 
\begin{align}
\normt{\nabla q^\lambda(x)}_* \le  \lambda^{-1}\normt{\nabla^2 \omega_\calX(x)}_{\rm op} \norm{x-\prox(q,x,\lambda)}\le \lambda^{-1}\beta_\calX \norm{x-\prox(q,x,\lambda)}\le \varepsilon. \tag*{\Halmos}
\end{align}
\endproof


\subsection{Approach I: Dual-then-primal smoothing.}\label{sec:perspective_I}

We analyze Algorithm~\ref{algo:PD_smoothing} by regarding it as an inexact proximal-point framework for finding a near-stationary point of $q_\rho$, namely the $\rho$-dually-smoothed $q$.  
Indeed, we first show that $x_K\in\calX$ is a near-stationary point of $q_\rho$, 
by bounding $\normt{x_K - \prox(q_\rho,x_K,\lambda)}$, 
and then bound the distance between the proximal points $\prox(q,x_K,\lambda)$ and $\prox(q_\rho,x_K,\lambda)$, namely  $\normt{\prox(q,x_K,\lambda) - \prox(q_\rho,x_K,\lambda)}$. 
 These two bounds together yield a bound on $\normt{x_K - \prox(q,x_K,\lambda)}$, thereby showing that  $x_K$ is a near-stationary point of $q$. These steps are formalized  below.
\begin{lemma}\label{lem:dual_then_primal1}
In Algorithm~\ref{algo:PD_smoothing}, for any $\eta>0$, we have $\norm{x_K - \prox(q_\rho,x_K,\lambda)}\le 2\sqrt{2\eta/(\lambda^{-1}-\gamma)}$. 
\end{lemma} 
 \proof{Proof.} 
By the fact that $\prox(q_\rho,x_K,\lambda) = \argmin_{x'\in\calX}\; Q_\rho^\lambda(x';x_K)$ (cf.~\eqref{eq:q_rho_Q_prox}) and the $(\lambda^{-1}-\gamma)$-strong-convexity of the function $Q_\rho^\lambda(\cdot;x_K)$ on $\calX$, we have 
\begin{align}
\frac{\lambda^{-1}-\gamma}{2}\norm{x_{K+1} - \prox(q_\rho,x_K,\lambda)}^2&\le Q_\rho^\lambda(x_{K+1};x_K) - Q_\rho^\lambda(\prox(q_\rho,x_K,\lambda);x_K)\le \eta,\label{eq:sc_Q_rho}
\end{align}
 where the second inequality follows from the condition in Step~\ref{item:inexact_PPM} and~\eqref{eq:q_Q_prox}.  
This implies that $\norm{x_{K+1} - \prox(q_\rho,x_K,\lambda)}\le \sqrt{2\eta/(\lambda^{-1}-\gamma)}$.  On the other hand, by the definition of $K$, we have $\normt{x_{K+1}-x_K}\le \sqrt{2\eta/(\lambda^{-1}-\gamma)}$.  This completes the proof. \Halmos
\endproof

\begin{lemma}\label{lem:dual_then_primal2}
For any  $\eta>0$ 
and any $x\in\calX$, we have
\begin{equation}
\normt{\prox(q_\rho,x,\lambda)-\prox(q,x,\lambda)}\le 2\sqrt{\rho R_\calY(\omega_\calY)/(\lambda^{-1}-\gamma)}. 
\end{equation}
In particular, if $x = x_K$ and $\rho = 2\eta/R_\calY(\omega_\calY)$ (as in Algorithm~\ref{algo:PD_smoothing}), then we have 
\begin{equation}
\normt{\prox(q_\rho,x_K,\lambda)-\prox(q,x_K,\lambda)}\le 2\sqrt{2\eta/(\lambda^{-1}-\gamma)}. 
\end{equation}
\end{lemma}
\proof{Proof.}
By the fact that $\prox(q,x,\lambda) = \argmin_{x'\in\calX}\; Q^\lambda(x';x)$ (cf.~\eqref{eq:def_prox_x_BPP}) and the $(\lambda^{-1}-\gamma)$-strong-convexity of $Q^\lambda(\cdot;x)$ on $\calX$, we have
\begin{align}
\frac{\lambda^{-1}-\gamma}{2} \normt{\prox(q_\rho,x,\lambda)-\prox(q,x,\lambda)}^2 \le Q^\lambda(\prox(q_\rho,x,\lambda);x) - Q^\lambda(\prox(q,x,\lambda);x).\label{eq:bound_two_prox}
\end{align}
In addition, 
\begin{align}
Q^\lambda(\prox(q_\rho,x,\lambda);x) - Q_\rho^\lambda(\prox(q_\rho,x,\lambda);x) &= f(\prox(q_\rho,x,\lambda)) - f_\rho(\prox(q_\rho,x,\lambda)) \nn\\
&\lea \rho\Omega_{\calY}(\omega_\calY),\label{eq:bound_Q_diff1}\\
Q_\rho^\lambda(\prox(q_\rho,x,\lambda);x) - Q^\lambda(\prox(q,x,\lambda);x) &\leb Q_\rho^\lambda(\prox(q,x,\lambda);x) - Q^\lambda(\prox(q,x,\lambda);x)\nn\\
& = f_\rho(\prox(q,x,\lambda)) - f(\prox(q,x,\lambda))\nn\\ 
&\lec \rho\Omega_{\calY}(\omega_\calY), \label{eq:bound_Q_diff2}
\end{align}
where (a) and (c) follow from~\eqref{eq:bound_f_rho_f} 
and (b) follows from~\eqref{eq:q_rho_Q_prox}.  
Now, by combining
~\eqref{eq:bound_two_prox},~\eqref{eq:bound_Q_diff1} and~\eqref{eq:bound_Q_diff2}, we complete the proof. \Halmos  
\endproof
Combining Lemmas~\ref{lem:dual_then_primal1} and~\ref{lem:dual_then_primal2}, we see that 
\begin{equation}
\norm{x_K - \prox(q,x_K,\lambda)}\le 4\sqrt{2\eta/(\lambda^{-1}-\gamma)},  \label{eq:dist_x_K_prox}
\end{equation}
and hence  
we have the following theorem.
\begin{theorem}\label{thm:conv_PPA}
In Algorithm~\ref{algo:PD_smoothing}, 
for any $\varepsilon>0$, if we set the accuracy parameter 
\begin{equation}
\eta = \varepsilon^2\lambda(1-\gamma\lambda)/(32\beta_\calX^2), \label{eq:choice_eta}
\end{equation}
then $\norm{x_K - \prox(q,x_K,\lambda)}\le \varepsilon\lambda/\beta_\calX$, 
meaning that 
$x_K$ is an $\varepsilon$-NS point of~\eqref{eq:main}. 
\end{theorem}

\subsection{Approach II: Primal-then-dual smoothing.}\label{sec:perspective_II}

Alternatively, we can directly view Algorithm~\ref{algo:PD_smoothing} as an inexact proximal-point framework for finding a near-stationary point of~\eqref{eq:main}. Specifically, we will show that Step~\ref{item:inexact_PPM} in Algorithm~\ref{algo:PD_smoothing} implies that $x_{k+1}$ is also an approximate solution for the optimization problem in~\eqref{eq:primal_sm_q}, namely  
\begin{equation}
Q^\lambda(x_{k+1};x_k)\le Q^\lambda(\prox(q,x_k,\lambda);x_k) + 5\eta, \quad \forall\,k\ge 1. \label{eq:prox_Q^lambda} 
\end{equation}
The reason that we (approximately) solve the ``dually smoothed'' optimization problem in~\eqref{eq:primal_sm_q_rho}, instead of that in~\eqref{eq:primal_sm_q}, is because the former enjoys certain smoothness properties and hence can be more efficiently solved using (primal-dual) first-order methods (cf.\ Section~\ref{sec:SPP}). Then, based on~\eqref{eq:prox_Q^lambda}, we can easily  arrive at~\eqref{eq:dist_x_K_prox}.  The details are shown below.


\begin{lemma}\label{lem:primal_then_dual}
In Algorithm~\ref{algo:PD_smoothing}, for all $k\ge 1$, we have~\eqref{eq:prox_Q^lambda} and 
\begin{align}
\norm{x_{k+1} - \prox(q,x_k,\lambda)}\le \sqrt{10\eta/(\lambda^{-1}-\gamma)}.\label{eq:dist_prox}
\end{align}
\end{lemma}
\proof{Proof.}
Note that by the definitions of $q^\lambda$ and $q_\rho^\lambda$ in~\eqref{eq:primal_sm_q} and~\eqref{eq:primal_sm_q_rho}, respectively, for any $x\in\calX$,  we have 
\begin{align}
& {\sup}_{x'\in\calX} \, \abst{Q^\lambda(x';x) - Q_\rho^\lambda(x';x)} = {\sup}_{x'\in\calX} \, \abst{f(x') - f_\rho(x')}\le \rho R_\calY(\omega_\calY), \label{eq:diff_qlam_qrholam}\\
&\abst{q^\lambda(x) - q_\rho^\lambda(x)} \le {\sup}_{x'\in\calX} \, \abst{Q^\lambda(x';x) - Q_\rho^\lambda(x';x)}\le \rho R_\calY(\omega_\calY),\label{eq:diff_qlam_qrholam2} 
\end{align}
where the inequality in~\eqref{eq:diff_qlam_qrholam} follows from~\eqref{eq:bound_f_rho_f}. By Step~\ref{item:inexact_PPM} and~\eqref{eq:diff_qlam_qrholam2}, we have
\begin{equation}
Q_\rho^\lambda(x_{k+1};x_k)\le q_\rho^\lambda(x_k) + \eta\le  q^\lambda(x_k) +\rho R_\calY(\omega_\calY)  + \eta.  \label{eq:Q_rho1}
\end{equation}
On the other hand, by~\eqref{eq:diff_qlam_qrholam}, we have
\begin{equation}
Q^\lambda(x_{k+1};x_k) - \rho R_\calY(\omega_\calY)  \le Q_\rho^\lambda(x_{k+1};x_k). \label{eq:Q_rho2}
\end{equation}
Combining~\eqref{eq:Q_rho1} and~\eqref{eq:Q_rho2}, and using that $\rho = 2\eta/R_\calY(\omega_\calY)$,  we have
\begin{equation}
Q^\lambda(x_{k+1};x_k)\le q^\lambda(x_k) +2\rho R_\calY(\omega_\calY)  + \eta = Q^\lambda\big(\prox(q,x_k,\lambda);x_k\big)+ 5\eta.\label{eq:bound_Q_lambda_q_lambda}
\end{equation}
By the fact that $\prox(q,x_k,\lambda) = \argmin_{x'\in\calX}\; Q^\lambda(x';x_k)$ (cf.~\eqref{eq:def_prox_x_BPP}) and the $(\lambda^{-1}-\gamma)$-strong-convexity of $Q^\lambda(\cdot;x)$ on $\calX$, we therefore have 
\begin{align}
\norm{x_{k+1} - \prox(q,x_k,\lambda)}\le \sqrt{10\eta/(\lambda^{-1}-\gamma)}.  \tag*{\Halmos} 
\end{align}
\endproof

Now, from~\eqref{eq:conv_crit}, we know that $x_K$ 
satisfies $\normt{x_{K+1}-x_K}\le \sqrt{2\eta/(\lambda^{-1}-\gamma)}$. This, together with Lemma~\ref{lem:primal_then_dual}, implies~\eqref{eq:dist_x_K_prox}, which then leads to Theorem~\ref{thm:conv_PPA}.

\subsection{Bound on the number of iterations of Algorithm~\ref{algo:PD_smoothing}.}\label{sec:finite_iterations}
The simple structure of Algorithm~\ref{algo:PD_smoothing} enables us to easily derive a bound on the number of iterations of Algorithm~\ref{algo:PD_smoothing}.

\begin{theorem}\label{thm:total_num_iter}
For any $\lambda\in(0.8\gamma^{-1},\gamma^{-1})$, Algorithm~\ref{algo:PD_smoothing}  terminates in at most 
\begin{equation}
\left\lceil \frac{q(x_1)-q^*}{(\frac{\gamma\lambda}{1-\gamma\lambda}-4)\eta}\right\rceil+1 
\end{equation}
iterations, where $q^*>-\infty$ is the optimal value of~\eqref{eq:main}.  
\end{theorem}
\proof{Proof.}
Note that if~\eqref{eq:conv_crit} is not satisfied at iteration $k$, then $\normt{x_{k+1}-x_k}> \sqrt{2\eta/(\lambda^{-1}-\gamma)}$. Therefore, using~\eqref{eq:lb_BD} and~\eqref{eq:diff_qlam_qrholam2}, we have
 \begin{equation}
 Q_\rho^\lambda(x_{k+1};x_k)\ge q_\rho(x_{k+1}) + (2\lambda)^{-1} \normt{x_{k+1}-x_k}^2 \ge  q(x_{k+1}) - \rho R_\calY(\omega_\calY) + \eta/(1-\gamma\lambda). \label{eq:outer_iter_1}
 \end{equation}
 On the other hand, using Step~\ref{item:inexact_PPM} and~\eqref{eq:diff_qlam_qrholam2}, we have 
\begin{align}
Q_\rho^\lambda(x_{k+1};x_k)\le q_\rho^\lambda(x_k) + \eta \le q^\lambda(x_k)  + \rho R_\calY(\omega_\calY)+ \eta\le  q(x_k)  + \rho R_\calY(\omega_\calY)+ \eta, \label{eq:outer_iter_2}
\end{align}
where the last step follows from $q^\lambda(x)\le Q^\lambda(x;x)= q(x)$ for all $x\in\calX$. 
Combining~\eqref{eq:outer_iter_1}  and~\eqref{eq:outer_iter_2}, and plugging in the value $\rho = 2\eta/R_\calY(\omega_\calY)$, we have
\begin{align}
q(x_{k+1}) \le q(x_k)  + 2\rho R_\calY(\omega_\calY) - \frac{\gamma\lambda}{1-\gamma\lambda}\eta = q(x_k)   - \left(\frac{\gamma\lambda}{1-\gamma\lambda}-4\right)\eta. 
\end{align}
Summing over $k = 1,\ldots,K-1$, we have 
\begin{align}
\left(\frac{\gamma\lambda}{1-\gamma\lambda}-4\right)\eta (K-1)\le q(x_{1}) - q(x_{K})\le q(x_{1}) - q^*. 
\end{align}
Rearranging and we complete the proof. 
\Halmos
\endproof

Based on Theorem~\ref{thm:conv_PPA} and Theorem~\ref{thm:total_num_iter}, we have the following corollary.

\begin{corollary}\label{cor:total_num_iter}
For any $\lambda\in(0.8\gamma^{-1},\gamma^{-1})$ and any $\varepsilon>0$, if we set $\eta$ as in~\eqref{eq:choice_eta}, then Algorithm~\ref{algo:PD_smoothing}  returns an $\varepsilon$-NS point of~\eqref{eq:main} in no more than $\barK$ iterations, where 
\begin{align}
\barK\defeq \left\lceil \frac{32\beta_\calX^2(q(x_1)-q^*)}{5\varepsilon^2\lambda(\gamma\lambda-0.8)}\right\rceil+1.\label{eq:barK}
\end{align}
\end{corollary}

\subsection{Choice of $\lambda$ and the primal and dual gradient complexities of Algorithm~\ref{algo:PD_smoothing}.} \label{sec:choice_lambda}

From Corollary~\ref{cor:total_num_iter}, we see that in order to reduce the  bound on the number of iterations of Algorithm~\ref{algo:PD_smoothing} (namely $\barK$), we should choose $\lambda$ to be as close as $\gamma^{-1}$ as possible. However, note that in Step~\ref{item:inexact_PPM} we need to solve the minimization problem in~\eqref{eq:primal_sm_q_rho} with accuracy $\eta$. Since the choice of the accuracy parameter $\eta$ in~\eqref{eq:choice_eta} is proportional to $\lambda(1-\gamma\lambda)$, and the strong-convexity parameter of the function $Q_\rho^\lambda(\cdot;x_k)$ is $\lambda^{-1}-\gamma$, if $\lambda$ is close to $\gamma^{-1}$, then the  strong-convexity parameter becomes very small (in other words, the problem becomes ill-conditioned) and the accuracy becomes very high. Intuitively, this should  result in a high computational cost in Step~\ref{item:inexact_PPM}, 
and we formalize this intuition below. 

\underline{\em Primal and dual gradient complexities for implementing Step~\ref{item:inexact_PPM}.} 
Indeed, since $\lambda<\gamma^{-1}$ and $\omega_\calX$ and $\omega_\calY$ are 1-strongly-convex on $\calX$ and $\calY$, respectively, we observe that the minimization problem~\eqref{eq:primal_sm_q_rho} in Step~\ref{item:inexact_PPM} is indeed a (strongly) convex-concave SPP:
\begin{align}
{\min}_{x\in\calX}\;{\max}_{y\in\calY} \;\; r(x) + \lambda^{-1}D_{\omega_\calX}(x;x_k) + \Phi(x,y) - g(y) - \rho \omega_\calY(y). \label{eq:cvx_ccv_SPP}
\end{align}
In the next two sections (namely Sections~\ref{sec:inexact_APG} and~\ref{sec:SPP}), we will develop an efficient first-order method for solving a general class of (strongly) convex-concave SPP that subsumes~\eqref{eq:cvx_ccv_SPP} as a special case. Specifically, this method finds $x_{k+1}\in\calX$ that satisfies the $\eta$-optimality condition in Step~\ref{item:inexact_PPM} with primal gradient complexity 
\begin{equation}
C_\rmp(\eta) = O\left(\bigg(\sqrt{\frac{L_{yy}}{\rho}}+\frac{L_{xy}}{\sqrt{(\lambda^{-1} - \gamma)\rho}}\bigg)\sqrt{\frac{L_{xx}+\gamma\beta_\calX}{\lambda^{-1} - \gamma}}\ln^2\left(\frac{1}{\eta\rho}\right)\right) \label{eq:primal_grad_comp}
\end{equation}
and dual gradient complexity
\begin{equation}
C_\rmd(\eta) = O\left(\bigg(\sqrt{\frac{L_{yy}}{\rho}}+\frac{L_{xy}}{\sqrt{(\lambda^{-1} - \gamma)\rho}}\bigg)\ln\left(\frac{1}{\eta\rho}\right)\right),\label{eq:dual_grad_comp}
\end{equation}
where the primal and dual gradient complexities are defined in Section~\ref{sec:complexity_measure}. 
From~\eqref{eq:primal_grad_comp} and~\eqref{eq:dual_grad_comp}, it is clear that if $\lambda$ is close to $\gamma^{-1}$, then both the primal and dual gradient complexities for solving~\eqref{eq:primal_sm_q_rho} with accuracy $\eta$ 
becomes very high. In fact, to reduce both $C_\rmp(\eta)$ and $C_\rmd(\eta)$, we wish to choose $\lambda$ as small as possible. 

The analysis above reveals a trade-off in the choice of $\lambda$, that is, between reducing the number of iterations of Algorithm~\ref{algo:PD_smoothing} and reducing the computational cost of implementing Step~\ref{item:inexact_PPM}. From Theorem~\ref{thm:conv_PPA} (and Corollary~\ref{cor:total_num_iter}), we know the legitimate range of $\lambda$ is $(0.8\gamma^{-1}, \gamma^{-1})$, and hence a natural choice of $\lambda$ would be the mid-point of this interval, namely  $\lambda=0.9\gamma^{-1}.$
Based on this choice of $\lambda$, as well as the choice of $\eta$ in Theorem~\ref{thm:conv_PPA},  we have the following result. 

\begin{corollary}
For any $\varepsilon>0$, if we choose $\lambda=0.9\gamma^{-1}$ and set $\eta$ as in~\eqref{eq:choice_eta}, then Algorithm~\ref{algo:PD_smoothing} finds an $\varepsilon$-NS point of~\eqref{eq:main} with primal gradient complexity
\begin{equation}
T_\rmp(\varepsilon)=O\left(\sqrt{\gamma L_{xx}}\big({\sqrt{L_{yy}\gamma}}+L_{xy}\big)\varepsilon^{-3}\ln(\varepsilon^{-1})^2 \right)
\end{equation}
and dual gradient complexity
\begin{equation}
T_\rmd(\varepsilon)=O\Big(\gamma\big(\textstyle{\sqrt{L_{yy}\gamma}}+L_{xy}\big)\varepsilon^{-3}\ln(\varepsilon^{-1})\Big).
\end{equation}
\end{corollary}
\proof{Proof.}
Indeed, if we substitute the value of $\rho$ as in Algorithm~\ref{algo:PD_smoothing}, the value of $\eta$  as in~\eqref{eq:choice_eta} and $\lambda=0.9\gamma^{-1}$ into the definitions of $C_\rmp(\eta)$ and $C_\rmd(\eta)$ in~\eqref{eq:primal_grad_comp} and~\eqref{eq:dual_grad_comp}, respectively, 
then we have
\begin{align}
\begin{split}
C_\rmp(\eta) &= O\left(\sqrt{L_{xx}/\gamma}\big({\sqrt{L_{yy}\gamma}}+L_{xy}\big)\varepsilon^{-1}\ln(\varepsilon^{-1})^2 \right)\quad\mbox{and}\\
 C_\rmd(\eta) &= O\Big(\big(\textstyle{\sqrt{L_{yy}\gamma}}+L_{xy}\big)\varepsilon^{-1}\ln(\varepsilon^{-1})\Big).
 \end{split}\label{eq:Cp_Cd}
\end{align}
Furthermore, Corollary~\ref{cor:total_num_iter} states that if $\eta$ is set as in~\eqref{eq:choice_eta}, then in order for Algorithm~\ref{algo:PD_smoothing} to return an $\varepsilon$-NS point of~\eqref{eq:main}, the number of iterations is bounded by 
$\barK=O(\varepsilon^{-2}\gamma)$ (since $\lambda=0.9\gamma^{-1}$). Combining this bound  with~\eqref{eq:Cp_Cd}, we complete the proof. \Halmos
\endproof

As promised above, in the next two sections (namely Sections~\ref{sec:inexact_APG} and~\ref{sec:SPP}), we will develop an efficient first-order method for solving the SPP in~\eqref{eq:cvx_ccv_SPP}. Specifically, Section~\ref{sec:inexact_APG} is devoted to a new non-Hilbertian inexact APG method for strongly convex composite optimization,  which forms the basis of the actual first-order method for  solving~\eqref{eq:cvx_ccv_SPP} that will be developed in Section~\ref{sec:SPP} . 


\section{A non-Hilbertian inexact APG method.}\label{sec:inexact_APG} 

Let us consider the following (strongly) convex optimization problem: 
\begin{align}
(P):\quad P^*:= {\min}_{u\in\calU} \;\big\{P(u) \defeq h(u) + \zeta(u) + \mu\omega_\calU(u)\big\}, \label{eq:primal_minimization}
\end{align}
where $\calU$ is a nonempty, convex and closed set in the normed space $\bbU$ as given in Section~\ref{sec:DGF_BD}, the function $h$ is $L_h$-smooth on $\calU$, namely it is differentiable on some open set $\calU'\supseteq\calU$ and $\nabla h$ is $L_h$-Lipschitz on $\calU$, and 
the function $\zeta$ has an easily computable BPP on $\calU$ with DGF $\omega_\calU$ (cf.\ Section~\ref{sec:BPP}). Both functions $h$ and $\zeta$ are convex on $\calU$. In addition, by the 1-strong-convexity of $\omega_\calU$ on $\calU$, the objective function $P$ is $\mu$-strongly convex on $\calU$, where  $\mu\ge 0$.  We assume that~\eqref{eq:primal_minimization} has an optimal solution $u^*$, which necessarily lies in $\calU^o (= \calU\cap\inter\dom \omega_\calU)$, and hence $P(u^*) = P^*$. (Note that if $\mu>0$, $u^*\in \calU^o$ is guaranteed to exist and is unique.)

In particular, we are interested in the case where $h$ has the following max-structure:
\begin{equation}
h(u) = {\max}_{v\in\calV}\; \barPsi(u,v),  \quad \forall u\in\calU', \label{eq:h_max_struct}
\end{equation}
where $\calV$ is a nonempty, compact and convex set contained in some open set $\calV'$ and $\barPsi:\calU'\times\calV'\to\bbR$ is jointly continuous on $\calU'\times\calV'$. 
In addition, for any $v\in\calV$, $\barPsi(\cdot,v)$ is convex on $\calU$ and differentiable on $\calU'$, and for any $u,u'\in\calU$ and $v,v'\in\calV$, we have  
\begin{align}
&\normt{\nabla_u \barPsi(u,v) - \nabla_u \barPsi(u',v)}_* \le L_{uu}\normt{u-u'},\label{eq:L_uu}\\
&\normt{\nabla_u \barPsi(u,v) - \nabla_u \barPsi(u,v')}_* \le L_{uv}\normt{v-v'}.\label{eq:L_uv}
\end{align} 
Also, $\barPsi(u,\cdot)$ is $\rho$-strongly concave on $\calV$ for any $u\in\calU$ and some $\rho > 0$.
Under these structural assumptions on $\barPsi$, from Lemma~\ref{lem:smooth_f_rho}, we see that $h$ is indeed $L_h$-smooth on $\calU$ with $L_h:= L_{uu} + L_{uv}^2/\rho$. Also, the max-structure of $h$ enables us to write the dual problem associated with~\eqref{eq:h_max_struct} as follows:
\begin{equation}
(D):\quad \Xi^*:=  {\max}_{v\in\calV}\;\big\{\Xi(v):= {\inf}_{u\in\calU}\; \barPsi(u,v)+\zeta(u) + \mu\omega_\calU(u)\big\}.
\end{equation}
(Note that from Sion's minimax theorem~\citep{Sion_58}, we know that strong duality holds between $(P)$ and $(D)$, namely $P^*=\Xi^*$.)
Accordingly, let us define the duality gap 
\begin{equation}
\barDelta(u,v):= P(u) - \Xi(v), \quad \forall\,(u,v)\in\calU\times\calV. \label{eq:barDelta}
\end{equation}
The usefulness of the max-structure in~\eqref{eq:h_max_struct} will become clear in Section~\ref{sec:SPP}.

Indeed, in our setting, a typical choice to solve~\eqref{eq:primal_minimization} is the non-Hilbertian proximal gradient methods and its accelerated variants (see e.g.,~\citet{Nest_05,Tseng_08}). These methods assume that the gradient of $h$ at any $u\in\calU$ 
can be obtained {\em exactly}. However, this can be restrictive in some scenarios where computing the gradient involves conducting certain simulations or solving certain optimization problems, which is precisely the case in Section~\ref{sec:SPP}.   In this section, we are interested in the scenario where the value and gradient of $h$ at $u\in\calU$ together has certain {non-zero} but {\em controllable} error $\delta$, and satisfies the $(\delta,\barL)$-inexact model as in~\citet{Dev_14}, which will be reviewed shortly in Section~\ref{sec:inexact_model}.  
Our purpose in this section is to develop an {\em inexact non-Hilbertian} APG method for finding an approximately optimal solution of~\eqref{eq:primal_minimization} under this inexact model. 

Before presenting our method, let us remark that although the inexact APG methods for strongly convex optimization problems  have been well studied in the Hilbertian setting (see e.g.,~\citet{Dev_14,Dev_13,Schmidt_11}), the study in the non-Hilbertian setting has been rather scarce. 
Indeed, when $\bbX$ is Hilbertian, the analyses in the various previous works critically leverage several special properties of the Hilbertian distance  $D_{\omega_\calX}(y,x) = (1/2)\normt{y-x}_\bbX^2$ (where $\omega_\calX = (1/2)\normt{\cdot}_\bbX^2$), including symmetry and inner-product inducibility. Therefore, these analyses cannot be straightforwardly extended to the non-Hilbertian setting, and different techniques have to be developed. Another attractive feature of our method is that convergence guarantees on the duality gap in~\eqref{eq:barDelta} can be obtained when $\calU$ is bounded (in addition to those on the primal optimality gap). To our knowledge, such guarantees have been rarely studied in the literature of inexact proximal gradient methods, even in the Hilbertian setting.

\subsection{$(\delta,\barL)$-inexact model.}\label{sec:inexact_model}
Before presenting our inexact APG method, let us first define the 
$(\delta,\barL)$-inexact model as in~\citet[Definition~1]{Dev_14}. 

\begin{definition}[$(\delta,\barL)$-inexact model] \label{def:delta_L_model}
For any $u\in\calU$, the pair $(\hath(u),\hat{\nabla} h(u))\in\bbR\times\bbU^*$ is called a $(\delta,\barL)$-first-order-approximation (abbreviated as $(\delta,\barL)$-FOA) of $h$ at $u$ if 
\begin{equation}
\hath(u) + \ipt{\hat{\nabla} h(u)}{u'-u}\le h(u')\le \hath(u) + \ipt{\hat{\nabla} h(u)}{u'-u} + (\barL/2)\normt{u'-u}^2 + \delta, \quad \forall\,u'\in\calU. \label{eq:delta_L}
\end{equation}
If we can compute such a $(\delta,\barL)$-FOA of $h$ at any $u\in\calU$, then  we say that $h$ is equipped with the $(\delta,\barL)$-inexact model on $\calU$.  \Halmos
\end{definition}

From the convexity and $L_h$-smoothness of $h$ on $\calU$, given any $u\in\calU$, we see that the exact first-order information $(h(u),\nabla h(u))$ satisfies~\eqref{eq:delta_L} with $\delta=0$ and $\barL=L_h$, and this is the best first-order approximation that we can obtain.  Therefore, we always have $\delta\ge 0$ and $\barL\ge L_h$. Although the $(\delta,\barL)$-inexact model appears a bit unnatural, it well suits the max-structure of $h$ in~\eqref{eq:h_max_struct}. 
Specifically, 
let $\hatv\in\calV$ be an approximate solution of the maximization problem in~\eqref{eq:h_max_struct}. As shown in the following lemma, $(\barPsi(u,\hatv),\nabla_u \barPsi(u,\hatv))$ is indeed 
a $(\delta,\barL)$-FOA of $h$ at $u$. 

\begin{lemma}\label{lem:1st_order_approx}
For any $u\in\calU$ and $\delta>0$, let $\hatv\in\calV$ satisfy that $h(u) - \barPsi(u,\hatv)\le \delta/2$, then we have
\begin{align}
&h(u')\ge \barPsi(u,\hatv) + \ipt{\nabla_u \barPsi(u,\hatv)}{u'-u}, \quad \forall\,u'\in\calU,  \label{eq:delta_L_Psi_lb} \\
&h(u')\le \barPsi(u,\hatv) + \ipt{\nabla_u \barPsi(u,\hatv)}{u'-u} + (\barL/2)\normt{u'-u}^2 + \delta, \quad \forall\,u'\in\calU, \label{eq:delta_L_Psi_ub}
\end{align}
where $\barL = 2L_h$. In words, $(\barPsi(u,\hatv),\nabla_u \barPsi(u,\hatv))$ is a $(\delta,2L_h)$-FOA of $h$ at $u$. 
\end{lemma}

\proof{Proof.}
The proof can be regarded as an extension of that in~\citet[Section~3.2]{Dev_14}, and is shown in Appendix~\ref{app:proof_1st_order_approx}. \Halmos
\endproof

Finally, let us remark that there exist many more scenarios where $h$ is equipped with the $(\delta,\barL)$-inexact model.  For details, we refer readers to~\citet[Section 2.3 and Section 3]{Dev_14}. 




\subsection{Algorithm statement.}\label{sec:subproblem}

\begin{algorithm}[t!]
\caption{Inexact non-Hilbertian APG method} \label{algo:APG}
\begin{algorithmic}
\State {\bf Input}: 
DGF $\omega_\calU$, approximation errors $\{\delta_t\}_{t\ge 0}$, approximate smoothness parameter $\barL\ge L_h$, and weights $\{\alpha_t\}_{t\ge 0}$
\State {\bf Define}: Accumulated weights $\{A_t\}_{t\ge 0}$ where $A_t := \sum_{i=0}^t \alpha_i$ for $t\ge 0$, and averaging sequence $\{\tau_t\}_{t\ge 0}$ where $\tau_t := \alpha_t/A_t$ for $t\ge 0$
\State {\bf Initialize}: $t:=0$, $u_0\in\calU^o$, $s_0 = \hat{\nabla} h(u_0)$ where $(\hath(u_0),\hat{\nabla} h(u_0))$ is a $(\delta_0,\barL)$-FOA of $h$ at $u_0$, and 
$z_0 := \argmin_{u\in\calU}\; \ipt{s_0}{u} + \zeta(u) + \mu\omega_{\calU}(u) + LD_{\omega_\calU}(u,u_0)$, 
\State {\bf Repeat} (until some convergence criterion is satisfied)
\vspace{-.4cm}
\begin{align}
&\baru_{t+1} := {\argmin}_{u\in\calU} \;  \lranglet{s_t}{u} + A_t(\zeta(u)+ \mu\omega_\calU(u))+  \barL D_{\omega_\calU}(u,u_0)\label{eq:baru}\\
&u_{t+1} := (1-\tau_{t+1}) z_t + \tau_{t+1}\baru_{t+1}\label{eq:u}\\
&\mbox{Compute $\hat{\nabla} h(u_{t+1})$ where $(\hath(u_{t+1}),\hat{\nabla} h(u_{t+1}))$ is a $(\delta_{t+1},\barL)$-FOA of $h$ at $u_{t+1}$} \label{eq:approx_FOA} \\
&s_{t+1} := s_t + \alpha_{t+1} \hat{\nabla} h(u_{t+1})\label{eq:s}\\
&w_{t+1} :=  {\argmin}_{u\in\calU} \;  \alpha_{t+1}\big(\lranglet{\hat{\nabla} h(u_{t+1})}{u} + \zeta(u)+ \mu\omega_\calU(u)\big)+  (A_t\mu+\barL) D_{\omega_\calU}(u,\baru_{t+1})\label{eq:w}\\
&z_{t+1} :=  (1-\tau_{t+1}) z_t + \tau_{t+1}w_{t+1}\label{eq:z}\\
&t := t+1
\end{align}
\vspace{-.8cm}
\end{algorithmic}
\end{algorithm}

Our inexact non-Hilbertian APG  method is shown in Algorithm~\ref{algo:APG}. The design of this method leverages Nesterov's famous {\em estimate sequence framework}~\citep{Nest_04}, and in particular the version proposed in~\citet{Nest_05} (which results in an {\em exact} non-Hilbertian APG method for solving~\eqref{eq:primal_minimization} with $\mu=0$).   

Now, let us make some comments on Algorithm~\ref{algo:APG}. First, let us focus on the choices of the input. Intuitively, the approximation errors $\{\delta_t\}_{t\ge 0}$ will accumulate along the iterations of Algorithm~\ref{algo:APG}, in a way that depends on the weights $\{\alpha_t\}_{t\ge 0}$.  
Therefore, the choices of $\{\delta_t\}_{t\ge 0}$ should depend on the accuracy of the approximate solution of~\eqref{eq:primal_minimization} that we wish to find. 
The choice of $\barL$ is made such that $(\hath(u_{t}),\hat{\nabla} h(u_{t}))$ is a $(\delta_{t},\barL)$-FOA of $h$ at $u_{t}$ for all $t\ge 0$, and hence depends on the way that we find the first-order approximation $(\hath(u_{t}),\hat{\nabla} h(u_{t}))$. Lastly, the weights $\{\alpha_t\}_{t\ge 0}$ control both the convergence rate of Algorithm~\ref{algo:APG} as well as the accumulation rate of the approximation errors, and their choices will be made clear in Section~\ref{sec:conv_results_APG}. 

Next, let us focus on the solving the sub-problems in Algorithm~\ref{algo:APG}. Indeed, the sub-problems occur in three places: i) finding $z_0$ in the initialization phase, ii) finding $\baru_{t+1}$ in~\eqref{eq:baru} and iii) finding $w_{t+1}$ in~\eqref{eq:w}. Using the definition of the Bregman divergence $D_{\omega_\calU}(\cdot,\cdot)$ in~\eqref{eq:Bregman_div}, we see that all of the three sub-problems share the same form below:  
\begin{equation}
u^*:= {\argmin}_{u\in\calU}\;\;\zeta(u) + \lranglet{\xi}{u} + \alpha^{-1}\omega_\calU(u) \quad \mbox{for some $\xi\in\bbU^*$ and $\alpha>0$,}  \label{eq:BPP2}
\end{equation}
which is a BPP problem associated with $\zeta$ and the DGF $\omega_\calU$ 
(cf.~\eqref{eq:Bregman_projection} in Section~\ref{sec:BPP}). As assumed at the beginning of  Section~\ref{sec:inexact_APG}, the solution of this problem is easily computable. In the following, let us provide some examples to justify this assumption. 

\underline{\em Some ``easy'' examples of~\eqref{eq:BPP2}.} First, note that if $\bbU$ is Hilbertian with inner product $\lranglet{\cdot}{\cdot}$ and its induced norm $\norm{\cdot}$, and $\calU=\bbU$, we can take $\omega_\calU=(1/2)\norm{\cdot}^2$, so that the problem in~\eqref{eq:BPP2} becomes the usual proximal minimization problem associated with the function $\zeta$. As such, we will provide two simple examples below where $\bbU$ is non-Hilbertian and $\zeta\equiv 0$ (or equivalently, $\zeta=\iota_\calU$ and $\calU = \bbU$). 
For more examples, we refer readers to~\citet{Nest_05,Juditsky_12a}. 

\begin{enumerate}[label={\rm (E\arabic*)}]
\item\label{item:ell_1_entropy} Let $\bbU=(\bbR^n,\norm{\cdot}_1)$ with $\norm{u}_1 \defeq \sum_{i=1}^n \abs{u_i}$, $\calU = \Delta_n$ and $\omega_\calU(u)=\sum_{i=1}^n u_i\ln u_i$. 
From~\citet{Nest_05}, we know that $\omega_\calU$ is 1-strongly convex on $\Delta_n$ with respect to $\norm{\cdot}_1$, and the minimization problem in~\eqref{eq:BPP2}, commonly referred to as the ``entropic projection'' problem, has the following closed-form solution: 
\begin{equation}
u^*_i =  {\exp(-\alpha \xi_i)}\big /\textstyle{\sum_{j=1}^n \exp(-\alpha \xi_j)}, \quad\;\forall\,i\in[n], 
\end{equation}
\item Let $\bbU=(\bbR^n,\norm{\cdot}_p)$, where $p\in(1,2]$ and $\norm{u}_p \defeq \left(\sum_{i=1}^n \abs{u_i}^p\right)^{1/p}$. Consequently, $\bbU^*=(\bbR^n,\norm{\cdot}_q)$, where $q \defeq 1/(1-p^{-1})\in[2,+\infty)$. Let $\calU = \bbR_+^n$ and $\omega_\calU(u) = (1/2)\norm{u}_p^2$, and~\eqref{eq:BPP2} becomes
\begin{equation}
u^*:= {\argmin}_{u\ge 0}\; \lranglet{\xi}{u} + (2\alpha)^{-1}\normt{u}_p^2.\label{eq:BPP_Lp}
\end{equation}
  Note that $\omega_\calU$ is $(p-1)$-strongly convex with respect to $\norm{\cdot}_p$ on $\bbR^n$ (cf.~\citet[Section~8]{BenTal_01}). In addition, let us observe that the minimization problem in~\eqref{eq:BPP_Lp} can be solved in closed-form. Indeed, from the KKT conditions, we easily see that if $\xi_i\ge 0$, then $u^*_i=0$, for all $i\in[n]$. Therefore, without loss of generality, let us assume $\xi<0$, and rewrite~\eqref{eq:BPP_Lp} as
  \begin{equation}
(t^*,\baru^*):= {\argmin}_{t\ge 0}\; {\argmin}_{\baru\ge 0, \normt{\baru}_p=1}\; t\lranglet{\xi}{\baru} + (2\alpha)^{-1}t^2.\label{eq:BPP_Lp2}
\end{equation}
Note that we can recover $u^*$ from $(t^*,\baru^*)$ by letting $u^* = t^*\baru^*$. 
Clearly, since $\xi<0$, we have $\baru_i^* = (\abst{\xi_i}/\normt{\xi}_q)^{q/p}$ for $i\in[n]$ and hence $\lranglet{\xi}{\baru^*} = - \normt{\xi}_q<0$. Based on this, we then have $$t^* = {\argmin}_{t\ge 0} \; (2\alpha)^{-1}t^2  - \normt{\xi}_q t = \alpha\normt{\xi}_q.$$ 
(As a side note, note that the approach above can also be used to derive the closed-form solution of~\eqref{eq:BPP_Lp} without the nonnegativity constraint on $u$.)
\end{enumerate}

Finally, let us observe that all the iterates generated in Algorithm~\ref{algo:APG} (including $\{u_t\}_{t\ge 0}$, $\{\baru_t\}_{t\ge 0}$, $\{z_t\}_{t\ge 0}$ and $\{w_t\}_{t\ge 0}$) lie in $\calU^o$. Indeed,  since $z_0$ is the output of a BPP problem,  we know that $z_0\in\calU^o$ (cf.\ Section~\ref{sec:BPP}), and  the aforementioned observation simply follows  from induction. 

\subsection{Convergence results of Algorithm~\ref{algo:APG}.} \label{sec:conv_results_APG}
Let us present the choices of $\{\alpha_t\}_{t\ge 0}$ and the associated convergence results of Algorithm~\ref{algo:APG} for both the non-strongly-convex ($\mu=0$) and strongly-convex ($\mu>0$) cases.  We will focus on analyzing the strongly-convex case since it is the one that will be used in Section~\ref{sec:SPP}. 
Based on the max-structure of $h$ in~\eqref{eq:h_max_struct}, we let $\{v_t\}_{t\ge 0}\subseteq \calV$ be any sequence that satisfies
\begin{equation}
h(u_t) - \barPsi(u_t,v_t)\le \delta_t/2, \quad\forall\,t\ge 0. \label{eq:v_t_cond}
\end{equation}
Indeed, from Lemma~\ref{lem:1st_order_approx}, we know that $(\hath(u_t),\hat{\nabla} h(u_t)) = (\barPsi(u_t,v_t), \nabla_u\barPsi(u_t,v_t))$ is a $(\delta_t,2L_h)$-FOA of $h$ at $u_t$, for all $t\ge 0$. In addition, define another sequence $\{\barv_t\}_{t\ge 0}\subseteq \calV$ such that
\begin{equation}
\barv_t:= A_t^{-1}\textstyle\sum_{i=0}^t\alpha_i v_i, \quad\forall\,t\ge 0. \label{eq:barv}
\end{equation}
Namely,  $\{\barv_t\}_{t\ge 0}$ is the weighted average of $\{v_t\}_{t\ge 0}$. Our results below not only concern the convergence of the primal optimality gap $\{P(z^t)-P^*\}_{t\ge 0}$, but also the convergence of the duality gap $\{\barDelta(z_t,\barv_t)\}_{t\ge 0}$ (cf.~\eqref{eq:barDelta}),  in the case where $\calU$ is bounded. 

\begin{theorem}[{The case $\mu=0$}]\label{thm:conv_iAPG_mu=0}
In Algorithm~\ref{algo:APG}, if $\mu=0$, then we can choose 
\begin{equation}
\alpha_0 = 1 \quad\mbox{and}\quad  \alpha_t = (2t+3)/4, \quad \forall\, t\ge 1, \label{eq:param_mu=0}
\end{equation}
and under such choices, we have 
\begin{equation}
P(z_t) - P^*\le \frac{4\barL D_{\omega_\calU}(u^*,u_0)}{(t+2)^2} + \frac{\sum_{i=0}^t (i+2)^2 \delta_i}{(t+2)^2}, \quad \forall\,t\ge 0. 
\end{equation}
In addition, if $\calU$ is bounded and  $(\hath(u_t),\hat{\nabla} h(u_t)) = (\barPsi(u_t,v_t), \nabla_u\barPsi(u_t,v_t))$ for all $t\ge 0$, then we can choose $\barL = 2L_h$ and obtain 
\begin{equation}
\barDelta(z_t,\barv_t)\le \frac{8L_h \Omega_{\omega_\calU}(u_0)}{(t+2)^2} + \frac{\sum_{i=0}^t (i+2)^2 \delta_i}{(t+2)^2}, \quad\forall\,t\ge 0, 
\end{equation}
where 
\begin{equation}
\Omega_{\omega_\calU}(u_0):= {\max}_{u\in\calU}\;D_{\omega_\calU}(u,u_0)<+\infty. \label{eq:def_Omega_U}
\end{equation}
\end{theorem}
\proof{Proof.}
See Appendix~\ref{app:proof_inexact_APG}. \Halmos
\endproof

\begin{theorem}[{The case $\mu>0$}]\label{thm:conv_iAPG_mu>0}
In Algorithm~\ref{algo:APG}, if $\mu>0$, then we can 
choose 
\begin{equation}
\alpha_0=1\quad\mbox{and}\quad  \alpha_t = \textstyle(1+\sqrt{\theta})^{t-1}\sqrt{\theta}, \quad \forall\, t\ge 1,
\quad \where \;\;\;\;\theta:= \mu/\barL. \label{eq:param_mu>0}
\end{equation}
Under such choices,  we have 
\begin{align}
P(z_t) - P^*\le (1+\sqrt{\theta})^{-t}{\barL D_{\omega_\calU}(u^*,u_0)} + (1+\sqrt{\theta})^{-t}\textstyle{\sum_{i=0}^t (1+\sqrt{\theta})^{i} \delta_i}, \quad \forall\,t\ge 0.\label{eq:conv_primal_gap}
\end{align}
In addition, if $\calU$ is bounded and  $(\hath(u_t),\hat{\nabla} h(u_t)) = (\barPsi(u_t,v_t), \nabla_u\barPsi(u_t,v_t))$ for all $t\ge 0$, then we can choose $\barL = 2L_h$ and obtain 
\begin{equation}
\barDelta(z_t,\barv_t)\le 2(1+\sqrt{\theta})^{-t}{L_h \Omega_{\omega_\calU}(u_0)} + (1+\sqrt{\theta})^{-t}\textstyle{\sum_{i=0}^t (1+\sqrt{\theta})^{i} \delta_i},  \quad\forall\,t\ge 0,\label{eq:conv_barDelta}
\end{equation}
where $\Omega_{\omega_\calU}(u_0)<+\infty$ is defined in~\eqref{eq:def_Omega_U}. 
\end{theorem}
\proof{Proof.}
See Appendix~\ref{app:proof_inexact_APG}.  \Halmos
\endproof

Let us consider the simple case where the sequence of errors  $\{\delta_t\}_{t\ge 0}$ is uniformly bounded by $\delta>0$ (namely, $\delta_t\le \delta$ for all $t\ge 0$). 
If  $\calU$ is bounded and  $(\hath(u_t),\hat{\nabla} h(u_t)) = (\barPsi(u_t,v_t), \nabla_u\barPsi(u_t,v_t))$ for all $t\ge 0$,  then~\eqref{eq:conv_barDelta} becomes
\begin{align}
\barDelta(z_t,\barv_t)\le {2(1+\sqrt{\theta})^{-t}{L_h \Omega_{\omega_\calU}(u_0)}} + (1+1/\sqrt{\theta})\delta. 
\label{eq:conv_iAPG_cst_error} 
\end{align}
We observe that the right-hand side of~\eqref{eq:conv_iAPG_cst_error} consists of two terms: the first term linearly decreases in $t$ at rate $(1+\sqrt{\theta})^{-1}$, and the second term, which represents the accumulated errors resulted from the approximate gradients $\{\hat{\nabla} h(u_t)\}_{t\ge 0}$, is a constant in $t$ and proportional to $\delta$. 
Consequently, to find the number of iterations $t$ needed to ensure $\barDelta(z_t,\barv_t)\le \epsilon$, we can properly choose $\delta$ such that the second term  $(1+1/\sqrt{\theta})\delta\le \epsilon/2$, and then find $t$ needed such that the first term $$2(1+\sqrt{\theta})^{-t}{L_h \Omega_{\omega_\calU}(u_0)}\le \epsilon/2.$$ 
Of course, if $\calU$ is unbounded, then based on~\eqref{eq:conv_primal_gap}, we can use the same reasoning to find the number of iterations $t$ needed to ensure $P(z_t)-P^*\le \epsilon$. 
This is formalized in the corollary below. 

\begin{corollary}\label{cor:cst_error_iAPG}
In Algorithm~\ref{algo:APG}, if $\mu>0$, choose $\{\alpha_t\}_{t\ge 0}$ as in~\eqref{eq:param_mu>0}.  Fix any $\epsilon>0$ and let $\{\delta_t\}_{t\ge 0}$ satisfy that 
\begin{equation}
\delta_t\le \delta:= \frac{\epsilon}{2(1+\sqrt{\barL/\mu})}, \quad\forall\,t\ge 0. \label{eq:def_delta}
\end{equation} 
Under such choices,  we have that $P(z_t)-P^*\le \epsilon$ for all $t\ge \bart_\rmp$, where 
\begin{equation}
\bart_\rmp\defeq \left\lceil\left(\sqrt{{\frac{\smash[b]{\barL}}{\mu}}}+1\right)\ln\left(\frac{2\barL D_{\omega_\calU}(u^*,u_0)}{\epsilon}\right)\right\rceil. \label{eq:def_bartp}
\end{equation}
In addition,  if $\calU$ is bounded and $(\hath(u_t),\hat{\nabla} h(u_t)) = (\barPsi(u_t,v_t), \nabla_u\barPsi(u_t,v_t))$ for all $t\ge 0$, 
then choose $\barL=2L_h$ and we have that $\barDelta(z_t,\barv_t)\le \epsilon$ for all $t\ge \bart_\rmd$, where 
\begin{align}
\bart_\rmd\defeq \left\lceil\left(\sqrt{{\frac{2L_h}{\mu}}}+1\right)\ln\left(\frac{4L_h \Omega_{\omega_\calU}(u_0)}{\epsilon}\right)\right\rceil. \label{eq:def_bartd}
\end{align}
\end{corollary}




\subsection{An adaptive stopping criterion when $\mu>0$.} \label{sec:stop_crit} 
We derive a sufficient condition to certify the $\epsilon$-optimality of the iterates $\{w_t\}_{t\ge 0}$ generated in Algorithm~\ref{algo:APG}, in the case where $\mu>0$. 
 The purpose of developing this condition is to provide an adaptive stopping criterion that allows us to terminate Algorithm~\ref{algo:APG} ``early'', which we explain below. 

Indeed, from Corollary~\ref{cor:cst_error_iAPG}, we know that if the sequence of approximation errors $\{\delta_t\}_{t\ge 0}$ satisfy~\eqref{eq:def_delta}, then in the {\em worst-case}, we have $P(z^t)-P^*\le \epsilon$ after $\bart_\rmp$ iterations, where $\bart_\rmp$ is defined in~\eqref{eq:def_bartp}. 
However, in certain cases, stopping Algorithm~\ref{algo:APG} after $\bart_\rmp$  iterations can be quite conservative, as some iterates in Algorithm~\ref{algo:APG} (such as $z_t$, $u_t$ or $w_t$)  may already  
be $\epsilon$-optimal for some $t\ll \bart_\rmp$.
Therefore, in this section, we will derive an easy-to-check condition so that as soon as it is satisfied at the $t$-th iteration, we can stop Algorithm~\ref{algo:APG} and conclude that $P(w_{t+1})-P^*\le \epsilon$. In addition, this condition also mitigates the situation where $\calU$ is unbounded and it is difficult to estimate the quantity $D_{\omega_\calU}(u^*,u_0)$ that appears in the definition of $\bart_\rmp$, 
which prohibits us from running Algorithm~\ref{algo:APG} for a fixed number of iterations. 

To state our stopping criterion, given any $u,\baru\in\calU^o$, let us first define
\begin{equation}
u^+ \defeq {\argmin}_{u'\in\calU}\; \lranglet{\nabla h(\baru)+e(\baru)}{u'} + \zeta(u') + \mu\omega_\calU(u') + \lambda^{-1}D_{\omega_\calU}(u',u),\label{eq:u_+}
\end{equation}
where $e(\baru)\in\bbU^*$ denotes the error on the gradient $\nabla h(\baru)$ and $\lambda>0$. Based on $u,\baru,u^+\in\calU^o$, 
define 
\begin{equation}
G:= L_h(u^+-\baru)\qquad\mbox{and}\quad \quad \barG:=  \lambda^{-1}(\nabla \omega_\calX(u)-\nabla \omega_\calX(u^+)).  \label{eq:def_G1_G2}
\end{equation} 
The following lemma is crucial to establish our stopping criterion. 

\begin{lemma}\label{lem:grad_map}
We have 
\begin{equation}
P(u^+)-P^*\le 3\big(\normt{\barG}_*^2 + \normt{G}^2+ \normt{e(\baru)}_*^2\big)/(2\mu).  \label{eq:ub_u+}
\end{equation}
\end{lemma}

\proof{Proof.}
See Appendix~\ref{app:proof_grad_map}. 
 \Halmos
\endproof

Now, let us observe that~\eqref{eq:u_+} has exactly the same form as~\eqref{eq:w}, by letting $\baru := u_{t+1}$, $u := \baru_{t+1}$, $u^+ := w_{t+1}$, $\lambda := \alpha_{t+1}/(A_t\mu+\barL)$ and 
\begin{equation}
e(u_{t+1}) := \hat{\nabla} h(u_{t+1})-\nabla h(u_{t+1}).\label{eq:def_grad_error}
\end{equation}
 Similar to~\eqref{eq:def_G1_G2}, let us define 
\begin{equation}
G_t:= L_h(w_{t+1}- u_{t+1})\quad\mbox{and}\quad \barG_t:= (A_t\mu+\barL)\alpha_{t+1}^{-1}(\nabla \omega_\calU(\baru_{t+1})-\nabla \omega_\calU(w_{t+1})).  \label{eq:def_G1_G2_t}
\end{equation}
Based on Lemma~\ref{lem:grad_map}, we easily have the following stopping criterion. 

\begin{theorem}\label{thm:adapt_stop}
In Algorithm~\ref{algo:APG}, for any $t\ge 0$, if 
\begin{equation}
\normt{\barG_t}_*^2 + \normt{G_t}^2\le \mu\epsilon/3\quad\mbox{and}\quad\normt{e(u_{t+1})}_*\le \sqrt{\mu\epsilon/3}, \label{eq:G_t_crit} 
\end{equation}
then $P(w_{t+1})-P^*\le \epsilon$.
\end{theorem}

Let us make several comments about the stopping criterion in Theorem~\ref{thm:adapt_stop}. First of all, the objects $e(\baru)$, $G$ and $\barG$ (cf.~\eqref{eq:u_+} and~\eqref{eq:def_G1_G2}) together can be regarded as an extension to the {\em proximal gradient mapping} proposed in~\citet{Nest_13}, in the following two senses. First, we allow the error term $e(\baru)$ to exist in the gradient $\nabla h(\baru)$, and second, the projection in~\eqref{eq:u_+} can be non-Hilbertian (namely, we do not restrict $\omega= (1/2)\normt{\cdot}^2$ for some Hilbertian norm $\normt{\cdot}$). Next, let us focus on the implementation of this criterion, which requires  i) computing two additional sequences $\{G_t\}_{t\ge 0}$ and $\{\barG_t\}_{t\ge 0}$ and ii) ensuring that $\normt{e(u_{t+1})}_*\le \sqrt{\mu\epsilon/3}$. For the first requirement, note that all the quantities appearing in the definitions of $G_t$ and $\barG_t$ (cf.~\eqref{eq:def_G1_G2_t}) have already been computed in Algorithm~\ref{algo:APG}, and so computing the additional sequences $\{G_t\}_{t\ge 0}$ and $\{\barG_t\}_{t\ge 0}$ only slightly increases the computational cost of Algorithm~\ref{algo:APG} at each iteration. The second requirement can be accomplished by properly choosing the approximation error $\delta_{t+1}$ and the approximate smoothness parameter $\barL$, which we will discuss below. 

\subsubsection{Estimating $\normt{e(u)}_*$.}\label{sec:est_error}
Let us illustrate two situations where we can ensure $\normt{e(u)}_*\le \epsilon_\rme:= \sqrt{\mu\epsilon/3}$ for some $u\in\calU^o$, where $e(u):= \hat{\nabla} h(u) - \nabla h(u)$ denotes the error on the approximate gradient $\hat{\nabla} h(u)$.  The first situation is more general, and includes any $(\delta,\barL)$-inexact model with ``extended domain''. The second situation is simpler but more restrictive, as it specifically makes use of the max-structure of $h$ as in~\eqref{eq:h_max_struct}. 

\underline{\em Situation I: $(\delta,\barL)$-inexact model with ``extended domain''.}
Let us slightly extend the definition of $(\delta,\barL)$-inexact model in Definition~\ref{def:delta_L_model}, in the sense that~\eqref{eq:delta_L} holds for all $u\in\bar{\calU}$, where $\bar{\calU}$ is a closed convex set with nonempty interior such that $\calU^o\subseteq\inter\bar{\calU}$ --- we shall call this the ``extended'' $(\delta,\barL)$-inexact model. Let $\calB_{\normt{\cdot}}(u,r):= \{u'\in\bbU:\normt{u'-u}\le r\}$ denote the $\normt{\cdot}$-ball centered at $u\in\calU^o$ with radius $r$, and  define the distance from $u$ to $\bdry\bar{\calU}$ (i.e., the boundary of $\bar{\calU}$) as
\begin{equation}
d(u):=\dist(u,\bdry\bar{\calU}) := \sup\big\{r\ge 0:\;\calB_{\normt{\cdot}}(u,r)\subseteq \bar{\calU}\big\}. 
\end{equation}
Using the same argument as in~\citet[Section 2.2]{Dev_14}, if $(\hath(u), \hat{\nabla} h(u))$ satisfies the ``extended'' $(\delta,\barL)$-inexact model, then 
\begin{align}
\normt{e(u)}_*\le 
\begin{cases}
(\barL/2) d(u) + \delta/d(u), \quad &\mbox{if}\quad 0<  d(u) \le \sqrt{2\delta/\barL}\\[1ex]
\sqrt{2\barL\delta}, \quad & \mbox{if}\quad d(u) > \sqrt{2\delta/\barL}
\end{cases}. \label{eq:bound_error} 
\end{align}
Based on~\eqref{eq:bound_error}, it is simple to ensure $\normt{e(u)}_*\le \epsilon_\rme$ by properly choosing $\delta$, as shown below. 

\begin{proposition}
Let $u\in\calU^o$ and $(\hath(u), \hat{\nabla} h(u))$ satisfies the ``extended'' $(\delta,\barL)$-inexact model.  Then $\normt{e(u)}_*\le \epsilon_\rme$ if 
\begin{equation}
\delta := (\epsilon_\rme/2)\min\{ d(u), \epsilon_\rme/\barL \}. \label{eq:def_delta2}
\end{equation}
In particular, if $\bar{\calU} = \bbU$, then $\delta = \epsilon_\rme^2/(2\barL) = \mu\epsilon/(6\barL)$. 
\end{proposition}

\proof{Proof.}
If $d(u)\le \epsilon_\rme/\barL $, then $\delta = \epsilon_\rme d(u)/2 \ge \barL d(u)^2/2$ (or $d(u) \le \sqrt{2\delta/\barL}$), and hence $\normt{e(u)}_*\le (\barL/2) d(u) + \delta/d(u)\le \epsilon_\rme/2 + \epsilon_\rme/2 = \epsilon_\rme.$ If $d(u)> \epsilon_\rme/\barL $, then $\delta  = \epsilon_\rme^2/(2\barL)< \barL d(u)^2/2$ (or $d(u) > \sqrt{2\delta/\barL}$), and hence $\normt{e(u)}_*\le \sqrt{2\barL\delta} = \epsilon_\rme$. If $\calU = \bbU$, then  $d(u)=+\infty$  and $\delta = \epsilon_\rme^2/(2\barL)$. \Halmos   
\endproof

\underline{\em Situation II: $h$ has the max-structure in~\eqref{eq:h_max_struct}.} In this case, let $\hatv\in\calV$ satisfy that $h(u) - \barPsi(u,\hatv)\le \delta/2$. From Lemma~\ref{lem:1st_order_approx}, we already know that $(\hath(u), \hat{\nabla} h(u))=(\barPsi(u,\hatv),\nabla_u \barPsi(u,\hatv))$ is a $(\delta,2L_h)$-FOA of $h$ at $u$.  In fact, the error $e(u) = \nabla_u \barPsi(u,\hatv) - \nabla h(u)$ can also be easily bounded. 

\begin{lemma}\label{lem:norm_bound_err}
If $\hatv\in\calV$ satisfies that $h(u) - \barPsi(u,\hatv)\le \delta/2$, then we have
\begin{equation}
\normt{e(u)}_* = \normt{\hat{\nabla} h(u) - \nabla h(u)}_*=\normt{\nabla_u \barPsi(u,\hatv) - \nabla h(u)}_*\le L_{uv}\sqrt{\delta/\rho}. 
\end{equation}
As a result, we have $\normt{e(u)}_*\le\epsilon_\rme$ if $\delta \le \epsilon_\rme^2\rho/L_{uv}^2$. 
\end{lemma}
\proof{Proof.}
See Appendix~\ref{app:proof_1st_order_approx}. \Halmos
\endproof




\section{An efficient first-order method for solving convex-concave SPPs.}\label{sec:SPP}
Based on the non-Hilbertian inexact APG method developed in Section~\ref{sec:inexact_APG}, we are ready to develop a first-order method for solving a class of (strongly) convex-concave SPPs that subsumes the one in~\eqref{eq:cvx_ccv_SPP} as a special case. We analyze the primal and dual gradient complexities of this method, which enable us to derive the primal and dual gradient complexities of Algorithm~\ref{algo:PD_smoothing} in Section~\ref{sec:choice_lambda}. 

\subsection{Problem setup.}\label{sec:SPP_prob}
Let us consider the following (strongly) convex-concave SPP: 
\begin{equation}
{\min}_{x\in\calX}{\max}_{y\in\calY}\; \big\{S(x,y)\defeq \mu \omega_\calX(x) + r(x) + \Psi(x,y) - g(y) - \rho \omega_\calY(y)\big\},  \label{eq:SPP}
\end{equation}
where $\mu>0$, $\rho> 0$, $\Psi:\calX'\times\calY'\to\bbR$ is jointly continuous on $\calX'\times\calY'$ and convex-concave on $\calX\times\calY$, namely $\Psi(\cdot,y)$ is convex on $\calX$ for any $y\in\calY$ and $\Psi(x,\cdot)$ is concave on $\calY$ for any $x\in\calX$. In addition, $\Psi$ shares the same smoothness assumptions as $\Phi$ in Assumptions~\ref{assump:smooth_Phi} and~\ref{assump:smooth_Phi2},  except that in~\eqref{eq:L_xx}, the smoothness parameter $L_{xx}$ is replaced by a larger one $L'_{xx}\ge L_{xx}$.  

Before proceeding further, let us mention that  the (strongly) convex-concave SPP in~\eqref{eq:cvx_ccv_SPP}, which is solved in Step~\ref{item:inexact_PPM} of Algorithm~\ref{algo:PD_smoothing}, is a special case of the SPP in~\eqref{eq:SPP}. Indeed, in~\eqref{eq:cvx_ccv_SPP},  using the definition of the Bregman divergence $D_{\omega_\calX}(\cdot,\cdot)$ (cf.~\eqref{eq:Bregman_div}), we can write 
\begin{align}
\lambda^{-1}D_{\omega_\calX}(x;x_k) + \Phi(x,y) \eqcst \underbrace{(\lambda^{-1}-\gamma) \omega_\calX(x)}_{:=\mu \omega_\calX(x)} +\underbrace{ \gamma \omega_\calX(x)  - \lambda^{-1} \lranglet{\nabla \omega_\calX(x_k)}{x} + \Phi(x,y)}_{:= \Psi(x,y)}, \label{eq:equiv_form}
\end{align}
where $\eqcst$ omits the terms that are constant w.r.t.\ the optimization variable $x$. Since $\omega_\calX$ is 1-strongly convex on $\calX$ and $\Phi(\cdot,y)$ is $\gamma$-weakly convex on $\calX$, we see that $\Psi(\cdot,y)$ is convex on $\calX$ . In addition, from the $\beta_\calX$-smoothness of $\omega_\calX$ and $L_{xx}$-smoothness of $\Phi(\cdot,y)$ on $\calX$, we have $L'_{xx} := L_{xx} + \gamma\beta_\calX$. Now, substitute~\eqref{eq:equiv_form} into~\eqref{eq:cvx_ccv_SPP} and we see that~\eqref{eq:cvx_ccv_SPP} falls under the problem class in~\eqref{eq:SPP}.

Next, let us write down the primal and dual problems associated with~\eqref{eq:SPP}: 
\begin{align}
\mbox{Primal}: \qquad p^*\defeq\, &{\min}_{x\in\calX} \;\big\{p(x) \defeq \barf(x) + r(x)+\mu\omega_\calX(x)\big\}, 
\label{eq:primal_problem} \\
\mbox{Dual}: \qquad d^*\defeq\, &{\max}_{y\in\calY} \;\big\{d(y) \defeq \pi(y)-g(y)-\rho\omega_\calY(y) \big\},  \label{eq:dual_problem} 
\end{align}
where the functions $\barf:\calX'\to\bbR$ and $\pi:\calY'\to\bbR$ are defined as 
\begin{align}
\barf(x)&\defeq{\max}_{y\in\calY}\; \big\{\psi^\rmD(x,y) \defeq \Psi(x,y) - g(y) - \rho \omega_\calY(y)\big\},\quad \forall\,x\in\calX',\label{eq:bar_f_rho}\\
\pi(y)&\defeq{\min}_{x\in\calX}\; \big\{\psi^\rmP(x,y) \defeq \Psi(x,y)  + r(x)+\mu\omega_\calX(x)\big\}, \quad \forall\, y\in\calY'. \label{eq:def_pi}
\end{align}
We call $p:\calX\to\bbR$ in~\eqref{eq:primal_problem}  and $d:\calY\to\bbR$ in~\eqref{eq:dual_problem} the {\em primal} and {\em dual} functions, respectively. 

Next, we state several facts about the SPP in~\eqref{eq:SPP} and its associated primal and dual problems in~\eqref{eq:primal_problem} and~\eqref{eq:dual_problem}, respectively. These facts will be useful in our algorithmic development. First, note that due to the 1-strong-convexity of $\omega_\calX$ and $\omega_\calY$ on $\calX$ and $\calY$, respectively, the primal function $p$ and dual function $d$ are $\mu$-strongly-convex and $\rho$-strongly-concave on $\calX$ and $\calY$, respectively. As such,  both the primal and dual problems (in~\eqref{eq:primal_problem} and~\eqref{eq:dual_problem}) have unique optimal solutions, which we denote by $x^*\in\calX$ and $y^*\in\calY$, respectively. Since  $S:\calX\times\calY\to\bbR$ in~\eqref{eq:SPP} is convex-concave and jointly continuous on $\calX\times\calY$, together with  the compactness of $\calY$, we can invoke Sion's minimax theorem~\citep{Sion_58} to conclude that $p^*=d^*$. Hence, the SPP in~\eqref{eq:SPP} has a unique saddle point $(x^*,y^*)\in\calX\times\calY$, which by definition satisfies that 
\begin{equation}
S(x^*,y)\le S(x^*,y^*)\le S(x,y^*), \qquad\forall\;(x,y)\in\calX\times\calY. \label{eq:def_sad_pt}
\end{equation} 
Also, notation-wise, let us denote the unique optimal solutions of~\eqref{eq:bar_f_rho} and~\eqref{eq:def_pi} as $y^*(x)$ and $x^*(y)$, respectively, namely,
\begin{equation}
y^*(x)\defeq {\argmax}_{y\in\calY}\; \psi^\rmD(x,y)\quad\mbox{and}\quad  x^*(y)\defeq {\argmin}_{x\in\calX}\;\psi^\rmP(x,y),\quad \forall\,(x,y)\in\calX\times\calY,   \label{eq:primal_dual_opt_sln}
\end{equation}
and from~\eqref{eq:def_sad_pt}, we easily see that 
\begin{equation}
y^*(x^*) = y^* \quad\mbox{and}\quad x^*(y^*) = x^*. \label{eq:y^*(x^*)} 
\end{equation}
In addition, by Lemma~\ref{lem:lips_y*_rho}, we see that the mapping $y^*:\calX'\to\calY$ is $(L_{xy}/\rho)$-Lipschitz on $\calX$ and similarly, the mapping $x^*:\calY'\to\calX$ is $(L_{xy}/\mu)$-Lipschitz on $\calY$. 

Finally, let us show that the function $\pi:\calY'\to\bbR$ in~\eqref{eq:def_pi} is smooth on $\calY$. At this point, it is tempting to conclude this property directly from Lemma~\ref{lem:smooth_f_rho}. However, this approach would require the boundedness of $\calX$, which does not necessarily hold. To circumvent this difficulty, let us 
first note that the function $\pi$ in~\eqref{eq:def_pi} can be equivalently written as 
\begin{align}
\pi(y)={\min}_{x\in\barcalX} \;\psi^\rmP(x,y),\quad \mbox{where}\quad \barcalX:=\clconv x^*(\calY) \label{eq:new_def_pi}
\end{align}
and $x^*(\calY)\defeq \{x^*(y):y\in\calY\}\subseteq\calX$ consists of all the optimal solutions $x^*(y)$ of~\eqref{eq:def_pi} given $y\in\calY$. Since $\calY$ is compact and $x^*$ is continuous on $\calY$, we see that $x^*(\calY)$ is compact, and hence its closed convex hull $\barcalX\subseteq\calX$ in~\eqref{eq:new_def_pi} is convex and compact. Now, based on the new definition of $\pi$ in~\eqref{eq:new_def_pi}, we can invoke Lemma~\ref{lem:smooth_f_rho} to conclude the following.  

\begin{lemma}\label{lem:smooth_pi}
The function $\pi$ is differentiable on $\calY'$ with $\nabla \pi(y) = \nabla_y \Psi(x^*(y),y)$. In addition, $\nabla \pi:\calY'\to\bbY^*$  is $L_\pi$-Lipschitz  on $\calY$, where $L_\pi\defeq L_{yy} + L_{xy}^2/\mu$.
\end{lemma}

Let us define the duality gap associated with~\eqref{eq:SPP} as
\begin{equation}
\Delta(x,y) = p(x) - d(y),\qquad \forall\;(x,y)\in\calX\times\calY.  
\end{equation}
In the following, we propose a non-Hilbertian accelerated dual inexact gradient method for solving~\eqref{eq:SPP}. Specifically, for any $\eta>0$, we aim to find a primal-dual pair $(x,y)\in\calX\times\calY$ such that $$\Delta(x,y)\le \eta.$$ 

\subsection{A dual non-Hilbertian inexact APG method.}\label{sec:dual_inexact}

\begin{figure}[t]\centering
  \includegraphics[width=.8\textwidth]{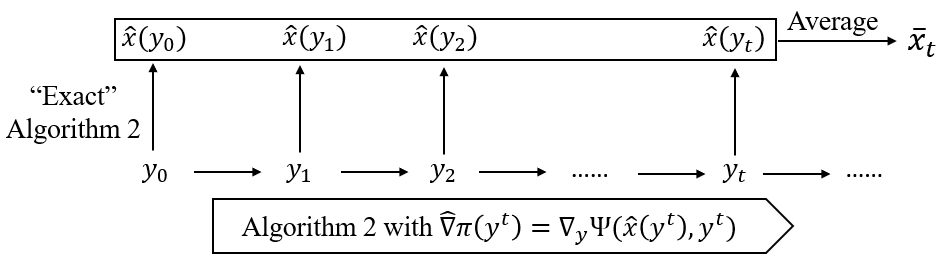}
  \caption{Illustration of the dual non-Hilbertian inexact APG method.}
  \label{fig:illustration_dual_method}
\end{figure}

The idea of this method is conceptually very simple: From Lemma~\ref{lem:smooth_pi}, we know that the function $\pi$ is $L_\pi$-smooth on $\calY$, and hence we can apply the non-Hilbertian accelerated APG method (namely Algorithm~\ref{algo:APG}) as developed in Section~\ref{sec:inexact_APG} to the dual maximization problem in~\eqref{eq:dual_problem}. 
(Note that 
the adaptation of Algorithm~\ref{algo:APG} to maximization problems are straightforward.) 
For each $t\ge 0$, given $y_t\in\calY$, to find the first-order approximation $(\hat{\pi}(y_t),\hat{\nabla} \pi(y_t))$ as in~\eqref{eq:approx_FOA}, we can first find $\hatx(y^t)\in{\calX}$ such that  
\begin{equation}
\psi^\rmP(\hatx(y_t),y_t) - \pi(y_t)\le \bareps/2 \quad \mbox{for some}\;\; \bareps>0, \label{eq:hatx_y_t}
\end{equation}
and then let $(\hat{\pi}(y_t),\hat{\nabla} \pi(y_t)) = (\Psi(\hatx(y_t),y_t), \nabla_y\Psi(\hatx(y_t),y_t))$. 
Indeed, observe that the problem in~\eqref{eq:def_pi} has the same form as the one in~\eqref{eq:primal_minimization}, 
and hence it can then be solved  by the ``exact'' version of Algorithm~\ref{algo:APG} (since the gradient $\nabla_x \Psi(x,y)$ can be computed exactly at any $x\in\calX$). 
In addition, let us define the weighted average of $\{\hatx(y^t)\}_{t\ge 0}$ as 
\begin{equation}
\barx_t:= A_t^{-1}\textstyle\sum_{i=0}^t\alpha_i \hatx(y^i), \quad\forall\,t\ge 0. 
\end{equation}
The structure of this method is illustrated in Figure~\ref{fig:illustration_dual_method}.
From Corollary~\ref{cor:cst_error_iAPG}, we immediately have the following result. 

\begin{corollary}\label{cor:iter_comp_dualAPG}
If we apply  Algorithm~\ref{algo:APG} to solve the dual problem in~\eqref{eq:dual_problem}, by choosing $\{\alpha_t\}_{t\ge 0}$ as in~\eqref{eq:param_mu>0} with $\theta = \rho/(2L_\pi)$, $(\hat{\pi}(y_t),\hat{\nabla} \pi(y_t)) = (\Psi(\hatx(y_t),y_t), \nabla_y\Psi(\hatx(y_t),y_t))$ for all $t\ge 0$ and 
\begin{equation}
\bareps = \frac{\eta}{2(1+\sqrt{2L_\pi/\rho})},  \label{eq:bareps}
\end{equation}
then for any starting point $y_0\in\calY^o$, Algorithm~\ref{algo:APG} generates a sequence $\{\tily_t\}_{t\ge 0}\subseteq \calY^o$ such that $\Delta(\barx_t,\tily_t)\le \eta$ for all $t\ge t_\rmd(\eta)$, where   
\begin{align}
t_\rmd(\eta)\defeq \left\lceil\left(\sqrt{{\frac{2L_\pi}{\rho}}}+1\right)\ln\left(\frac{4L_\pi \Omega_{\omega_\calY}(y_0)}{\eta}\right)\right\rceil 
\label{eq:def_td}
\end{align}
and $\Omega_{\omega_\calY}(y_0): =\max_{y\in\calY}\;D_{\omega_\calY}(y,y_0)<+\infty$. 
\end{corollary}

\begin{remark} \label{rmk:dual_inexact_APG}
Two remarks are in order. First, since each iteration of Algorithm~\ref{algo:APG} only involves computing one dual gradient $\nabla_y\Psi(x,y)$, from Corollary~\ref{cor:iter_comp_dualAPG}, 
we see that to find a primal-dual pair $(x,y)\in\calX\times\calY$ such that $\Delta(x,y)\le \eta$, the number of computed dual gradients 
in the dual inexact APG method is no more than $t_\rmd(\eta)$. 
Second, if the structure of either $\calY$ or $\omega_\calY$ (or both) is relatively simple (e.g., Example~\ref{item:ell_1_entropy}), the quantity  $\Omega_{\omega_\calY}(y_0)$ can be easily estimated. 
In addition, we can stop the method early (i.e., before $t_\rmd(\eta)$ iterations) as soon as  the adaptive stopping criterion as described in Theorem~\ref{thm:adapt_stop} is satisfied. 
\end{remark}

Next, from Corollary~\ref{cor:cst_error_iAPG}, we know that for any $t\ge 0$, to find $\hatx(y_t)\in\calX$ that satisfies~\eqref{eq:hatx_y_t}, the number of iterations of the ``exact'' version of Algorithm~\ref{algo:APG} does not exceed 
\begin{equation}
\left\lceil\left(\sqrt{{\frac{L'_{xx}}{\mu}}}+1\right)\ln\left(\frac{4L'_{xx} D_{\omega_\calX}(x^*(y_t),x_0)}{\bareps}\right)\right\rceil, \label{eq:upper_bound_primal_subprob}
\end{equation}
where $x_0\in\calX^o$ denotes the starting point. Note that in the above, the quantity $D_{\omega_\calX}(x^*(y_t),x_0)$ depends on $y_t\in\calY$, which is inconvenient for our analysis. Therefore, let us upper bound it by 
\begin{equation}
\Gamma_{\omega_\calX}(x_0)\defeq  {\sup}_{y\in\calY}\;\;D_{\omega_\calX}(x^*(y),x_0) = {\sup}_{x\in x^*(\calY)}\;\;D_{\omega_\calX}(x,x^0)<+\infty, \label{eq:upper_bound_BDiam} 
\end{equation}
which is independent of $y\in\calY$. (Note that the finiteness of $\Gamma_{\omega_\calX}(x^0)$ follows from the compactness of $x^*(\calY)\subseteq\calX$ and the continuity of $\omega_\calX$ on $\calX$.) 
By substituting the value of $\bareps$ in~\eqref{eq:bareps} and using $D_{\omega_\calX}(x^*(y_t),x_0)\le \Gamma_{\omega_\calX}(x_0)$, the quantity in~\eqref{eq:upper_bound_primal_subprob} can be upper bounded by 
\begin{equation}
t_\rmp(\eta):= \left\lceil\left(\sqrt{{\frac{L'_{xx}}{\mu}}}+1\right)\ln\left(\frac{8(1+\sqrt{2L_\pi/\rho})L'_{xx}\Gamma_{\omega_\calX}(x^0)}{\eta}\right)\right\rceil. \label{eq:def_tp}
\end{equation}
Consequently, we have the following corollary. 

\begin{corollary}\label{cor:comp_primal_oracle}
Under the setting of Corollary~\ref{cor:iter_comp_dualAPG}, to find a primal-dual pair $(x,y)\in\calX\times\calY$ such that $\Delta(x,y)\le \eta$,  the number of primal gradients $\nabla_x\Psi(x,y)$ computed 
does not exceed $t_\rmd(\eta) t_\rmp(\eta)$, where $t_\rmd(\eta)$ and $t_\rmp(\eta)$ are defined in~\eqref{eq:def_td} and~\eqref{eq:def_tp}, respectively. 
\end{corollary}

\begin{remark}
To find  $\hatx(y_t)\in\calX$ satisfying~\eqref{eq:hatx_y_t}, we can stop the ``exact'' version of Algorithm~\ref{algo:APG} as soon as the adaptive stopping criterion in Theorem~\ref{thm:adapt_stop} is satisfied. Note that in this case, since there are no gradient errors, the stopping criterion in~\eqref{eq:G_t_crit} simplifies to $\normt{\barG_t}_*^2 + \normt{G_t}^2\le \mu\epsilon/3$, and we no longer need to estimate $\normt{e(u)}_*$ as in Section~\ref{sec:est_error}. 
\end{remark}

Finally, let us state the dual and primal gradient complexities of the dual inexact APG method. 
 
\begin{corollary}\label{cor:comp_dual_inexact_APG}
The dual inexact APG method finds a primal-dual pair $(x,y)\in\calX\times\calY$ that satisfies $\Delta(x,y)\le \eta$ with dual gradient complexity
\begin{equation}
\barC_\rmd(\eta) = O\left(\bigg(\sqrt{\frac{L_{yy}}{\rho}}+\frac{L_{xy}}{\sqrt{\mu\rho}}\bigg)\ln\left(\frac{1}{\eta\rho}\right)\right)
\end{equation}
and primal gradient complexity
\begin{equation}
\barC_\rmp(\eta) = O\left(\bigg(\sqrt{\frac{L_{yy}}{\rho}}+\frac{L_{xy}}{\sqrt{\mu\rho}}\bigg)\sqrt{\frac{L'_{xx}}{\mu}}\ln^2\left(\frac{1}{\eta\rho}\right)\right). 
\end{equation}
\end{corollary}

\begin{remark}\label{rmk:comp_dual_iAPG}
As mentioned at the beginning of Section~\ref{sec:SPP}, the SPP in~\eqref{eq:SPP} encompasses that in~\eqref{eq:cvx_ccv_SPP} as special case, with 
\begin{equation}
\mu:=\lambda^{-1}-\gamma \quad \mbox{and}\quad L'_{xx} := L_{xx} + \gamma\beta_\calX. \label{eq:mu_L'_xx}
\end{equation}
Therefore, if we apply the dual inexact APG method to the SPP in~\eqref{eq:cvx_ccv_SPP}, then the dual and primal gradient complexities would be the same as $\barC_\rmd(\eta)$ and $\barC_\rmp(\eta)$ in Corollary~\ref{cor:comp_dual_inexact_APG} above, with $\mu$ and $L'_{xx}$ are replaced by their values in~\eqref{eq:mu_L'_xx}. 
\end{remark}

\section{Extensions and discussions.}\label{sec:extension}

Indeed, our primal-dual smoothing framework in Algorithm~\ref{algo:PD_smoothing} is fairly flexible, and depending on different problem assumptions, it can easily accommodate several variants and extensions. 
Let us discuss two of them in Sections~\ref{sec:simple_dual} and~\ref{sec:stoc}. 
In addition, in Section~\ref{sec:simple_dual}, we also provide a simple variant of the non-Hilbertian inexact APG method in Algorithm~\ref{algo:APG} that only involves solving one BPP problem at each iteration.  

\subsection{``Simple'' dual structure.} \label{sec:simple_dual}

Indeed, in the  cases where the dual maximization problem in the definition of 
$f$ in~\eqref{eq:main} has a simple form,  we may be  able to find a suitable DGF $\omega_\calY$ such that the  dual maximization problem that defines $\barf$ in~\eqref{eq:bar_f_rho} can be easily solved, without appealing to (iterative) first-order methods. For example, in Example~\ref{eg:max_smooth} in Section~\ref{sec:applications}, if we fix any $x\in\calX$ and let $c:= (\ell_i(x))_{i=1}^n$, then the maximization problem in~\eqref{eq:min_smooth_func} can be written as $\max_{p\in\Delta_n}\, c^\top p$. In this case, it is natural to choose $\omega_\calY(p) = \sum_{i=1}^n p_i\ln p_i$ for $p\in\bbR_+^n$, and then the maximization problem in~\eqref{eq:bar_f_rho} becomes ${\max}_{p\in\Delta_n}\; c^\top p - \rho \omega_\calY(p).$ As described in~\ref{item:ell_1_entropy} in Section~\ref{sec:subproblem}, this ``entropic projection'' problem has a simple closed-form solution that can be computed in $O(n)$ time.   

Indeed, from Lemma~\ref{lem:smooth_f_rho}, we know that in the cases above, the gradient of $\barf$ at any $x\in\calX$ can be easily computed, and $\barf$ is $(L'_{xx} + L_{xy}^2/\rho)$-smooth on $\calX$. Therefore, we can directly apply the ``exact'' version of Algorithm~\ref{algo:APG} 
to the primal problem in~\eqref{eq:primal_problem} and obtain an $\eta$-optimal solution $x\in\calX$. 
From Corollary~\ref{cor:cst_error_iAPG}, we know that the primal gradient complexities 
of this scheme is 
\begin{equation}
O\left(\bigg(\sqrt{\frac{L'_{xx}}{\mu}}+\frac{L_{xy}}{\sqrt{\mu\rho}}\bigg)\ln\left(\frac{1}{\eta\rho}\right)\right).\label{eq:primal_g_comp}
\end{equation}
Based on this, we can easily analyze the primal gradient complexity of Algorithm~\ref{algo:PD_smoothing} for finding an $\varepsilon$-NS point of~\eqref{eq:main}. 

\begin{corollary}
For any $\varepsilon>0$, if we choose $\lambda=0.9\gamma^{-1}$ and set $\eta$ as in~\eqref{eq:choice_eta}, then Algorithm~\ref{algo:PD_smoothing} finds an $\varepsilon$-NS point of~\eqref{eq:main} with primal gradient complexity
\begin{equation}
O\left(\big(\sqrt{\gamma L_{xx}}+L_{xy}\gamma\varepsilon^{-1}\big)\varepsilon^{-2}\ln(\varepsilon^{-1}) \right).
\end{equation}
\end{corollary}

Note that due to the ``simple'' dual structure, we neither need to assume any differentiability and smoothness properties of $\Phi(x,\cdot)$ (or equivalently, $\Psi(x,\cdot)$) on $\calY$, 
and nor need to analyze the dual gradient complexity of Algorithm~\ref{algo:PD_smoothing}. 

\subsection{The stochastic setting.} \label{sec:stoc}

Our primal-dual smoothing framework (i.e., Algorithm~\ref{algo:PD_smoothing}) can be easily extend to the stochastic setting, where at any $(x,y)\in\calX\times\calY$, we only have access to the primal gradient $\nabla_x \Phi(x,y)$ and the dual gradient $\nabla_y \Phi(x,y)$ via their unbiased stochastic estimators, denoted by $\tilnabla_x \Phi(x,y)$ and $\tilnabla_y \Phi(x,y)$, respectively. Specifically,  for any $(x,y)\in\calX\times\calY$, we assume that $\tilnabla_x \Phi(x,y)$ and $\tilnabla_y \Phi(x,y)$ satisfy the following conditions: 
\begin{alignat}{2}
&\bbE\big[\tilnabla_x \Phi(x,y)\big] = \nabla_x \Phi(x,y), \qquad &&\bbE\big[\tilnabla_y \Phi(x,y)\big] = \nabla_y \Phi(x,y),\\
&\bbE\big[\normt{\tilnabla_x \Phi(x,y)-\nabla_x \Phi(x,y)}_*^2\big] \le \sigma_x^2, \qquad &&\bbE\big[\normt{\tilnabla_y \Phi(x,y)-\nabla_y \Phi(x,y)}_*^2\big] \le \sigma_y^2,\label{eq:bounded_var}
\end{alignat}
where both $\sigma_x^2,\sigma_y^2<+\infty$. 
Note that the conditions in~\eqref{eq:bounded_var} indicate that both (stochastic) gradient estimators $(x,y)\mapsto\tilnabla_x \Phi(x,y)$ and $(x,y)\mapsto\tilnabla_y \Phi(x,y)$ have uniformly bounded variances over $\calX\times\calY$. 

Similar to the deterministic setting, our goal in the stochastic setting is to find an $\varepsilon$-NS point of~\eqref{eq:main}, but in the sense of {\em expectation}. Specifically,  we aim to find a random point $x\in\calX^o$ such that $\bbE[\normt{x-\prox(q,x,\lambda)}] \le \varepsilon\lambda/\beta_\calX$ (cf.~\eqref{eq:eps_near_stat_pt_normed}). Under this goal, we provide a stochastic extension of Algorithm~\ref{algo:PD_smoothing}, which is shown in Algorithm~\ref{algo:PD_smoothing_stoc}. Compared to Algorithm~\ref{algo:PD_smoothing}, there are two major differences. First, we modify the inexact criterion in Step~\ref{item:inexact_PPM} such that it holds in expectation. Second, instead of terminating Algorithm~\ref{algo:PD_smoothing_stoc} adaptively using the criterion in~\eqref{eq:conv_crit}, we run Algorithm~\ref{algo:PD_smoothing_stoc} for a pre-determined number of iterations $K$. This is because a single realization of $x_k$, in general, do not provide useful information in bounding $\bbE[\normt{x_k-\prox(q,x_k,\lambda)}]$. Consequently, in Algorithm~\ref{algo:PD_smoothing_stoc}, we do not output the last iterate $x_K$ as in Algorithm~\ref{algo:PD_smoothing}, but rather a random point uniformly sampled from $\{x_1,\ldots,x_K\}$.
The convergence guarantee of Algorithm~\ref{algo:PD_smoothing_stoc} is shown in the following theorem. 

\begin{algorithm}[t!]
\caption{Stochastic primal dual smoothing framework 
} \label{algo:PD_smoothing_stoc}
\begin{algorithmic}
\State {\bf Input}: Accuracy parameter $\eta>0$, smoothing parameters $\lambda \in (0,\gamma^{-1})$ and $\rho = \eta/(4R_\calY(\omega_\calY))$
\State {\bf Initialize}: $x_1\in\calX$
\State {\bf For} $k = 1,\ldots,K$: 
\State {\hspace{25pt} Find (a random point) $x_{k+1}\in\calX^o$ such that $\bbE[Q_\rho^\lambda(x_{k+1};x_k)\,|\,x_k]\le q_\rho^\lambda(x_k) + \eta$.}
\State {\bf Output}: $x_{\rm out} = x_k$, where $k\sim{\sf Unif}\{1,\ldots,K\}$ (namely the uniform distribution over $\{1,\ldots,K\}$)
\end{algorithmic}
\end{algorithm}

\begin{theorem}\label{thm:stoc_conv}
Fix any $\varepsilon>0$. 
In Algorithm~\ref{algo:PD_smoothing_stoc}, given any $\Delta_q(x_1)\ge q(x_1)-q^*$ and $\lambda \in (0,\gamma^{-1})$, if we let  
\begin{align}
\eta = \varepsilon^2\lambda/(12\beta_\calX^2) \quad \mbox{and}\quad K = 8\beta_\calX^2\Delta_q(x_1)/(\lambda\varepsilon^2). 
\label{eq:eta_K}
\end{align}
then we have $\bbE[\norm{x_{\rm out} - \prox(q,x_{\rm out},\lambda)}]\le \varepsilon\lambda/\beta_\calX$. 
\end{theorem}

\proof{Proof.}
See Appendix~\ref{app:proof_stoc_conv}. 
 \Halmos
\endproof

\begin{remark}
Note that in Theorem~\ref{thm:stoc_conv}, we need to estimate an upper bound of $q(x_1)-q^*$, i.e., $\Delta_q(x_1)$, which amounts to estimating a lower bound of $q^*$. Note that in many real-words problems, the function $q$, as a loss or cost function, is nonnegative on $\calX$. Therefore, we can simply let $\Delta_q(x_1) = q(x_1)$ in this case. 
\end{remark}

Let us briefly analyze the primal and dual gradient complexities for Algorithm~\ref{algo:PD_smoothing_stoc} to find an $\varepsilon$-NS point of~\eqref{eq:main} in expectation. For simplicity we only focus on the dependence of these complexities on the accuracy parameter $\varepsilon$. Indeed, in each iteration of Algorithm~\ref{algo:PD_smoothing_stoc}, there exist many stochastic first-order methods (e.g.,~\citet[Algorithm~1]{Chen_17} and~\citet[Algorithm~1]{Zhao_19}) that we can use to find the desired $x_{k+1}\in\calX^o$,  
and the primal and dual gradient complexities of these methods all share the same order, i.e.,  $O((\rho\eta)^{-1})$. From the choices of $\rho$ in Algorithm~\ref{algo:PD_smoothing_stoc} and $\eta$ in Theorem~\ref{thm:stoc_conv}, we have $O((\rho\eta)^{-1})=O(\eta^{-2}) = O(\varepsilon^{-4})$. In addition, since the total number of iterations $K=O(\varepsilon^{-2})$ (cf.~Theorem~\ref{thm:stoc_conv}), the primal and dual gradient complexities of Algorithm~\ref{algo:PD_smoothing_stoc} are of order $O(\varepsilon^{-6})$, which  indeed match the state-of-the-art (see e.g.,~\citet{Raf_18}). 

\subsection{A simple variant of Algorithm~\ref{algo:APG}.} \label{eq:simple_APG}

Let us observe that at each iteration in Algorithm~\ref{algo:APG}, we need to solve two BPP problems associated with $\zeta$ and the DGF $\omega_\calU$ in~\eqref{eq:baru} and~\eqref{eq:w}, respectively. One may naturally wonder if it is possible to solve only one BPP problem at each iteration, %
and this leads to Algorithm~\ref{algo:APG2}. Indeed, this algorithm is simpler than Algorithm~\ref{algo:APG}, 
in the sense that it does not involve the sequence $\{w_t\}_{t\ge 0}$, and hence only involves solving one BPP problem at each iteration. The design and analysis of Algorithm~\ref{algo:APG2} are almost identical to those of Algorithm~\ref{algo:APG}, since both algorithms can be  derived from the estimate sequence framework in~\citet{Nest_05}, and they differ  only in one step of the derivation. As such, the convergence guarantees  of Algorithm~\ref{algo:APG2} are similar to those of Algorithm~\ref{algo:APG} in Section~\ref{sec:conv_results_APG}, and we leave the details to the readers. 

\begin{algorithm}[t!]
\caption{A simple variant of Algorithm~\ref{algo:APG}} \label{algo:APG2}
\begin{algorithmic}
\State {\bf Input, Define \& Initialize}: Same as Algorithm~\ref{algo:APG}. Additionally, let $\baru_0 = z_0$.  
\State {\bf Repeat} (until some convergence criterion is satisfied)
\vspace{-.4cm}
\begin{align}
&u_{t+1} := (1-\tau_{t+1}) z_t + \tau_{t+1}\baru_{t}\label{eq:u2}\\
&\mbox{Compute $\hat{\nabla} h(u_{t+1})$ where $(\hath(u_{t+1}),\hat{\nabla} h(u_{t+1}))$ is a $(\delta_{t+1},\barL)$-FOA of $h$ at $u_{t+1}$} \label{eq:approx_FOA2} \\
&s_{t+1} := s_t + \alpha_{t+1} \hat{\nabla} h(u_{t+1})\label{eq:s2}\\
&\baru_{t+1} := {\argmin}_{u\in\calU} \;  \lranglet{s_{t+1}}{u} + A_{t+1}(\zeta(u)+ \mu\omega_\calU(u))+  \barL D_{\omega_\calU}(u,u_0)\label{eq:baru2}\\
&z_{t+1} :=  (1-\tau_{t+1}) z_t + \tau_{t+1}\baru_{t+1}\label{eq:z2}\\
&t := t+1
\end{align}
\vspace{-.8cm}
\end{algorithmic}
\end{algorithm}

\section{Conclusion and future work.}\label{sec:conclusion} 
In this work, we have proposed a primal-dual smoothing framework for finding an $\varepsilon$-NS point of a class of non-smooth non-convex optimization problems in~\eqref{eq:main}. As a contribution of independent interest, we have developed a non-Hilbertian inexact APG method for the strongly convex composite optimization problems in~\eqref{eq:primal_minimization}. 
There are some  problems left open and we wish to consider them in future work. 

First, the lower complexity bound for finding an $\varepsilon$-NS point of~\eqref{eq:main} is not known yet. Establishing this bound with dependence on the problem parameters (including $L_{xx}$, $L_{xy}$, $L_{yy}$ and $\gamma$) and the accuracy $\varepsilon$ can be useful to understand the ``optimality'' of the existing methods (including ours).  


Second, as detailed in Section~\ref{sec:stoc}, a straightforward extension of our framework to the stochastic setting can find an $\varepsilon$-NS point of~\eqref{eq:main} in expectation with primal and dual gradient complexities both of order $O(\varepsilon^{-6})$. It seems that this result can be further improved using ``cleverer'' 
strategies, and developing these strategies would be an interesting direction for future research. 

Third, in Section~\ref{sec:BPP}, the additional assumptions~\ref{item:X_domain} and~\ref{item:grad_lips} that we place on the DGF $\omega_\calX$ 
appear to be somewhat stringent. In fact, the well-known example  in~\ref{item:ell_1_entropy} does not satisfy either assumption. Note that assumption~\ref{item:grad_lips} amounts to assuming that $\normt{\nabla^2 \omega_\calX(x)}$ (namely, the operator norm of $\nabla^2 \omega_\calX(x)$) is uniformly bounded by $\beta_\calX<+\infty$ over $x\in\calX$. Without assumption~\ref{item:X_domain}, this may fail even when $\calX$ is bounded, which is precisely due to the potential ``blow-up'' behavior of  $\nabla^2 \omega_\calX(\cdot)$ near the  boundary of $\dom\omega_\calX$. The failure of assumption~\ref{item:grad_lips} poses serious challenges in designing optimization algorithms for finding an $\varepsilon$-NS point of~\eqref{eq:main} with complexity guarantees. We believe that addressing this problem will have far-reaching impact 
in the broader context of non-Euclidean non-convex optimization, and it is very worthwhile to pursue this problem in the future. 

\section*{Acknowledgments.}
The author would like to thank the two anonymous referees for their many constructive suggestions that have significantly improved the exposition of the current manuscript. The author would also like to thank Yangyang Xu for inspirational discussions, and Robert M.\ Freund for his constructive feedback during the preparation and revision of this manuscript. The author's research is supported by AFOSR Grant No.\ FA9550-22-1-0356.

\begin{APPENDICES}

\section{Proof of Lemma~\ref{lem:Danskin}.}\label{app:proof_Danskin}


Fix any $x\in\calX$ and any $d\in\bbX$.  
Consider any sequences $\{t_n\}_{n\ge 0}\subseteq\bbR$ and $\{d_n\}_{n\ge 0}\subseteq\bbX$ such that $t_n\downarrow 0$ and $d_n\to d$. Since $\calX\subseteq\calX'$ and $\calX'$ is open, then there exists some $N \ge 0$ such that for all $n\ge N$, $x+t_nd_n\in\calX'$. Note that due to the compactness of $\calY$, for all $x\in\calX'$, the set of dual optimal solutions $\calY^*(x)$ (as defined in~\eqref{eq:dir_deriv_f}) is nonempty.  
By the definition of $f$ in~\eqref{eq:main}, 
for any $y\in\calY^*(x)$ and $n\ge N$, we have
\begin{equation}
\frac{f(x+t_nd_n)-f(x)}{t_n}\ge \frac{\Phi(x+t_nd_n,y)-\Phi(x,y)}{t_n}. 
\label{eq:finite_diff_ge}
\end{equation}
Therefore, we have for any $y\in\calY^*(x)$, 
\begin{align*}
\liminf_{n\to+\infty}\frac{f(x+t_nd_n)-f(x)}{t_n}&\ge \liminf_{n\to+\infty}\frac{\Phi(x+t_nd_n,y)-\Phi(x,y)}{t_n}\\
&=  \liminf_{n\to+\infty}\frac{t_n\lranglet{\nabla_x \Phi(x,y)}{d_n} + o(t_n\norm{d_n})}{t_n}= \lranglet{\nabla_x \Phi(x,y)}{d}.
\end{align*}
As such, we have 
\begin{equation}
\liminf_{n\to+\infty}\frac{f(x+t_nd_n)-f(x)}{t_n}\ge {\sup}_{y\in\calY^*(x)}\;\lranglet{\nabla_x \Phi(x,y)}{d}. \label{eq:direc_deriv_ge}
\end{equation}

Next, let $\{x_n\}_{n\ge 0}\subseteq \calX'$ be any sequence such that $x_n\to x$. We aim to show that if $y_n\in \calY^*(x_n)$ for all $n\ge 0$ and $y_n\to y\in\calY$, then $y\in \calY^*(x).$
Indeed, by definition, 
\begin{align*}
\limsup_{n\to+\infty} f(x_n) &= \limsup_{n\to+\infty} \Phi(x_n,y_n) - g(y_n)\eqa \Phi(x,y) - g(y)\leb f(x),\nt\label{eq:y_inc_a}
\end{align*}
where (a) follows from the joint continuity of $\Phi$ on $\calX'\times\calY$ (cf.\ Section~\ref{sec:Problem}) 
and the continuity of $g$ on $\calY$,  and (b) follows from the definition of $f$ in~\eqref{eq:main} and that $y\in\calY$. 
On the other hand, since $\Phi(\cdot,y)$ is continuous on $\calX$ for any $y\in\calY$, $f$ is clearly lower semicontinuous. As a result, we have 
\begin{equation}
\limsup_{n\to+\infty} f(x_n)\ge \liminf_{n\to+\infty}  f(x_n) \ge f(x).\label{eq:y_inc_b}
\end{equation}
Combining~\eqref{eq:y_inc_a} and~\eqref{eq:y_inc_b}, we have $f(x)=\Phi(x,y)- g(y)$, implying that $y\in\calY^*(x)$. 
As a result, if we let $x_n \defeq x+t_nd_n$ (so that $x_n\to x$) and $y_n\in \calY^*(x+t_nd_n)$, then any limit point of $\{y_n\}_{n\in\bbN}$ (which exists since $\calY$ is compact) belongs to $\calY^*(x)$. 
Consequently, we have 
\begin{align*}
\limsup_{n\to+\infty}\frac{f(x+t_nd_n)-f(x)}{t_n}
&\le \limsup_{n\to+\infty}\frac{\Phi(x+t_nd_n,y_n)-\Phi(x,y_n)}{t_n}\\
&= \limsup_{n\to+\infty}\frac{t_n\lranglet{\nabla_x \Phi(x,y_n)}{d_n} + o(t_n\norm{d_n})}{t_n}\\
& \le {\sup}_{y\in\calY^*(x)}\;\lranglet{\nabla_x \Phi(x,y)}{d}.\nt\label{eq:direc_deriv_le}
\end{align*}
Combining~\eqref{eq:direc_deriv_ge} and~\eqref{eq:direc_deriv_le}, we see that for all $x\in\calX'$ and $d\in\bbX$, $f'(x;d)$ exists and 
\begin{equation}
f'(x;d)=\lim_{n\to+\infty}\frac{f(x+t_nd_n)-f(x)}{t_n}={\sup}_{y\in\calY^*(x)}\;\lranglet{\nabla_x \Phi(x,y)}{d}.\tag*{\Halmos}
\end{equation}


\begin{remark}
Our proof of Lemma~\ref{lem:Danskin} can be regarded as a simplified version of that for~\citet[Theorem~D1]{Bernhard_95}. This is because we assume the Fr\'echet differentiablility of $\Phi(\cdot,y)$, which is stronger 
than the notion of G\^ateaux differentiablility of $\Phi(\cdot,y)$ assumed in~\citep[Theorem~D1]{Bernhard_95}. However, note that our result is also stronger, namely we show that $f$ is (Hadamard) directionally differentiable, whereas~\citep[Theorem~D1]{Bernhard_95} only shows that $f$ is G\^ateaux directionally  differentiable, a notion weaker than (Hadamard) directional differentiablility. 
\end{remark}

\section{Proof of Lemma~\ref{lem:Lipschitz} } \label{app:proof_local_Lips}

Fix any  $x\in\calX$ and 
consider its compact neighborhood $\calV(x)$ in $\calX$, namely $x\in\calV(x)\subseteq\calX$. 
Define $M_{\calV(x)}\defeq {\sup}_{(z,y)\in\calV(x)\times\calY}\; \normt{\nabla_x\Phi(z,y)}_*$, and note that $M_{\calV(x)}<+\infty$ since $\nabla_x\Phi(\cdot,\cdot)$ is jointly (Lipschitz) continuous on $\calX\times\calY$ (which follows from Assumption~\ref{assump:smooth_Phi}) and $\calV(x)$ is compact. 
For any $x',x''\in\calV(x)$, we have
\begin{align*}
\abst{f(x') - f(x'')}&=\abs{{\sup}_{y\in\calY}\;[\Phi(x',y) - g(y)] - {\sup}_{y\in\calY}\; [\Phi(x'',y) - g(y)]} \\
&\le {\sup}_{y\in\calY}\abs{\Phi(x',y) - \Phi(x'',y)}\\
&\le {\sup}_{y\in\calY}\int_{0}^1 \abs{\lranglet{\nabla_x\Phi(x''+t(x'-x''),y)}{x'-x''}}\; \rmd t\\
&\le {\sup}_{y\in\calY}\int_{0}^1 \normt{\nabla_x\Phi(x''+t(x'-x''),y)}_*\;\rmd t\;\; \normt{x'-x''}\\
&\le M_{\calV(x)}\normt{x'-x''}.
\nt \label{eq:local_Lips}
\end{align*}
This shows that $f$ is $M_{\calV(x)}$-Lipschitz on $\calV(x)$. 

Next, let us show $ \partial f(x)=\clconv\{\nabla_x \Phi(x,y):y\in\calY^*(x)\}$. For notational convenience, let 
\begin{equation}
\calA(x): = \{\nabla_x \Phi(x,y):y\in\calY^*(x)\}. \label{eq:calA}
\end{equation}
Fix any $x\in\calX'$ and any $d\in\bbX$. 
We first show that 
$\clconv\calA(x)\subseteq \partial  f(x)$. 
To see this, if $y\in\calY^*(x)$, then 
\begin{align*}
f(x+d) - f(x) - \lranglet{\nabla_x \Phi(x,y)}{d}\ge \Phi(x+d,y) - \Phi(x,y) -  \lranglet{\nabla_x \Phi(x,y)}{d} = o(\normt{d}),
\end{align*}
and hence $\nabla_x \Phi(x,y)\in\partial f(x)$. Since $\partial f(x)$ is closed and convex, we have $\clconv\calA(x)\subseteq \partial  f(x)$. 
Then, by~\eqref{eq:dir_deriv_f}, we have
\begin{align}
f'(x;d) = {\sup}_{\xi\in\calA(x)}\;\lranglet{\xi}{d} = {\sup}_{\xi\in\clconv\calA(x)}\;\lranglet{\xi}{d}\le {\sup}_{\xi\in\partial f(x)}\;\lranglet{\xi}{d}. \label{eq:direc_deriv_2}
\end{align}
On the other hand, for any $\xi\in\partial f(x)$, 
by the definition in~\eqref{eq:dir_deriv_f}, we have  
\begin{align*}
f'(x;d)=\lim_{t\downarrow 0,\;d'\to d} \frac{f(x+td')-f(x)}{t} & \ge \lim_{t\downarrow 0,\;d'\to d} \frac{t\lranglet{\xi}{d'} + o(t\normt{d'})}{t} = \lranglet{\xi}{d}. 
\end{align*}
This shows that $f'(x;d)\ge  \sup_{\xi\in\partial f(x)} \lranglet{\xi}{d}$. This, together with~\eqref{eq:direc_deriv_2}, implies that 
\begin{equation}
f'(x;d)=  {\sup}_{\xi\in\partial f(x)}\; \lranglet{\xi}{d}={\sup}_{\xi\in\clconv\calA(x)}\;\lranglet{\xi}{d}. 
\end{equation} 
Since both $\partial f(x)$ and $\clconv\calA(x)$ are closed and convex sets and share  the same support function  $d\mapsto f'(x;d)$, they share the same indicator function (cf.~\citet[Theorem~13.2]{Rock_70}) and hence $ \partial f(x)=\clconv\calA(x)$. \Halmos

\section{Proof of Lemma~\ref{lem:lips_y*_rho}.}\label{app:proof_lips_y*_rho}

Consider any $x,x'\in\calX$. By the $\rho$-strong concavity of $\phi^\rmD_\rho(x,\cdot)$ on $\calY$ (for any $x\in\calX'$), we have
\begin{align*}
(\rho/2) \normt{y_\rho^*(x') - y_\rho^*(x)}^2&\le \phi^\rmD_\rho(x, y_\rho^*(x))-\phi^\rmD_\rho(x,y_\rho^*(x')),\\
(\rho/2) \normt{y_\rho^*(x') - y_\rho^*(x)}^2&\le \phi^\rmD_\rho(x', y_\rho^*(x'))-\phi^\rmD_\rho(x',y_\rho^*(x)). 
\end{align*}
As a result, 
\begin{align*}
\normt{y_\rho^*(x') - y_\rho^*(x)}^2 &\le \rho^{-1}\big(\phi^\rmD_\rho(x, y_\rho^*(x))-\phi^\rmD_\rho(x,y_\rho^*(x'))+\phi^\rmD_\rho(x', y_\rho^*(x'))-\phi^\rmD_\rho(x', y_\rho^*(x))\big)\\
&\le \rho^{-1}\int_{0}^1\lrangle{\nabla_x \phi^\rmD_\rho(x'+t(x-x'), y_\rho^*(x))-\nabla_x \phi^\rmD_\rho(x'+t(x-x'), y_\rho^*(x'))}{x-x'}\rmd t\\
&\le \rho^{-1}\int_{0}^1 \normt{\nabla_x \phi^\rmD_\rho(x'+t(x-x'), y_\rho^*(x))-\nabla_x \phi^\rmD_\rho(x'+t(x-x'), y_\rho^*(x'))}_*\normt{x-x'}\rmd t\\
&\le (L_{xy}/\rho)\normt{y_\rho^*(x)-y_\rho^*(x')}\normt{x-x'}. \nt\label{eq:y_rho_Lips}
\end{align*}
If $y_\rho^*(x)=y_\rho^*(x')$, then we trivially have $\normt{y_\rho^*(x)-y_\rho^*(x')}\le (L_{xy}/\rho)\normt{x-x'}$. Otherwise, this can be obtained by dividing both sides of~\eqref{eq:y_rho_Lips} by $\normt{y_\rho^*(x)-y_\rho^*(x')}$.

\section{Proof of Lemma~\ref{lem:smooth_f_rho}.}\label{app:proof_sm_f_rho}

 The differentiability of $f_\rho$ on $\calX'$ directly follows from the $\rho$-strong concavity of $\phi^\rmD_\rho(x,\cdot)$ on $\calY$ (for any $x\in\calX'$) and Lemma~\ref{lem:Danskin}, from which we also see that $\nabla f_\rho (x) = \nabla_x \Phi(x,y^*_\rho(x))$. Consequently, we have 
\begin{align*}
\normt{\nabla f_\rho(x')-\nabla f_\rho(x)}_* &= \normt{\nabla_x \phi^\rmD_\rho(x',y^*(x'))-\nabla_x \phi^\rmD_\rho(x,y^*(x))}_*\\
&\le \normt{\nabla_x \phi^\rmD_\rho(x',y^*(x'))-\nabla_x \phi^\rmD_\rho(x,y^*(x'))}_* + \normt{\nabla_x \phi^\rmD_\rho(x,y^*(x'))-\nabla_x \phi^\rmD_\rho(x,y^*(x))}_*\\
&\le L_{xx}\normt{x-x'} + L_{xy} \normt{y^*(x)-y^*(x')}\\
&\le (L_{xx}+ L^2_{xy}/\rho) \normt{x-x'}, 
\end{align*}
where the last step follows from Lemma~\ref{lem:lips_y*_rho}. 

\section{Proof of Lemma~\ref{lem:dist_grad_subdiff}.} \label{app:proof_dist_grad_subdiff}

From Lemmas~\ref{lem:Lipschitz} and~\ref{lem:smooth_f_rho}, we know that $\calA(x)\subseteq \partial f(x)$ (where $\calA(x)$ is defined in~\eqref{eq:calA}) and $\nabla f_\rho (x) = \nabla_x \Phi(x,y^*_\rho(x))$, and hence 
\begin{align*}
\dist(\nabla f_\rho(x),\partial f(x))&\le \dist\big(\nabla_x \Phi(x,y^*_\rho(x)),\calA(x)\big) \nn\\
&={\inf}_{y\in\calY^*(x)}\; \normt{\nabla_x \Phi(x,y^*_\rho(x)) - \nabla_x \Phi(x,y)}_*\nn\\
&\lea L_{xy} {\inf}_{y\in\calY^*(x)}\;\normt{y^*_\rho(x)-y}\\
&= L_{xy}\dist(y^*_\rho(x),\calY^*(x)), 
\end{align*}
where (a) follows from Assumption~\ref{assump:smooth_Phi}. This shows the first inequality in~\eqref{eq:bound_dist_smoothed_grad}. Also, as both $y^*_\rho(x)\in\calY$ and $\calY^*(x)\subseteq\calY$, we have $\dist(y^*_\rho(x),\calY^*(x))\le D_\calY$, which shows the second inequality.

\section{Proof of Lemmas~\ref{lem:1st_order_approx} and~\ref{lem:norm_bound_err}.}\label{app:proof_1st_order_approx}
Given any $u\in\calU$, let $v^*(u)\in\calV$ be the optimal solution of the maximization problem in~\eqref{eq:h_max_struct}, and hence $h(u) = \barPsi(u,v^*(u))$. From the $\rho$-strong concavity of $\barPsi(u,\cdot)$, we have 
\begin{equation}
(\rho/2)\normt{v^*(u) - \hatv}^2 \le h(u) - \barPsi(u,\hatv)\le \delta/2 \quad\Longrightarrow\quad \normt{v^*(u) - \hatv}^2\le \delta/\rho. \label{eq:dist_v^*u}
\end{equation}
 Therefore, we have
\begin{align*}
h(u') = \barPsi(u',v^*(u')) \gea \barPsi(u',\hatv)\geb \barPsi(u,\hatv) + \ipt{\nabla_u \barPsi(u,\hatv)}{u'-u},
\end{align*}
where (a) follows from $\hatv\in\calV$ and (b) follows from the convexity of $\barPsi(\cdot,\hatv)$. 
This shows~\eqref{eq:delta_L_Psi_lb}. On the other hand, from the $L_h$-smoothness of $h$ on $\calU$, we have 
\begin{align}
h(u')&\le h(u) + \ipt{\nabla h(u)}{u'-u} + (L_h/2)\normt{u'-u}^2\nn\\ 
&\le \barPsi(u,\hatv) + \delta/2 + \ipt{\nabla_u \barPsi(u,\hatv)}{u'-u} + \ipt{\nabla h(u)-\nabla_u \barPsi(u,\hatv)}{u'-u}+ (L_h/2)\normt{u'-u}^2.\label{eq:bound_sm_h}
\end{align}
In addition, from Lemma~\ref{lem:smooth_f_rho}, we have $\nabla h(u) = \nabla_u \barPsi(u,v^*(u))$, and hence
\begin{align}
\normt{\nabla h(u)-\nabla_u \barPsi(u,\hatv)}_*^2=\normt{\nabla_u \barPsi(u,v^*(u))-\nabla_u \barPsi(u,\hatv)}_*^2\lea L_{uv}^2\normt{v^*(u)-\hatv}^2\leb (L_{uv}^2/\rho)\delta,\label{eq:norm_bound_err}
\end{align}
where (a) follows from~\eqref{eq:L_uv} and (b) follows from~\eqref{eq:dist_v^*u}. This proves Lemma~\ref{lem:norm_bound_err}.  Consequently, 
\begin{align*}
\ipt{\nabla h(u)-\nabla_u \barPsi(u,\hatv)}{u'-u}&\le ({2L_h})^{-1}\normt{\nabla h(u)-\nabla_u \barPsi(u,\hatv)}_*^2+({L_h}/{2})\normt{u'-u}^2\\
&\lea (\delta/2) (L_{uv}^2/\rho)/L_h + ({L_h}/{2})\normt{u'-u}^2\\
&\leb  \delta/2 + ({L_h}/{2})\normt{u'-u}^2, \nt\label{eq:bound_ip}
\end{align*}
where (a) follows from~\eqref{eq:norm_bound_err} and (b) follows  from $L_h = L_{uu} + L_{uv}^2/\rho\ge L_{uv}^2/\rho$. Now, by substituting~\eqref{eq:bound_ip} into~\eqref{eq:bound_sm_h}, we have~\eqref{eq:delta_L_Psi_lb}.

\section{Proof of Theorems~\ref{thm:conv_iAPG_mu=0} and~\ref{thm:conv_iAPG_mu>0}.}  \label{app:proof_inexact_APG}
Let us begin our proof by defining the sequence of functions $\{\psi_t:\calU\to\bbR\}_{t\ge 0}$ such that for all $t\ge 0$,  
\begin{equation}
\psi_t(u):= \textstyle\sum_{i=0}^t \;\alpha_i(\hath(u_i) + \ipt{\hat{\nabla} h(u_i)}{u-u_i} + \zeta(u)+ \mu\omega_\calU(u)) + \barL D_{\omega_\calU}(u,u_0), \quad\forall\,u\in\calU. \label{eq:def_psi_t}
\end{equation}
Indeed, the functions $\{\psi_t\}_{t\ge 0}$ play pivotal roles in analyzing Algorithm~\ref{algo:APG}. In addition, let us define 
\begin{equation}
\psi_t^*:= {\min}_{u\in\calU}\; \psi_t(u), \quad \forall\,t\ge 0. 
\end{equation}
Our proof can be streamlined into the following three lemmas. 
The first lemma below establishes a lower bound of $\psi^*_t$, for all $t\ge 0$. 

\begin{lemma} \label{lem:psi_t^*_lb}
If $\alpha_0 = 1$ and 
\begin{equation}
(A_t\mu+\barL)A_{t+1}\ge \barL \alpha_{t+1}^2,\quad \forall\, t\ge 0,  \label{eq:alpha_A_cond}
\end{equation}
then we have that
for all $t\ge 0$, 
\begin{equation}
\psi^*_t\ge A_tP(z_t) - E_t, \quad \where\quad  E_t:= \textstyle \sum_{i=0}^t \, A_i\delta_i. \label{eq:psi_t_lb} 
\end{equation}
\end{lemma}

\proof{Proof.}
Let us show~\eqref{eq:psi_t_lb} using induction. When $t=0$, from the initialization in Algorithm~\ref{algo:APG}, we see that $z_0:= \argmin_{u\in\calU}\,\psi_0(u)$ and hence
\begin{align*}
\psi^*_0 = \psi_0(z_0) &= \alpha_0(\hath(z_0) + \ipt{\hat{\nabla} h(z_0)}{u-z_0} + \zeta(z_0)+ \mu\omega_\calU(z_0)) + \barL D_{\omega_\calU}(u,z_0)\\
&\gea \hath(z_0) + \ipt{\hat{\nabla} h(z_0)}{u-z_0} + (\barL/2) \normt{u-z_0}^2 + \zeta(z_0)+ \mu\omega_\calU(z_0)\\
&\geb h(z_0) - \delta_0 + \zeta(z_0)+ \mu\omega_\calU(z_0) \eqc P(z_0) - A_0\delta_0, 
\end{align*}
where in (a) we use $\alpha_0=1$ and the 1-strong convexity of $\omega_\calU$ on $\calU$, in (b) we use the second inequality in~\eqref{eq:delta_L} and in (c) we use $A_0 = \alpha_0=1$. 
Now, suppose that~\eqref{eq:psi_t_lb} holds for some $t\ge 0$, and let us show that~\eqref{eq:psi_t_lb} holds for $t+1$. First, observe that for any $u\in\calU$, we have 
\begin{align*}
\psi_{t+1}(u) &= \psi_t(u) + \alpha_{t+1}(\hath(u_{t+1}) + \ipt{\hat{\nabla} h(u_{t+1})}{u-u_{t+1}} + \zeta(u)+ \mu\omega_\calU(u))\\
&\gea \psi_t^* + (A_t\mu+\barL)D_{\omega_\calU}(u,\baru_{t+1}) + \alpha_{t+1}(\hath(u_{t+1}) + \ipt{\hat{\nabla} h(u_{t+1})}{u-u_{t+1}} + \zeta(u)+ \mu\omega_\calU(u))\\
&\geb \psi_t^* + (A_t\mu+\barL)D_{\omega_\calU}(w_{t+1},\baru_{t+1})\\
&\qquad\qquad + \alpha_{t+1}(\hath(u_{t+1}) + \ipt{\hat{\nabla} h(u_{t+1})}{w_{t+1}-u_{t+1}} + \zeta(w_{t+1})+ \mu\omega_\calU(w_{t+1})),\nt \label{eq:psi_t+1_lb}
\end{align*}
where (a) follows from $\baru_{t+1} = \argmin_{u\in\calU}\, \psi_{t}(u)$, since $s_t = \textstyle\sum_{i=0}^t \alpha_i \hat{\nabla} h(u_t)$  (cf.~\eqref{eq:baru}) and (b) follows from the definition of $w_{t+1}$ in~\eqref{eq:w}.  
Using the induction hypothesis, we have
\begin{align}
\psi_t^*&\ge A_t(h(z_t) + \zeta(z_t) + \mu\omega_{\calU}(z_t)) - E_t\nn\\
& \gea A_t(\hath(u_{t+1}) + \ipt{\hat{\nabla} h(u_{t+1})}{z_t-u_{t+1}} + \zeta(z_t) + \mu\omega_{\calU}(z_t))- E_t, \label{eq:psi_t*_lb}
\end{align} 
where (a) follows from the first inequality in~\eqref{eq:delta_L}. 
Combining~\eqref{eq:psi_t+1_lb} and~\eqref{eq:psi_t*_lb}, we have 
\allowdisplaybreaks
\begin{align*}
\psi_{t+1}(u)&\ge A_{t+1}\hath(u_{t+1}) + \ipt{\hat{\nabla} h(u_{t+1})}{A_tz_t+\alpha_{t+1} w_{t+1} -A_{t+1}u_{t+1}}\\
&\qquad  + \big\{A_t(\zeta(z_t) + \mu\omega_{\calU}(z_t)) + \alpha_{t+1}(\zeta(w_{t+1})+ \mu\omega_\calU(w_{t+1}))\big\}  + (A_t\mu+\barL)D_{\omega_\calU}(w_{t+1},\baru_{t+1}) - E_t\\
&\gea  A_{t+1} \big\{ \hath(u_{t+1}) + \ipt{\hat{\nabla} h(u_{t+1})}{(1-\tau_{t+1})z_t+\tau_{t+1} w_{t+1} -u_{t+1}}\\
&\qquad + (1-\tau_{t+1}) (\zeta(z_t) + \mu\omega_{\calU}(z_t)) + \tau_{t+1}(\zeta(w_{t+1})+ \mu\omega_\calU(w_{t+1}))\\
&\qquad  + ((A_t\mu+\barL)/(2A_{t+1}))\normt{w_{t+1}-\baru_{t+1}}^2\big\} - E_t\\
&\geb A_{t+1} \big\{ \hath(u_{t+1}) + \ipt{\hat{\nabla} h(u_{t+1})}{z_{t+1} -u_{t+1}} + \zeta(z_{t+1}) + \mu\omega_\calU(z_{t+1})\\
&\qquad +((A_t\mu+\barL)/(2A_{t+1}\tau_{t+1}^2))\normt{z_{t+1}-u_{t+1}}^2\big\} - E_t\\
&\gec A_{t+1} \big\{ \hath(u_{t+1}) + \ipt{\hat{\nabla} h(u_{t+1})}{z_{t+1} -u_{t+1}}+(\barL/2)\normt{z_{t+1}-u_{t+1}}^2 + \zeta(z_{t+1}) + \mu\omega_\calU(z_{t+1})\big\} - E_t\\
&\ged A_{t+1} \big\{ h(z_{t+1}) - \delta_{t+1} + \zeta(z_{t+1}) + \mu\omega_\calU(z_{t+1})\big\} - E_t\\
&= A_{t+1} P(z_{t+1}) - E_{t+1},\quad \forall\,u\in\calU,\nt \label{eq:psi_t+1_lb2}
\end{align*}
where in (a) we use $\tau_{t+1} = \alpha_{t+1}/A_{t+1}$ and the 1-strong convexity of $\omega_\calU$ on $\calU$, in (b) we use the definition of $z_{t+1}$ in~\eqref{eq:z}, in (c) we use $z_{t+1} - u_{t+1} = \tau_{t+1}(w_{t+1} - \baru_{t+1})$, which follows from~\eqref{eq:u} and~\eqref{eq:z}, and 
$$\frac{A_t\mu+\barL}{A_{t+1}\tau_{t+1}^2} = \frac{(A_t\mu+\barL)A_{t+1}}{\alpha_{t+1}^2}\ge \barL,$$
which follows from~\eqref{eq:alpha_A_cond}, and in (d) we use the second inequality in~\eqref{eq:delta_L}.
We finish the induction by minimizing $\psi_{t+1}$  over $\calU$ on the left-hand side of~\eqref{eq:psi_t+1_lb2}. \Halmos
\endproof

At this point, a natural question would be whether the choices of $\{\alpha_t\}_{t\ge 0}$ in~Theorems~\ref{thm:conv_iAPG_mu=0} and~\ref{thm:conv_iAPG_mu>0} 
satisfy the condition~\eqref{eq:alpha_A_cond}. This is confirmed in the next lemma. 

\begin{lemma}\label{lem:step-size}
The choices of $\{\alpha_t\}_{t\ge 0}$ in~\eqref{eq:param_mu=0} when $\mu=0$ and in~\eqref{eq:param_mu>0} when $\mu>0$ lead to the following values of $\{A_t\}_{t\ge 0}$:
\begin{equation}
A_t = \begin{cases}
(t+2)^2/4, \quad & \mbox{when}\;\;\; \mu = 0,\\[1ex]
 (1+\sqrt{\theta})^t, \quad & \mbox{when}\;\;\; \mu > 0,
\end{cases} \qquad \forall\,t\ge 0. \label{eq:val_A_t}
\end{equation}
In addition, the condition in~\eqref{eq:alpha_A_cond} is satisfied in both cases (i.e., $\mu=0$ and $\mu>0$). 
\end{lemma}

\proof{Proof.}
Let us first focus on the case where $\mu=0$, and show  $A_t = (t+2)^2/4$ for all $t\ge 0$ via induction. 
Clearly, this holds when $t=0$ as $A_0 = \alpha_0 = 1$. Suppose this is true for some $t\ge 0$. Then 
\begin{equation}
A_{t+1} = A_t + \alpha_{t+1} = (t+2)^2/4 + (2t+5)/4 = (t+3)^2/4. 
\end{equation} 
This completes the induction. Note that when $\mu=0$, the condition in~\eqref{eq:alpha_A_cond} simplifies to $\alpha_t^2\le A_t$ for all $t\ge 1$, which clearly holds as $\alpha_t = (t+1.5)/2$ and $A_t = (t+2)^2/4$ for all $t\ge 1$. Next, let us show $A_t = (1+\sqrt{\theta})^t$ for all $t\ge 0$ when $\mu>0$.  Again, we prove this using induction.  Note that this clearly holds when $t=0$ as $A_0 = \alpha_0 = 1$. Suppose this is true for some $t\ge 0$. Then we have 
\begin{equation}
A_{t+1} = A_t + \alpha_{t+1} = (1+\sqrt{\theta})^t + (1+\sqrt{\theta})^t\sqrt{\theta} = (1+\sqrt{\theta})^{t+1}, 
\end{equation}
and this finishes the induction. 
Using~\eqref{eq:val_A_t} and the monotonicity of $\{A_t\}_{t\ge 0}$, we have 
\begin{equation}
(A_t\mu+\barL)A_{t+1} =\barL(1+\theta A_t)A_{t+1}\ge \barL \theta A_tA_{t+1} \ge  \barL\theta A_t^2 = \barL\theta(1+\sqrt{\theta})^{2t} = \barL \alpha_{t+1}^2, \quad \forall\,t\ge 0, \nn
\end{equation}
which means that the condition in~\eqref{eq:alpha_A_cond} is satisfied. We hence complete the proof. \Halmos
\endproof

Lastly,  the third lemma below establishes two upper bounds of $\psi^*_t$, for all $t\ge 0$. The first one involves the optimal value $P^*$ and $D_{\omega_\calU}(u^*,u_0)$, which is the Bregman ``distance'' from any optimal solution $u^*$ to $u^0$. The second one involves the dual function $\Xi$ and $\Omega_{\omega_\calU}(u^0)$ (cf.~\eqref{eq:def_Omega_U}), which is the Bregman ``distance'' from  the furthest point $u\in\calU$ to $u_0$, and is finite when $\calU$ is bounded. 

\begin{lemma} \label{lem:psi_t^*_ub}
For any $t\ge 0$, we have
\begin{equation}
\psi_t^*\le A_t P^*+\barL D_{\omega_\calU}(u^*,u_0). \label{eq:psi_t_P*}
\end{equation}
In addition, if $\calU$ is bounded and $(\hath(u_t),\hat{\nabla} h(u_t)) = (\barPsi(u_t,v_t), \nabla_u\barPsi(u_t,v_t))$ for some $\{v_t\}_{t\ge 0}\subseteq \calV$ (which need not satisfy~\eqref{eq:v_t_cond}), then for any $t\ge 0$, we have
\begin{equation}
\psi_t^*\le A_t \Xi(\barv_t)+\barL \Omega_{\omega_\calU}(u_0),
\label{eq:psi_t_dual}
\end{equation}
where $\barv_t$ is defined in~\eqref{eq:barv} and $\Omega_{\omega_\calU}(u_0)<+\infty$ is defined in~\eqref{eq:def_Omega_U}. 
\end{lemma}

\proof{Proof.}
Since $(\hath(u_t),\hat{\nabla} h(u_t))$ is a $(\delta_t,\barL)$-FOA of $h$ at $u_t$ for all $t\ge 0$, by using the first inequality in~\eqref{eq:delta_L}, we have that for all $u\in\calU$, 
\begin{align}
\psi_t^*\le \psi_t(u)\le \textstyle\sum_{i=0}^t \;\alpha_i(h(u) + \zeta(u)+ \mu\omega_\calU(u)) + \barL D_{\omega_\calU}(u,u_0) = A_t P(u) + \barL D_{\omega_\calU}(u,u_0), \label{eq:psi_t_ub} 
\end{align}
Substitute $u=u^*$ into~\eqref{eq:psi_t_ub} and we obtain~\eqref{eq:psi_t_P*}. 
If $(\hath(u_t),\hat{\nabla} h(u_t)) = (\barPsi(u_t,v_t), \nabla_u\barPsi(u_t,v_t))$ for all $t\ge 0$, then using the convexity of $\barPsi(\cdot,v_t)$, we have
\begin{equation}
\hath(u_t) + \ipt{\hat{\nabla} h(u_t)}{u-u_t} = \barPsi(u_t,v_t) + \ipt{\nabla_u\barPsi(u_t,v_t)}{u-u_t}\le \barPsi(u,v_t), \quad \forall\,t\ge 0. 
\end{equation}
Based on this, and by using the concavity of $\barPsi(u,\cdot)$, we have
\begin{align*}
\psi_t(u)&\le \textstyle\sum_{i=0}^t \;\alpha_i\barPsi(u,v_i) + A_t(\zeta(u)+ \mu\omega_\calU(u)) + \barL D_{\omega_\calU}(u,u_0)\\
&\le A_t(\barPsi(u,\barv_t) + \zeta(u)+ \mu\omega_\calU(u)) + \barL\Omega_{\omega_\calU}(u_0).  \nt\label{eq:psi_t_ub_dual}
\end{align*}
Now, by taking infimum of $u$ over $\calU$ on both sides of~\eqref{eq:psi_t_ub_dual}, we obtain~\eqref{eq:psi_t_dual}. \Halmos
\endproof

Based on the three lemmas above, the proofs of Theorems~\ref{thm:conv_iAPG_mu=0} and~\ref{thm:conv_iAPG_mu>0} are immediate. Indeed, by combining~\eqref{eq:psi_t*_lb} and~\eqref{eq:psi_t_ub}, we obtain the convergence rate of the primal optimality gap $P(z_t) - P^*$, namely
\begin{equation}
P(z_t) - P^*\le \frac{\barL D_{\omega_\calU}(u^*,u_0)}{A_t} + \frac{\sum_{i=0}^t A_i \delta_i}{A_t}. \label{eq:primal_gap_bound}
\end{equation}
In addition, if we let the sequence $\{v_t\}_{t\ge 0}$ in Lemma~\ref{lem:psi_t^*_ub} satisfy~\eqref{eq:v_t_cond}, then from Lemma~\ref{lem:1st_order_approx}, we know that $(\hath(u_t),\hat{\nabla} h(u_t)) = (\barPsi(u_t,v_t), \nabla_u\barPsi(u_t,v_t))$ is a $(\delta_t,\barL)$-FOA of $h$ at $u_t$ with $\barL = 2L_h$. Therefore, we can combine~\eqref{eq:psi_t*_lb} and~\eqref{eq:psi_t_dual} to obtain the convergence rate of the duality gap $\barDelta(z_t,\barv_t) = P(z_t) - \Xi(\barv_t)$ (cf.~\eqref{eq:barDelta}), 
namely
\begin{equation}
\barDelta(z_t,\barv_t)\le \frac{2L_h \Omega_{\omega_\calU}(u_0)}{A_t} + \frac{\sum_{i=0}^t A_i \delta_i}{A_t}.\label{eq:duality_gap_bound}
\end{equation}
Now, to show Theorem~\ref{thm:conv_iAPG_mu=0}, we can simply substitute the values of $\{A_t\}_{t\ge 0}$ when $\mu=0$ (cf.~\eqref{eq:val_A_t}) into~\eqref{eq:primal_gap_bound} and~\eqref{eq:duality_gap_bound} above. We can also show Theorem~\ref{thm:conv_iAPG_mu>0} in the same way, except that the values of $\{A_t\}_{t\ge 0}$ when $\mu>0$ (cf.~\eqref{eq:val_A_t}) are substituted. 

\section{Proof of Lemma~\ref{lem:grad_map}.} \label{app:proof_grad_map}
The first-order optimality condition of~\eqref{eq:u_+} yields  
\begin{align*}
&0\in \partial( \zeta +\iota_\calU)(u^+) + \nabla h(\baru)+e(\baru) + \mu\nabla \omega_\calU(u^+)- \barG \\
\Longleftrightarrow \quad & \barG - \nabla h(\baru) - e(\baru)\in \partial( \zeta +\iota_\calU)(u^+)+ \mu\nabla \omega_\calU(u^+)\\
\Longleftrightarrow \quad & [\xi\defeq \barG + \nabla h(u^+) - \nabla h(\baru) - e(\baru)]\in \nabla h(u^+) + \partial( \zeta +\iota_\calU)(u^+)+ \mu\nabla \omega_\calU(u^+) = \partial P(u^+), 
\end{align*}
By the $\mu$-strong convexity of $P$ on $\calU$, we have
\begin{align}
P^* = P(u^*) \ge P(u^+) + \lranglet{\xi}{u^*-u^+} + (\mu/2)\normt{u^*-u^+}^2 \ge P(u^+) - \normt{\xi}_*^2/(2\mu). \label{eq:grad_map_bound1}
\end{align}
On the other hand, we have $\normt{\nabla h(u^+) - \nabla h(\baru)}_*\le L_h \normt{u^+-\baru} = \normt{G}$, and hence 
\begin{align}
\normt{\xi}_*^2 &\le 3\big(\normt{\barG}_*^2 + \normt{\nabla h(u^+) - \nabla h(\baru)}_*^2 + \normt{e(\baru)}_*^2\big)\le 3\big(\normt{\barG}_*^2 + \normt{G}^2+ \normt{e(\baru)}_*^2\big). \label{eq:grad_map_bound2}
\end{align}
By combining~\eqref{eq:grad_map_bound1} and~\eqref{eq:grad_map_bound2}, we arrive at~\eqref{eq:ub_u+}.

\section{Proof of Theorem~\ref{thm:stoc_conv}.} \label{app:proof_stoc_conv}
Using the same arguments that lead to~\eqref{eq:bound_Q_lambda_q_lambda},  we have that for any $k\ge 1$, 
\begin{align}
 \bbE[Q^\lambda(x_{k+1};x_k)\,|\,x_k]&\le q^\lambda(x_k) +2\rho R_\calY(\omega_\calY)  + \eta= Q^\lambda(\prox(q,x_k,\lambda);x_k)  + 3\eta/2, \label{eq:proof_stoc1}
\end{align}
where the equality follows from~\eqref{eq:q_Q_prox} and $\rho = \eta/(4R_\calY(\omega_\calY))$ (cf.\ Algorithm~\ref{algo:PD_smoothing_stoc}).  
By the $(2\lambda)^{-1}$-strong-convexity of $Q^\lambda(\cdot;x_k)$ on $\calX$, we have 
\begin{align}
  Q^\lambda(x_k;x_k) - Q^\lambda(\prox(q,x_k,\lambda);x_k) \ge (4\lambda)^{-1} \norm{\prox(q,x_k,\lambda)-x_k}^2. \label{eq:proof_stoc2}
\end{align}
Combining~\eqref{eq:proof_stoc1} and~\eqref{eq:proof_stoc2}, we have 
\begin{equation}
\bbE[Q^\lambda(x_{k+1};x_k)\,|\,x_k]  + (4\lambda)^{-1} \norm{\prox(q,x_k,\lambda)-x_k}^2\le Q^\lambda(x_k;x_k) + 3\eta/2. 
\end{equation}
Since $Q^\lambda(x_{k+1};x_k) = q(x_{k+1}) + \lambda^{-1}D_{\omega_\calX}(x_{k+1},x_k) \ge q(x_{k+1})$ 
and $Q^\lambda(x_k;x_k) = q(x_k)$, we have
\begin{align}
\bbE[q(x_{k+1})\,|\,x_k] + (4\lambda)^{-1} \norm{\prox(q,x_k,\lambda)-x_k}^2\le q(x_k) + 3\eta/2. \label{eq:proof_stoc3}
\end{align}
If we telescope~\eqref{eq:proof_stoc3} over $k = 1,\ldots,K$, then we have
\begin{align}
\bbE[q(x_{K+1})] + (4\lambda)^{-1} \textstyle{\sum_{k=1}^K}\bbE[\norm{\prox(q,x_k,\lambda)-x_k}^2]\le q(x_{1}) + 3K\eta/2. 
\end{align}
Using the definition of $x_{\rm out}$ in Algorithm~\ref{algo:PD_smoothing_stoc} and the fact that $\bbE[q(x_{K+1})]\ge q^*$, we have
\begin{align}
\begin{split}
\bbE\big[\norm{\prox(q,x_{\rm out},\lambda)-x_{\rm out}}^2\big]&= (1/K)\textstyle{\sum_{k=1}^K}\bbE\big[\norm{\prox(q,x_k,\lambda)-x_k}^2\big]\\
&\le 4\lambda(q(x_{1}) - q^*)/K + 6\lambda\eta. 
\end{split}\label{eq:proof_stoc4}
\end{align}
Taking square root on both sides of~\eqref{eq:proof_stoc4} and using the choices of $\eta$ and $K$ in~\eqref{eq:eta_K}, 
we have 
\begin{align}
\textstyle\sqrt{\bbE\big[\norm{\prox(q,x_{\rm out},\lambda)-x_{\rm out}}^2\big]}\le \sqrt{4\lambda(q(x_{1}) - q^*)/K + 6\lambda\eta}\le \varepsilon\lambda/\beta_\calX. \label{eq:proof_stoc5}
\end{align}
Since the function $a\mapsto\sqrt{a}$ is concave on $\bbR_+$, we have 
\begin{equation}
\bbE\big[\norm{\prox(q,x_{\rm out},\lambda)-x_{\rm out}}\big]\le \sqrt{\bbE\big[\norm{\prox(q,x_{\rm out},\lambda)-x_{\rm out}}^2\big]}. \label{eq:proof_stoc6}
\end{equation}
Combining~\eqref{eq:proof_stoc5} and~\eqref{eq:proof_stoc6}, we complete the proof.

\end{APPENDICES}


\bibliographystyle{informs2014} 
\bibliography{math_opt,mach_learn}

\end{document}